\documentclass[reqno]{amsart}
\usepackage{amsmath,amsthm,amssymb,a4wide}
\usepackage{times}
\usepackage{enumerate}
\usepackage{mathrsfs,epsf,a4,color,graphicx,eucal,cite}

\theoremstyle{definition}
\newtheorem{theorem}{Theorem}[section]

\newtheorem{lemma}[theorem]{Lemma}
\newtheorem{definition}[theorem]{Definition}
\newtheorem{remark}[theorem]{Remark}

\newtheorem{problem}[theorem]{Problem}

\numberwithin{equation}{section}

\definecolor{ddmagenta}{rgb}{0.7,0,0.9}
\definecolor{ddcyan}{rgb}{0,0.2,1.0}
\definecolor{dred}{rgb}{.8,0,0}

\newcommand{\ITEM}[2]{\parbox[t]{.05\textwidth}{#1}\hfill\parbox[t]{.95\textwidth}{#2}\vspace*{.8mm}}

\newcommand{\DDD}[3]{\begin{array}[t]{c}#1\vspace*{-1em}\\_{#2}\vspace*{-.5em}\\_{#3}\end{array}}
\newcommand{\ddd}[3]{\DDD{\begin{array}[t]{c}\underbrace{#1}\vspace*{.6em}\end{array}}{\text{\footnotesize #2}}{\text{\footnotesize #3}}}

\newcommand\eps{}
\newcommand\R{\mathbb R}
\newcommand\N{\mathbb N}
\newcommand\Q{\mathcal Q}

\newcommand\DT[1]{\mathchoice
                 {{\buildrel{\hspace*{.1em}\text{\LARGE.}}\over{#1}}}
                 {{\buildrel{\hspace*{.1em}\text{\Large.}}\over{#1}}}
                 {{\buildrel{\hspace*{.1em}\text{\large.}}\over{#1}}}
                 {{\buildrel{\hspace*{.1em}\text{\large.}}\over{#1}}}}
\newcommand\DDT[1]{\mathchoice
   {{\buildrel{\hspace*{.1em}\text{\LARGE.\hspace*{-.1em}.}}\over{#1}}}
   {{\buildrel{\hspace*{.1em}\text{\Large.\hspace*{-.1em}.}}\over{#1}}}
   {{\buildrel{\hspace*{.1em}\text{\large.\hspace*{-.1em}.}}\over{#1}}}
   {{\buildrel{\hspace*{.1em}\text{\large.\hspace*{-.1em}.}}\over{#1}}}}

\newcommand\JUMP[2]{\mathchoice
                   {\big[\hspace*{-.3em}\big[#1\big]\hspace*{-.3em}\big]_{#2}}
                   {[\hspace*{-.15em}[#1]\hspace*{-.15em}]_{#2}}
                   {[\![#1]\!]_{#2}}
                   {[\![#1]\!]_{#2}}}
\newcommand\bbA{\mathbb A}
\newcommand\bbC{\mathbb C}
\newcommand\bbD{\mathbb D}
\newcommand\bbE{\mathbb E}
\newcommand\bbK{\mathbb K}
\newcommand\bbB{\mathbb B}
\newcommand\bbH{\mathbb H}

\newcommand\bbG{\mathbb G}

\renewcommand\d{\mathrm d}
\newcommand\w{\vartheta}
\newcommand\ent{h}

\newcommand\GC{\Gamma_{\mbox{\tiny\rm C}}}
\newcommand\SC{\Sigma_{\mbox{\tiny\rm C}}}
\newcommand\overlineGC{\overline{\Gamma}_{\mbox{\tiny\rm C}}}
\newcommand\overlineSC{\overline{\Sigma}_{\mbox{\tiny\rm C}}}
\newcommand{\dt}{\mathrm{D}_t}

\newcommand{\FRM}{F}
\newcommand{\GRM}{G}
\newcommand{\fRM}{f}
\newcommand{\gRM}{g}

\newcommand{\calD}{\mathcal{R}}

\newcommand{\aein}{\text{a.e. in}}
\newcommand{\foraa}{\text{for a.a.}}

\newcommand{\testu}{v}
\newcommand{\testw}{w}
\newcommand{\minus}{-}
\newcommand{\Gdir}{\Gamma_{\mbox{\tiny\rm D}}}
\newcommand{\Gnew}{\Gamma_{\mbox{\tiny\rm N}}}
\newcommand{\Sdir}{\Sigma_{\mbox{\tiny\rm D}}}
\newcommand{\Snew}{\Sigma_{\mbox{\tiny\rm N}}}

\newcommand{\dela}{\bbA}
\newcommand{\delam}{\bbA}
\newcommand{\het}{\eta}
\newcommand{\ind}{I}
\newcommand{\dd}{\mathrm{d}}
\newcommand{\piecewiseConstant}[2]{\overline{#1}_{\kern-1pt#2}}
\newcommand{\pwc}{\piecewiseConstant}
\newcommand{\underpiecewiseConstant}[2]{\underline{#1}_{\kern-1pt#2}}
\newcommand{\upwc}{\underpiecewiseConstant}

\newcommand{\piecewiseLinear}[2]{#1_{\kern-1pt#2}}
\newcommand{\pwl}{\piecewiseLinear}

\newcommand{\weaksto}{\stackrel{*}{\rightharpoonup}}
\newcommand{\weakto}{{\rightharpoonup}\,}
\newcommand{\GE}{\succeq}
\newcommand{\GEstar}{\succeq\hspace{-.8em}^{^*}\hspace{.3em}}

\newcommand{\pairing}[4]{ \sideset{_{#1 }}{_{ #2}}  {\mathop{\langle #3 , #4  \rangle}}}

\newcommand{\wt}{\w_{\eps\tau}}
\newcommand{\zt}{z_{\eps\tau}}
\newcommand{\ue}{u_{\eps}}
\newcommand{\we}{\w_{\eps}}
\newcommand{\ze}{z_{\eps}}

\newcommand{\ude}{u}
\newcommand{\wde}{\w}
\newcommand{\zde}{z}

\newcommand{\uk}{u_{\tau}^k}
\newcommand{\wk}{\w_{\tau}^k}
\newcommand{\zk}{z_{\tau}^k}
\newcommand{\hk}{\mathfrak{h}_{\tau}^k}

\newcommand{\ukm}{u_{\tau}^{k-1}}
\newcommand{\wkm}{\w_{\tau}^{k-1}}
\newcommand{\zkm}{z_{\tau}^{k-1}}

\newcommand{\cone}{K}

\newcommand{\Omegaone}{\Omega_+}
\newcommand{\Omegatwo}{\Omega_-}
\newcommand{\indabs}{\partial \abs}
\newcommand{\abs}{\mathcal{I}_{\cone}}
\newcommand{\norm}{\nu_{\mbox{\tiny\rm C}}}

\newcommand{\calZ}{\mathcal{Z}}

\newcommand{\calR}{\mathcal{R}}
\newcommand\Meas{\mathcal{M}}

\newcommand\Colon{{:}}

\def\neweta{n}  % <- TEMPORARY DEFINITION
\def\Nn{N_\neweta}      % <- TEMPORARY DEFINITION

\newenvironment{newtom}{\color{dred}}{\color{black}}
\newcommand\BT{\begin{newtom}}
\newcommand\ET{\end{newtom}}

\begin{document}

%%%%% To ease editing, add:

\title{Adhesive contact delaminating at mixed mode,
its thermodynamics and analysis}

\author{Riccarda Rossi}
\address{Dipartimento di Matematica, Universit\`a di
  Brescia, via Valotti 9, I--25133 Brescia, Italy.}
\email{\ttfamily riccarda.rossi\,@\,ing.unibs.it}

\author{Tom\'a\v s Roub\'\i\v cek}

\address{Mathematical Institute, Charles University,
Sokolovsk\'a 83, CZ--186~75~Praha 8,
and
Institute of Thermomechanics of the ASCR,
Dolej\v skova~5, CZ--182 00~Praha 8, Czech Republic.}
\email{\ttfamily tomas.roubicek\,@,mff.cuni.cz}
\thanks{This research was initiated during a visit of R.R. at the
%Mathematical Institute of
Charles University in Prague, whose hospitality is gratefully
acknowledged, with the support of   the ``Ne\v cas center for
mathematical modeling'' LC 06052 (M\v SMT \v CR), and partially of
the PRIN 2008 project ``Optimal mass transportation, geometric and
functional inequalities and applications". T.R.~acknowledges the
hospitality of the University of Brescia supported from Inst.\
Nazionale di Alta Matematica F.Severi (GNAMPA project), as well as
partial support from the grants A~100750802 (GA~AV~\v CR),
%106/09/1573
201/09/0917 and 201/10/0357 (GA \v CR), and MSM~21620839 (M\v SMT \v
CR), and from the research plan AV0Z20760514 ``Complex dynamical
systems in thermodynamics, mechanics of fluids and solids'' (\v CR).}

\begin{abstract}
An adhesive unilateral contact between visco-elastic heat-conductive
bodies in linear Kelvin-Voigt rheology is scrutinised. The flow-rule
for debonding the adhesive is considered rate independent,
unidirectional, and non-associative due to dependence on the mixity
of modes of delamination, namely of Mode I (opening) and of Mode II
(shearing). Such mode-mixity  dependence of  delamination is a very
pronounced (and experimentally confirmed) phenomenon typically
considered in engineering models.  An anisothermal,
thermodynamically consistent model is derived, considering a
heat-conductive viscoelastic material and the  coupling via thermal
expansion and adhesion-depending heat transition through the contact
surface.  We prove the  existence of weak solutions by passing to
the limit in a carefully designed semi-implicit time-discretization
scheme.

 \noindent {\bf AMS Subject Classification}:  35K85;  74A15.

  \noindent {\bf Key words:} Contact mechanics;
rate-independent processes;
%Adhesive contact; nonlinear
heat equation;
%rate-independence; non-associative model;
energetic solution.
%; existence.
\end{abstract}
%%%%%%%%%%%%%%%%%%%

  %%%%%%%%%%%%%%%%%%%%%%%%%%%%%%%%%%%%%%%%%%%%%%%%%
  %%%%%%%%%%%%%%%%%%%%%%%%%%%%%%%%%%%%%%%%%%%%%%%%%%%%%%
\maketitle

\section{Introduction}\label{sec-intro}

%\begin{rcomm} to modify and possibly shorten \end{rcomm}{\tiny

Nonlinear contact mechanics is an important part of  mechanical
engineering and receives still growing attention due to its numerous
applications. We focus on the  modelling and analysis of an
inelastic process called delamination (sometimes also  referred to
as debonding), between two elastic bodies, glued together along a
prescribed delamination surface.
 ``Microscopically'' speaking, some
macromolecules in the adhesive may break upon loading and we assume
that they can never be glued back, i.e., no ``healing'' is possible.
This makes the process {\it unidirectional},  viz.\ irreversible. On
the glued surface, we consider the {\it delamination process} as
{\it rate-independent} and, in the bulk, we also take into account
rate-dependent {\it inertial, viscous}, and {\it heat-transfer
effects}. The ultimate phenomenon counted in engineering modelling
(and so far mostly ignored in the mathematical literature), is  the
dependence of this process  on the  modes under which it proceeds.
Indeed, Mode I (=opening) usually dissipates much less energy than
Mode II (=shearing). The difference may be tens or even hundreds of
percents, cf.\
\cite{BanAsk00NFCI,LieCha92ASIF,Mant08DRLM,SwLiLo99ITAM}. Moreover,
the delamination process rarely  follows
 such pure modes:   in general, the
%so-called
{\it mixed mode} is favoured. Microscopically, this difference is
explained either by some roughness of the glued surface (to be
overcome in Mode II but not in Mode I, cf.\ \cite{ERDC90FEBI}) or by
some plastification, either in the adhesive or in a narrow strip
around the delaminating surface, before the delamination itself is
triggered (and, as usual in plasticity involving trace-free plastic
strain, again only shearing in Mode II  can trigger this
plastification, not opening in Mode I,
cf.~\cite{LieCha92ASIF,TveHut93IPMM}).

In this article, we focus  on a standard engineering model which,
however, was never rigorously analysed so far. Even for the
isothermal variant of this model, existence of solutions  has not
yet been proved, although computational simulations are  routinely
launched and successfully used in applications.  In particular,
%in this paper
we  extend the analysis of \cite{KoMiRo,rr+tr} to the
mixity-sensitive case. For other results on  models for
rate-independent adhesive contact, we refer to
\cite{MiRoTh??DDNE,tr-LS-CZ}. The analysis of models featuring a
\emph{rate-dependent} evolution for the delamination variable was
carried out in \cite{Point,Raous,BBR1,BBR2,BBR3,BBR4,BBR5}, cf.\ also the
monographs~\cite{SoHaSh06AACP,eck} for further references.

The initial-boundary-value problem we are going to analyse is
written down in its classical formulation in Section \ref{sec2}. In
the following lines, we just highlight the main features of the
model, in particular focusing on the mixity of delamination modes.
We confine ourselves to {\it small strains} and, just for the sake
of notational simplicity, we restrict the analysis to the case of
two (instead of several) bodies $\Omegaone$ and $\Omegatwo$ glued
together along the {\it contact surface} $\GC$. The material in the
bulk is taken to be heat conductive, and thus the system is
completed by the nonlinear heat equation in a thermodynamically
consistent way. The contact surface is considered infinitesimally
thin, so that the thermal capacity of the adhesive is naturally
neglected. The coupling of the mechanical and thermal effects thus
results from thermal expansion, dissipative/adiabatic heat
production/consumption (depending, in particular, on the mode mixity
on $\GC$), and here also from the
%\CHECK{????  I don't understand what the following means????}
possible dependence of the heat transfer through the contact surface
$\GC$ on the
%stage of
delamination itself, and on the possible slot between the bodies if
the contact is debonded.
%\CHECK{$<\!\!--$IF WE GO TO THIS GENERALIZATION}

We consider an elastic response of the adhesive, and then one speaks
about  {\it adhesive contact} (in contrast to  brittle contact where
a mixity-dependent model seems to be particularly difficult to
analyse). The elastic response in the adhesive will be assumed
linear, determined by the (positive-definite) matrix of elastic
moduli $\bbA$. At a current time, the ``volume fraction'' of
debonded molecular links will be ``macroscopically'' described by
the {\it scalar delamination parameter} $z:\GC\to[0,1]$, which  can
be referred to the modelling approach by
%\textsc
{M.\,Fr\'emond}, see \cite{Fre82,Fre87}.
 The state $z(x)=1$ means that the adhesive is
still $100\%$ undestroyed and thus fully effective, while the
intermediate state $0<z(x)<1$ means that there are some molecular
links which have been broken but the remaining ones are effective,
and eventually $z(x)=0$ means that the surface is already completely
debonded at $x\in\GC$.
%As already pointed out in~\cite{KoMiRo},
In some simplification, based on the {\it Griffith concept}
\cite{Grif20PRFS}, it is assumed that a specific phenomenologically
prescribed energy $a$ (in J/m$^2$, in 3-dimensional situations) is
needed to break the macromolecular structure of the adhesive,
independently of the rate of this process. Thus, delamination is a
\emph{rate-independent} and activated phenomenon, governed by the
maximum dissipation principle, and we shall accordingly consider a
rate-independent flow rule for $z$.
%Activating the delamination process in the adhesive
%contact at a given point $x\in\GC$ again needs the
%(phenomenologically prescribed) energy $a(x)$.

Let us now comment on the main new feature of the model presently
analysed, i.e.\ its mixity-sensitivity.  An immediate reflection of
the standard engineering approach as e.g.~in
\cite{HutSuo92MMCL,TMGCP10ACTA,TMGP??BEMA} is to make the activation
energy $a=a(\psi_{\rm G})$ depend on the so-called \emph{mode-mixity
angle} $\psi_{\rm G}$.  For instance, if $\norm=(0,0,1)$
 at some $x \in \GC$ (with $\norm$ the unit normal to
$\GC$, oriented from $\Omegaone$ to $\Omegatwo$), and $\bbA={\rm
diag}(\kappa_{\rm n},\kappa_{\rm t},\kappa_{\rm t})$, the
mode-mixity  angle
  is defined as
\begin{align}\label{mixity-angle}
\psi_{\rm G}=\psi_{\rm G}(\JUMP{u}{}):={\rm arc\,tan}
\sqrt{\frac{\kappa_{\rm t}|\JUMP{u}{\rm t}|^2}{\kappa_{\rm
n}|\JUMP{u}{\rm n}|^2}}
\end{align}
where $\JUMP{u}{\rm t}$ and $\JUMP{u}{\rm n}$ stand for
 the tangential and the normal traction. They give  the
decomposition of the jump of displacement across the boundary $\GC$
as $\JUMP{u}{}= \JUMP{u}{\rm n}\norm+\JUMP{u}{\rm t}$, with
$\JUMP{u}{\rm n}=\JUMP{u}{}\cdot\norm$.
% as in \eqref{E-delam-small0}.
In fact, to satisfy natural analytical
 requirements (viz.\ the continuity of the function $a_1$ below, cf.\
\eqref{a-mixity-dependence+}), one should rather consider a suitable
regularization of \eqref{mixity-angle} to avoid discontinuity at
$0$.
 For example, it is sufficient to take
\[
\psi_{\rm G}(\JUMP{u}{})={\rm arc\,tan} \sqrt{\frac{\kappa_{\rm
t}|\JUMP{u}{\rm t}|^2} {\kappa_{\rm n}|\JUMP{u}{\rm n}|^2+\epsilon}}
\quad \text{with a small $\epsilon>0$.}
\]
The coefficient  $\kappa_{\rm t}$ is often smaller than $\kappa_{\rm
n}$, and a typical phenomenological form of $a(\cdot)$ used in
engineering \cite{HutSuo92MMCL} is, e.g.,
\begin{align}\label{a-mixity-dependence}
a(\psi_{\rm G}):=a_{\rm I}\big(1+{\rm tan}^2((1{-}\lambda)\psi_{\rm
G})\big).
\end{align}
In \eqref{a-mixity-dependence},  $a_{\rm I}=a(0)$ is the activation
threshold for
%fracture
delamination mode I and $\lambda$ is the so-called
\emph{delamination-mode-sensitivity} parameter.  Note that  a
moderately strong
%fracture
delamination-mode sensitivity occurs when the ratio $a_{\rm
II}/a_{\rm I}$ is about 5-10 (where $a_{\rm II}=a(90^\circ)$ is the
activation threshold for the pure
%fracture
delamination mode II). Then,
 one has $\lambda$ about 0.2-0.3; cf.~\cite{TMGP??BEMA}.

In the thermodynamical context, the energy $a$ needed for
delamination is dissipated by the system in two ways: one part $a_1$
is spent to the chaotic vibration of the atomic lattice of both
sides of the delaminating surface $\GC$, which leads
``macroscopically'' to heat production (cf.~also \cite[Remark
4.2]{tr-LS-CZ}), while another part $a_0$ is spent to create a new
delaminated surface (or, ``microscopically'' speaking, to break the
macromolecules of the adhesive). Thus $a=a_0+a_1$. Consistently with
the dissipation, Mode II also heats up considerably more than Mode
I, as experimentally documented in  \cite{Ritt99TADC}. In view of
\eqref{a-mixity-dependence}, an option suggested already in
\cite{tr+mk+jz} is
\begin{align}\label{a-mixity-dependence+}
a_0\big(\JUMP{u}{}\big):=a_{\rm I},\qquad a_1\big(\JUMP{u}{}\big):=
a_{\rm I}{\rm tan}^2\big((1{-}\lambda)\psi_{\rm G}(\JUMP{u}{})\big),
\end{align}
meaning that  plain delamination does not  trigger heat production
at all, and only the additional dissipation related with Mode II
contributes to the heat production on the delaminating surface. We
summarize the features of these particular modes in
Table~\ref{table-of-modes}.

\begin{table}
{\small\sffamily
\begin{center}
\begin{tabular}{|c||c|c|c|}
\cline{1-4} \ \ \ \ \ \ \ \ features: & energy & heat & mode-mixity
\\[-.2em]
mode:\ \ \ \ \ \ \ \ \ \ & dissipation & generation & angle
$\psi_{\rm G}$
\\[-.0em]\cline{1-4}
& & &  \\[-.9em]
\cline{1-4}
Mode I $^{^{^{^{}}}}$ & small & small & 0$^\circ$ \\[-.3em]
(opening) &  &   or none  &
\\[-0mm]\cline{1-4}
mixed\ \ \ $^{^{^{^{}}}}$ & moderate & moderate & in between
\\[-.3em]
\ \ \ mode &  & & 0$^\circ$ and 90$^\circ$
\\[-0mm]\cline{1-4}
Mode II $^{^{^{^{}}}}$ & large & large & 90$^\circ$
\\[-.3em]
(shearing) & & &
\\[0mm]
\cline{1-4}
\end{tabular}\\[-5.5em]
\end{center}
} \vspace*{-1.5em} \caption{Schematic summary of particular
modes.}\label{table-of-modes}
\end{table}

%In the present anisothermal case, we will be able
%to prove such a result only in rather special cases. To this goal,
%it is natural to formulate also the brittle delamination problem in
%terms of the delamination parameter, requiring the jump of
%displacement on $\GC$ to be zero wherever $z>0$. In fact, it
%generalizes the Griffith concept which, in terms of $z$, implicitly
%relies on $z$ taking the values only 0 or 1.

The mathematical difficulties  attached to the analysis of the PDE
system for the present mixity-sensitive delamination model  arise
both from the proper thermodynamical coupling, and from hosting an
inelastic rate-independent process on $\GC$.
  Models combining  thermal and rate-independent effects
have already been successfully analysed in \cite{tr1} for inelastic
processes in the bulk, and in \cite{rr+tr} for surface delamination.
 The essential ingredient for the analysis  is the satisfaction of
the energy balance.  In this direction,  the concept of energetic
solutions to rate-independent systems recently developed in
\cite{Miel05ERIS,MieThe99MMRI,MieThe04RIHM,MiThLe02VFRI} and adapted
to systems with inertia and viscosity in \cite{tr2} appears truly
essential. Here, additional difficulties are related with the
mixity-dependence of the dissipation, which makes it {\it
non-associative}, in contrast to the mixity-insensitive case and to
another model recently devised and analysed in an isothermal case in
\cite{tr+mk+jz,tr+vm+cp}. This analytical feature has led us to
resort to a higher-order gradient in the momentum equation via the
concept of the so-called hyperstresses, already justified and used
in  the theoretical-mechanical literature,  cf.\ e.g.\
\cite{toupin,fried-gurtin,ppg-mv}.  Such a regularization brings
various inevitable technicalities into the classical formulation of
the problem, cf.\ \eqref{eq6:adhes-class-form} and \eqref{eqsystem}
later on.

The main result of this paper is the existence of solutions to the
initial-boundary value problem associated with the mixity-sensitive
model under investigation. The proof is performed by passing to the
limit in a suitably devised semi-implicit discretization scheme,
cf.\ \eqref{semi-implicit}. Let us mention that such a kind of
scheme (already applied in \cite[Sect.3.1]{LaOrSu10ESRM} for a
special dynamic isothermal fracture problem) leads to considerable
analytical simplifications, in comparison with the fully implicit
scheme used in \cite{rr+tr}.  In the existence proof we shall
distinguish the \emph{dynamic} case, involving inertial terms  in
the momentum equation, and the \emph{quasistatic} one, where inertia
is neglected. In the latter situation, we will be able to tackle a
fairly general contact conditions for the displacement variable $u$,
in particular including (a generalization of) the Signorini
frictionless contact law. For further explanations and comments, we
refer to Remark \ref{rmk:4.3}.

The plan of the paper is as follows:  in Section~\ref{sec2} we
formulate the initial-boundary-value problem in its classical
formulation.  We also briefly sketch its derivation, referring to
\cite{rr+tr} for more details. Hence, in Sec.\ \ref{sec-ent-trans},
after introducing a suitable transformation for the heat equation,
we advance a suitably devised weak formulation of our  PDE system,
and comment on its relation to the classical formulation of
Sec.~\ref{sec2}. We state our main existence result in
Sec.~\ref{sec-main-res}, and set up the approximation via a
  semi-implicit time-discretization scheme in Sec.~\ref{s:new}.
 For the discrete solutions, suitable a priori estimates are
 obtained, which allow us to pass to the limit
 in the time-discrete approximation,
 and conclude the existence of solutions in Sec.~\ref{ss:4.4}.

\section{The model and its derivation}\label{sec2}
%        ~~~~~~~~~~~~~~~~~~~~~~~~~~~
Hereafter, we  suppose that the elastic body occupies a  bounded
Lipschitz domain $\Omega\subset\R^3$.
%$$\text{$\, \Omega\subset\R^3\,$,  bounded and with a
%Lipschitz boundary $\partial\Omega$.}$$
We assume that
\[
\Omega = \Omegaone \cup \GC \cup \Omegatwo\,,
\]
 with $\Omegaone$ and $\Omegatwo$
disjoint Lipschitz subdomains and $\GC$ their common boundary, which
represents a prescribed
%delamination two-dimensional
surface where delamination may occur. We  denote by $\nu$ the
outward unit normal to $\partial \Omega$, and by $\norm$
 the unit normal to $\GC$, which we consider
oriented from $\Omegaone$ to $\Omegatwo$. Moreover, given  $v \in
W^{1,2} (\Omega {\setminus} \GC)$, $v^+$ (respectively, $v^-$)
signifies the restriction of $v$ to $\Omegaone$ (to $\Omegatwo$,
resp.).
 We further suppose  that the boundary of $\Omega$ splits as
\[
\partial \Omega = \Gdir\cup \Gnew\,,
\]
with $\Gdir$ and $\Gnew $ open subsets in the relative topology of
$\partial\Omega$,  disjoint one from each other and each of them
with a smooth (one-dimensional) boundary. Considering $T>0$ a fixed
time horizon, we  set
\begin{displaymath}
Q:=(0,T)  \times \Omega, \quad \Sigma : = (0,T)  \times \partial
\Omega, \quad \SC:= (0,T)  \times\GC , \quad \Sdir:= (0,T) \times
\Gdir , \quad \Snew:= (0,T) \times\Gnew.
\end{displaymath}
For readers' convenience, let us summarize the basic notation used
in what follows:

\vspace{.7em}

\!\!\!\begin{minipage}[t]{0.49\linewidth}
%\small

$u:\Omega{\setminus}\GC \to\R^3$ displacement,

$\theta:\Omega{\setminus}\GC\to(0,+\infty)$ absolute temperature,

$z:\GC\to[0,1]$ delamination variable,

$e=e(u)=\frac12\nabla u^\top\!+\frac12\nabla u$ small strain tensor,

$\JUMP{u}{}= u^+|_{\GC} - u^-|_{\GC}\ $ jump of $u$ across $\GC$,

$\sigma$  stress tensor,

$\w$ rescaled temperature (enthalpy),

$\bbC\in\R^{3^4}$ elasticity constants,

$\bbD\in\R^{3^4}$ viscosity constants,

$\bbH\in\R^{3^6}$ elasticity constants for hyperstress,

$\bbG\in\R^{3^6}$ viscosity constants for hyperstress,

$\bbK{=}\bbK(e,\theta){\in}\R^{3\times 3}$
%matrix of
heat-conduction coefficients,

$\bbE\in\R^{3\times 3}$
%matrix of
thermal-expansion coefficients,

$\bbB:=\bbC\bbE$,

\end{minipage}\ \
\begin{minipage}[t]{0.49\linewidth}\small

$\bbA\in\R^{3\times 3}$
%matrix of
elastic coefficients of the adhesive,

$\varrho>0$  mass density,

$ c_{\rm v}=c_{\rm v}(\theta)$  heat capacity,

$a_0=a_0(\JUMP{u}{})$
%specific
energy
%(per area) which is
stored
%by disintegrating the adhesive
on $\GC$,

$a_1=a_1(\JUMP{u}{})$
%specific
energy
%(per area) which is
dissipated
%by debonding
%disintegrating
%the adhesive,
on $\GC$,

%$\het{=}\het(x,\JUMP{u}{},z)$
$\het{=}\het(\JUMP{u}{},z)$ heat-transfer coefficient on $\GC$,

$\FRM:Q\to\R^3$  applied bulk force,

$w_{\rm D}
%:\Sdir\to\R^3
$  prescribed
%time-dependent
boundary displacement,

$\fRM:\Snew\to\R^3$  applied traction,

$\GRM:Q\to\R$  bulk heat source,

$\gRM:\Sigma\to\R$  external heat flux,

$\psi$'s (bulk and surface) free energies,

$\zeta$'s
%(bulk and surface)
(pseudo)potentials of dissipative forces,

$\xi$'s (bulk and surface) rates of dissipation.
%(or, equivalently, heat production rates).

\end{minipage}\medskip

\noindent The {\it state} is formed by the triple $(u,\theta,z)$.
 We use {\it
Kelvin-Voigt's rheology} and {\it thermal expansion}.
 As a further contribution to the \emph{stress}
$\sigma:(0,T)\times\Omega\rightarrow\R^{3\times 3}$, we also
consider the {\it hyperstress} $\mathfrak
h:(0,T)\times\Omega\rightarrow\R^{3\times 3\times 3}$, which
accounts for ``capillarity-like'' effects. Mimicking
 Kelvin-Voigt's rheology, we  incorporate in the hyperstress contribution to $\sigma$
the corresponding dissipation mechanisms.  Hence we assume the
\emph{stress} $\sigma:(0,T)\times\Omega\rightarrow\R^{3\times 3}$ in
the form:
\begin{align}
\sigma=\sigma(u,\DT{u},\theta) :=\!\!\!\ddd{\bbD
e(\DT{u})}{viscous}{stress}\!\!\!
+\!\!\!\ddd{\bbC\big(e(u){-}\bbE\theta\big)}{elastic}{stress}\!\!\!
-{\rm div}\big(\!\!\!\ddd{\bbH\nabla e(u)+ \bbG\nabla
e(\DT{u})}{$=:\mathfrak h$}{hyperstress}\!\!\!\big),
\end{align}
The $\bbG$-term will ensure the acceleration $\DDT u$ to be in
duality with the velocity $\DT u$, which is needed if $\varrho>0$ is
considered.
%\begin{rcomm}
%REMARK: it seems to me that we NEED the term $\bbG\nabla e(\DT{u})$ in the hype%rstress, NOT because of
%compactness problems, but because we need higher spatial regularity of $\DT{u}$% to make all the arguments RIGOROUS.
%\end{rcomm}
Furthermore, we shall denote by $T=T(u,v,\theta)$  the traction
stress on some two-dimensional surface $\Gamma$
 (later, we shall take either $\Gamma= \GC$ or
$\Gamma=\Gnew$), i.e.
\begin{align}\label{T-stress}
T(u,\DT{u},\theta):=\sigma(u,\DT{u},\theta)\big|_{\Gamma}\nu\,,
\end{align}
where of course we take as $\nu$ the unit normal  $\norm$ to $\GC$,
if $\Gamma=\GC$.

We address a generalization of the standard frictionless Signorini
conditions on $\GC$ for the displacement $u$.  This is rendered
through a closed, convex cone $\cone(x)\subset\R^3$, possibly
depending on $x\in\GC$. In terms of   the multivalued, cone-valued
mapping $\cone: \GC \rightrightarrows \R^3$,
 the boundary conditions
on $\GC$  can be given in the complementarity form as
\begin{align}\label{b.c.}
\left.\begin{array}{rll}
\JUMP{u}{}&\GE&0,\ \ \\[0em]
T(u,\DT{u},\theta)&\GEstar&0,\ \ \\[0em]
T(u,\DT{u},\theta)\cdot\JUMP{u}{}&=&0
\end{array}\right\}
\ \text{ on }\GC.
\end{align}
In \eqref{b.c.}, $\GE $ is the ordering induced by the mapping
$\cone: \GC \rightrightarrows \R^3$,  in the sense that, for
$v_1,v_2: \GC \to \R^3$,
\begin{equation}\label{e:ordering}
\text{$v_1\GE v_2\ $ if and only if $v_1(x){-}v_2(x)\in \cone(x)$
for a.a.~$x\in\GC$.}
\end{equation}
Likewise,  $\GEstar$ is the dual ordering induced by the negative
polar cone to $\cone$, viz.\ for $\zeta_1,\zeta_2: \GC \to \R^3$,
\[
\text{$\zeta_1 \GEstar\zeta_2\ $ if and only if
$\zeta_1(x){\cdot}v\ge\zeta_2(x){\cdot}v$ for all $v\in \cone(x)$,
for a.a.~$x\in\GC$.}
\]
Possible  choices for the cone-valued  mapping $\cone:\GC
\rightrightarrows \R^3$ are
\begin{subequations}\label{CC}\begin{align}
%\label{CC1}
%&&&&&&&&&\cone(x)=\cone=\R^3 \quad \foraa\, x \in \GC,&&\text{or}&&&&&&&&\\
\label{CC2} &&&&&&&&&\cone(x)=\{v\in\R^3;\ v\cdot\norm(x)\ge0\}\quad
\foraa\, x \in \GC,&&\text{or }&&&&&&&&
\\\label{CC3}
&&&&&&&&&\cone(x)=\{v\in\R^3;\ v\cdot\norm(x)=0\}\quad \foraa\, x
\in \GC.&&&&&&&&
\end{align}\end{subequations}
 In the case \eqref{CC2},
\eqref{b.c.} reduce to the standard unilateral frictionless
Signorini contact conditions for
 the normal displacement.  The
%last
case \eqref{CC3} prescribes  the normal jump of the displacement,
variable at $x\in\GC$, to zero. Thus, it only allows   for a
tangential slip along $\GC$. This may be a relevant model under high
pressure, when no cavity of $\GC$ can be expected anyhow. Such a
situation occurs, e.g., on lithospheric faults deep under the earth
surface. Note that, in \eqref{CC3}, $\cone(x)$ is a linear manifold
for a.a. $x \in \GC$.
 Indeed, later on we shall have to assume the latter
property in the case $\varrho>0$. %%%%%%%%%%%%%%%%%%
%%%%%%%%%%%%%%%%%%%%%%%%%
\paragraph{Classical formulation of the adhesive contact problem.}
 Beside the force equilibrium coupled with the heat equation inside
   $\Q{\setminus}\SC$ and supplemented  with standard boundary
  conditions, we have two complementarity problems on $\SC$. Altogether, we
  have the boundary-value problem
\begin{subequations}\label{eq6:adhes-class-form}
\begin{align}
\label{eq6:adhes-class-form1} & \varrho\DDT{u}
-\mathrm{div}\big(\bbD e(\DT{u})+\bbC
e(u)-\bbB\theta-\mathrm{div}\,\mathfrak h\big)= \FRM, \qquad
\mathfrak h=\bbH\nabla e(u){+}\bbG\nabla e(\DT{u}) &&\text{in
}Q{\setminus}\SC,
\\[-.5em]
\label{eq6:adhes-class-form1bis}
 & c_{\rm v}(\theta)\DT{\theta}
-\mathrm{div}\big(\bbK(e(u),\theta)\nabla\theta\big)= \bbD
e(\DT{u}){:}e(\DT{u}) -\theta\bbB{:}e(\DT{u})+ \bbG\nabla
e(\DT{u}){\vdots}\nabla e(\DT{u})
%+\GRM
&&\text{in }Q{\setminus}\SC,
\\
\label{eq6:adhes-class-form2} &u=0 &&\text{on }\Sdir,\hspace{1.2em}
\\\label{eq6:adhes-class-form3-bis}
&T(u,\DT{u},\theta)-\mathrm{div}_{{}_{\mathrm{S}}} (\mathfrak{h}
\cdot \nu)=\fRM &&\text{on }\Snew,\hspace{1.2em}
\\\label{eq6:adhes-class-form3}
&(\bbK(e(u),\theta)\nabla\theta)\nu=\gRM &&\text{on
}\Sigma,\hspace{1.7em}
\\
\label{eq6:adhes-class-form-hot-BC} & \mathfrak{h} \Colon (\nu
\otimes \nu)=0 &&\text{on }\Sdir\cup\Snew,\hspace{-.15em}
\\\label{adhes-form-d1}
& \JUMP{\bbD e(\DT{u})+ \bbC e(u)-\bbB\theta
-\mathrm{div}\,\mathfrak h}{}\norm - \mathrm{div}_{{}_{\mathrm{S}}}
(\JUMP{\mathfrak{h}}{}\norm) =0
 &&\text{on
}\SC,\hspace{1.2em}
\\
\label{eq6:new-jump-BC} & \mathfrak{h}^+ \Colon (\norm \otimes
\norm)=\mathfrak{h}^-\Colon(\norm\otimes\norm)=0 &&\text{on
}\SC,\hspace{-.15em}
\\
\label{adhes-form-d2} &\JUMP{u}{}\GE0  &&\text{on
}\SC,\hspace{1.2em}
\\\label{adhes-form-d3}
& T(u,\DT{u},\theta) -\mathrm{div}_{{}_{\mathrm{S}}}
(\mathfrak{h}\cdot\nu) +z\bbA\JUMP{u}{}{-}za_0'(\JUMP{u}{})\GEstar0
&&\text{on }\SC,\hspace{1.2em}
\\
\label{adhes-form-d4} &
\big(T(u,\DT{u},\theta){-}\mathrm{div}_{{}_{\mathrm{S}}}(\mathfrak{h}\cdot\nu)
{+}z\bbA\JUMP{u}{}{-}za_0'(\JUMP{u}{})\big)\cdot\JUMP{u}{}=0
&&\text{on }\SC,\hspace{1.2em}
\\\label{adhes-form-d6}
&\DT{z}\le0&&\text{on }\SC,\hspace{1.2em}
\\\label{adhes-form-d7}
& d\le a_0(\JUMP{u}{})+a_1(\JUMP{u}{}) &&\text{on
}\SC,\hspace{1.2em}
\\\label{adhes-form-d8-bis}
& \DT{z} \left( d - a_0(\JUMP{u}{})-a_1(\JUMP{u}{})\right) =0
&&\text{on }\SC,\hspace{1.2em}
\\\label{adhes-form-d8}
& d\in \partial I_{[0,1]}(z)+
%\partial\big[\IRM_{\cdot\JUMP{u}{}=0}\big](z)
%\partial_zJ\big(\JUMP{u}{},z\big)
\mbox{$\frac12$}\bbA\JUMP{u}{}{\cdot}\JUMP{u}{} &&\text{on
}\SC,\hspace{1.2em}
\\\label{adhes-form-d9}
& \mbox{$\frac12$}\big(\bbK(e(u),\theta)\nabla\theta|_{\GC}^+
+\bbK(e(u),\theta)\nabla\theta|_{\GC}^-\big)\cdot\norm
+\het(\JUMP{u}{},z)\JUMP{\theta}{}=0 &&\text{on }\SC,\hspace{1.2em}
\\
\label{adhes-form-d10} &
\JUMP{\bbK(e(u),\theta)\nabla\theta}{}\cdot\norm ={-}
a_1(\JUMP{u}{})\DT{z} &&\text{on }\SC,\hspace{1.2em}
\end{align}
\end{subequations}
where $\vdots$ denotes the tensorial product involving summation
over $3$ indices, and
 we have used the notation
$\mathfrak{h}^+=\mathfrak{h}|_{\Omegaone}$ and
$\mathfrak{h}^-=\mathfrak{h}|_{\Omegatwo}$. In
\eqref{eq6:adhes-class-form3-bis}, we denoted by
$\mathrm{div}_{{}_{\mathrm{S}}}$  the two-dimensional ``surface
divergence", defined by $\mathrm{div}_{{}_{\mathrm{S}}}:=
\mathrm{tr}(\nabla_{{}_{\mathrm{S}}}),$
%\begin{subequations}
%\label{surf-div}
%\begin{equation}
%\label{surface-div}
%\end{equation}
 where $\mathrm{tr}$ is the trace
operator (of a $2\times 2$ matrix), and $\nabla_{{}_{\mathrm{S}}}$
 denotes the tangential derivative, defined by
$ \nabla_{{}_{\mathrm{S}}} v = \nabla v -\frac{\partial v}{\partial
\nu} \nu.$
% \end{subequations}
%%%%%%%%%%%%%
\begin{subequations}\label{ass-on-K-C-D-H-G}
As to the involved tensorial symbols,
%In addition, we assume
\begin{align}
&\bbK=\bbK(e,\theta)\quad\text{ is a $2$nd-order positive definite
tensor,} \intertext{i.e.\ a $3{\times}3$-matrix, while}
\label{posit}&\bbC,\,\bbD:\R_{\mathrm{sym}}^{3\times 3}\to
\R_{\mathrm{sym}}^{3\times 3}\quad\text{ are $4$th-order positive
definite and symmetric tensors, $\bbC$ potential,} \intertext{in the
sense that both $\bbC e(u)$ and $\bbD e(u)$ are
$\R_{\mathrm{sym}}^{3\times 3}$-valued, and the operator div$\bbC
e(u)$ has a potential. Analogously, for the  higher-order terms we
suppose that
%we assume
%i.e. $\bbC_{ijkl}= \bbC_{jikl}= \bbC_{klij}$, and the same for
%$\bbD$,
} &\bbH,\bbG:\R^{3\times 3\times 3}\to \R^{3\times 3\times 3}\
\text{are  $6$th-order posit. def. and symm. tensors, $\bbH$
potential,}
%Eventually, we assume}
%\begin{rcomm} do we REALLY need that they have a potential?\end{rcomm}
\end{align}\end{subequations}
in the sense that both $\mathrm{div}\bbH\nabla e(u)$ and
$\mathrm{div}\bbG\nabla e(u)$ are $\R_{\mathrm{sym}}^{3\times
3}$-valued, and the operator div$^2\bbH\nabla e(u)$ has a potential.

The complementarity problem
\eqref{adhes-form-d2}--\eqref{adhes-form-d4} describes  general,
possibly unilateral (depending on the choice of the mapping
$\cone:\GC \rightrightarrows \R^3$) contact. Indeed, for
$\mathfrak{h}=0$ and $z=0$ (i.e., no hyperstress contribution, and
complete delamination), \eqref{adhes-form-d2}--\eqref{adhes-form-d4}
reduce to relations \eqref{b.c.}, which in turn generalize the
Signorini conditions. For later reference, we  point out that the
complementarity conditions
\eqref{adhes-form-d2}--\eqref{adhes-form-d4} are equivalent to the
subdifferential inclusion
\begin{equation}
\label{subdiff-adhes}
\partial \ind_{\cone}(\JUMP{u}{}) +
T(u,\DT{u},\theta) +z\bbA
\JUMP{u}{}-za_0'(\JUMP{u}{})\ni0\qquad\text{on $\SC$,}
\end{equation}
which features the indicator functional
 $\ind_{\cone}: L^2 (\GC;\R^3) \to [0,+\infty]$
 associated with the
multivalued $\cone : \GC \rightrightarrows \R^3$, viz.
\begin{equation}
\label{ind_D} \ind_{\cone}(v) = \int_{\GC}\ind_{\cone(x)}(v(x))\,\dd
S \quad \text{for all $v \in L^2 (\GC;\R^3)$,}
\end{equation}
and its subdifferential (in the sense of convex analysis)
 $\partial
\ind_{\cone}: L^2(\GC;\R^3) \rightrightarrows L^2(\GC;\R^3)$.

In turn,  adhesive contact results from  the complementarity
conditions~\eqref{adhes-form-d6}--\eqref{adhes-form-d8}, which
 can be reformulated as
the \emph{flow rule}
\begin{equation}
\label{e:reactivation-2}
\partial \ind_{(-\infty,0]}(\DT{z}) + \partial\ind_{[0,1]}(z) +
\mbox{$\frac12$}\bbA\JUMP{u}{}{\cdot}\JUMP{u}{}-a_0(\JUMP{u}{})-a_1(\JUMP{u}{})\ni0
\qquad\text{in $\SC$,}
\end{equation}
with the indicator functions $I_{(-\infty,0]},\, I_{[0,1]}: \R \to
[0,+\infty]$ and their (convex analysis) subdifferentials $\partial
I_{(-\infty,0]},\,
\partial I_{[0,1]}: \R \rightrightarrows \R$.
%%%%%%%%%%%%%%%%%%%%%5
 \paragraph{Some comments on the derivation of the model.}
In \cite{rr+tr}, a thorough derivation of the analogue of system
\eqref{eq6:adhes-class-form}, in the case when the dependence on the
delamination modes is included in the model,   has been developed.
Therefore, we refer the reader to \cite[Sec.\ 3]{rr+tr} for all
details, and in the next lines we just highlight the main
differences between the system considered in \cite{rr+tr} and the
present \eqref{eq6:adhes-class-form}. Namely,  here
  the free
energy is enhanced by the higher-order term $\frac12\bbH\nabla
e(u){\vdots}\nabla e(u)$, the dissipation energy in the bulk is
enhanced by $\bbG\nabla e(\DT u){\vdots}\nabla e(\DT u)$, and on the
delamination surface $\GC$ we have $a_0$ and $a_1$ depending on
$\JUMP{u}{}$.

More specifically,  we consider the  free energy, the dissipation
rate, and the (pseudo)potential of dissipative forces in the bulk
given by
\begin{subequations}
\begin{align}
%\nonumber
&\psi^{\mathrm{bulk}} (e,\nabla e,\theta)=
\frac12\bbC(e{-}\bbE\theta){:}(e{-}\bbE\theta) +\frac12\bbH\nabla
e{\vdots}\nabla e -\frac{\theta^2}2\bbB{:}\bbE-\psi_0(\theta),
%\\
\label{psi-bulk}
%&=\frac12\bbC e{:}e+\frac12\bbH\nabla e(u){\vdots}\nabla e(u)
%-\theta \bbB{:}e-\psi_0(\theta)\quad\text{ with }\ \bbB:=\bbC\bbE,
\\&
\xi^{\mathrm{bulk}}(\DT{e},\nabla\DT{e})=2\zeta_2(\DT{e},\nabla\DT{e})
\ \ \text{ with }\ \ \zeta_2(\DT{e},\nabla\DT{e}):=\frac12
\bbD\DT{e} {:}\DT{e} + \frac12 \bbG\nabla\DT{e} {\vdots}\nabla
\DT{e},
\end{align}
\end{subequations}
where  $\psi_0: (0,+\infty)\to\R$ is a strictly convex function.
  The free energy and the dissipation rate  on the contact
 surface are
\begin{subequations}
\begin{align}
\label{psi-surf} &\psi^{\mathrm{surf}}(v,z)=
\begin{cases}
%\displaystyle{\int_{\GC}\!\!
z\big(\frac12{\dela}v{\cdot}v-a_0(v)\big)
%\,\d S}
&\mbox{if $v\GE0$ and $\ 0{\le}z{\le}1$,
%a.e.~on $\GC$
}\\+\infty&\text{otherwise},\end{cases}
\\
\label{psi-surf-b}
&\xi^{\mathrm{surf}}\big(v,\DT{z}\big)=\zeta_1\big(v,\DT{z}\big)
:=\begin{cases} \displaystyle{a_1(v)\big|\DT{z}\big|}& \mbox{ if
$\DT{z}\le0$},
%a.e. in $\GC$},
\\+\infty&\text{ otherwise}.\end{cases}
\end{align}
\end{subequations}
The overall free energy and the (pseudo)potential of dissipative
forces are then
\begin{align}
\label{Phi-epos} &\Psi(u,z,\theta)= \int_{\Omega {\setminus}
\GC}\!\!\psi^{\mathrm{bulk}} (e(u),\nabla e(u),\theta)\,\d x +
\int_{\GC}\!\!\psi^{\mathrm{surf}}(\JUMP{u}{}, z)\,\d S\,,
\\
\label{dissip-rate} &\Xi(u;\DT{u},\DT{z})
%:=\int_{\overline{\Omega}}\big[\xi(\DT{e},\DT{z})\big](\d x)
=\int_{\Omega{\setminus}\GC}\xi^{\mathrm{bulk}}
\big(e(\DT{u}),\nabla e(\DT{u})\big)\,\d x
+\int_{\GC}\!\xi^{\mathrm{surf}}\big(\JUMP{u}{},\DT{z}\big)\,\d S\,,
\end{align}
respectively. Considering the specific kinetic energy
$\frac12\varrho\,|v|^2$ (with $\varrho>0$ the mass density), for all
$v\in L^2 (\Omega)$ we define the overall kinetic energy
$T_\mathrm{kin}$ and the external mechanical load $L$ by
\begin{align}
\label{extr-load}
T_\mathrm{kin}(v):=\frac{1}{2}\int_\Omega\varrho\,|v|^2\,\d x
\qquad\text{ and }\qquad\big\langle L(t),u\big\rangle=\int_\Omega
F(t){\cdot}u\d x+\int_{\Gnew}\!\!f{\cdot}u\d S.
\end{align}
% OK letter $\Lambda$??
The mechanical part of  system \eqref{eq6:adhes-class-form}, i.e.\
equations (\ref{eq6:adhes-class-form}a,c,d,f-o), is  then just the
classical formulation of the abstract evolutionary system
\begin{align}
\label{abstract-subdiff} \!\!\!\!\partial T_\mathrm{kin}(\DDT
u(t))+\partial_{(\DT u,\DT z)}\Xi(u(t);e(\DT u(t)),\DT z(t))+
\partial_{(u,z)}\Psi(u(t),z(t),\theta(t))\ni L(t) \  \text{ for
} t \in (0,T),
\end{align}
where $\partial$ denotes the (convex analysis) subdifferential of
the functionals $\Xi$ and $\Psi$, w.r.t.\ suitable topologies which
we do not specify.  The remaining equations in
\eqref{eq6:adhes-class-form} yield the heat-transfer problem, i.e.\
(\ref{eq6:adhes-class-form}b,e,p,q). Its derivation proceeds
standardly from  postulating the entropy $s$ by the so-called Gibbs'
relation $s=-\Psi'_\theta(u,z,\theta)$,  viz.\ $ \langle
s,\tilde\theta\rangle=-\Psi'_\theta(u,z,\theta;\tilde\theta)$ for
all $\tilde\theta$, where $\Psi'_\theta(u,z,\theta;\tilde\theta) $
is the directional derivative of $\Psi$ at $(u,z,\theta)$ in the
direction $\tilde\theta$. This yields the entropy in the bulk as
\begin{align}\label{spec-entropy}
s=s(\theta,e)=-\frac{\partial\psi^{\mathrm{bulk}}}{\partial\theta}(e(u),\theta)
=\bbB{:}e(u)+\psi_0'(\theta).
\end{align}
Further, we  use the so-called {\it entropy equation}
\begin{align}\label{8-3-*}
\theta\DT{s}=\xi^{\mathrm{bulk}}(e(\DT{u}))-{\rm div}\,j.
\end{align}
Substituting $\DT{s} =\bbB{:}e(\DT{u})
-\psi_0''(\theta)\DT{\theta}$, cf.~\eqref{spec-entropy}, into the
entropy equation \eqref{8-3-*} yields the {\it heat equation}
\begin{align}
%\label{heat-eq}
\label{cap} c_{\rm v}(\theta)\DT{\theta}+{\rm div}\,j=
%\xi(e(\DT{u}),\DT{z})
2\zeta_2(e(\DT{u}))-\theta\bbB{:}e(\DT{u})\quad\text{ with }\ \
c_{\rm v}(\theta)=\theta\psi_0''(\theta).
\end{align}
%with the {\it heat capacity}
%\begin{align}\label{cap}
% c_{\rm v}(\theta)=\theta\psi_0''(\theta).
%\end{align}
Hence, assuming   the constitutive relation
$j:=-\bbK(e(u),\theta)\nabla\theta$  for the heat flux, i.e.~{\it
Fourier's law} in an anisotropic medium, one obtains  the heat
equation in the form \eqref{eq6:adhes-class-form1bis}.

Similar, but simpler thermodynamics can  also be seen  on the
contact boundary by involving $\psi^{\mathrm{surf}}$ and
$\xi^{\mathrm{surf}}$. Since \eqref{psi-surf} is independent of
temperature, the ``boundary
entropy''=$-\frac{\partial}{\partial\theta}\psi^{\mathrm{surf}}$ is
simply zero, and the corresponding entropy equation reduces to
$0=\xi^{\mathrm{surf}}(\DT{z})-\JUMP{j}{}$ (as an analog of
\eqref{8-3-*}), which then results in \eqref{adhes-form-d10}.
Incorporating the analog of the phenomenological Fourier's law,
%\eqref{Fourier-law},
we arrive at~\eqref{adhes-form-d9}.

%\INSERT{HERE SOME ABSTRACT HEAT EQ. THRU AN OPERATOR $G$ AND A POTENTIAL
%$R$ VIA A TERM $G^*R'(u,\theta)G\frac1\theta$ YIELDING ENTROPY PRODUCTION
%$\langle R'(u,\theta)G\frac1\theta,G\frac1\theta\rangle$,
%ESP. IF WE USE IT LATER IN Sect.~\ref{s:new}.....? I WILL THINK ABOUT IT}

 We emphasize  that the model is thermodynamically consistent, in the sense
that it conserves the  total energy, i.e.\ here
\begin{align}\nonumber
%\begin{aligned}
  &\frac{\d}{\d t}\!\!\!\!\ddd{\int_{\Omega{\setminus}\GC}\!
\frac\varrho2|\DT{u}|^2+\frac12\bbC e(u){:}e(u) +\frac12\bbH\nabla
e(u){\vdots}\nabla e(u)+h(\theta)\,\d x}{kinetic, elastic, and
thermal energies}{}
\\&\label{total-energy}
 + \!\!\! \!\!\!\ddd{\int_{\GC}\!\! \frac z2{\dela}
\JUMP{u}{}{\cdot}\JUMP{u}{}-a_0(\JUMP{u}{}) z\,\d S}{mechanical
energy}{in the adhesive}\!\!\!
   =\!\!\!\!\ddd{\int_{\Omega}\!
%\GRM+
\FRM{\cdot}\DT{u}\,\d x}
%{power of bulk heat}{and mechanical load}
{power of bulk}{mechanical load} \!\!\!\!
+\!\!\!\ddd{\int_{\Gamma}\gRM\,\d
S+\int_{\Gnew}\!\!\fRM{\cdot}\DT{u}\,\d S}{power of surface
heat}{and mechanical load}\!\!\!,
%\end{aligned}
\end{align}
and it  satisfies the Clausius-Duhem's entropy inequality:
\begin{equation}
\label{e:clausius-duhem}
\begin{aligned}
\frac{\d}{\d t}\int_\Omega\!s\,\d x=
%\int_\Omega\! \left(\frac{{\rm div}\left(\!\bbK \nabla\theta\right)}{\theta}+
%\frac{\GRM}{\theta}\right)\,\d x =\int_\Omega\!
%\left(\frac{\bbK\nabla\theta\cdot\nabla\theta}{\theta^2}
%+\frac{\GRM}{\theta}\right)\,\d x
\int_\Omega\!\frac{{\rm div}(\bbK\nabla\theta)}{\theta}\,\d x
=\int_\Omega\!\frac{\bbK\nabla\theta\cdot\nabla\theta}{\theta^2}\,\d
x +\int_{\partial \Omega}\frac { \gRM }\theta\,\d S \ge0,
\end{aligned}
\end{equation}
as well as non-negativity of the temperature under suitable natural
conditions, cf.\ Theorem \ref{th:3.0}.

\section{Enthalpy transformation and energetic solution}\label{sec-ent-trans}
%%%%%%%%%%%%%%%%%%%%%%%%%%%%
\noindent In what follows we are going to tackle a reformulation of
the PDE system \eqref{eq6:adhes-class-form}, in which we replace the
heat equation \eqref{eq6:adhes-class-form1bis} with an
\emph{enthalpy} equation. Namely, we switch from the absolute
temperature $\theta$, to the enthalpy $\w$, defined via
 the  so-called \emph{enthalpy transformation}, viz.
\begin{align}\label{hat-c}
\w=\ent(\theta):=\int_0^\theta\!\!c_{\rm v}(r)\,\d r.
\end{align}
Thus, $\ent$ is a primitive function of $c_{\mathrm v}$, normalized
in such a way that $\ent(0)=0$.  We will assume in \eqref{30a} below
that
 $c_{\rm v}$ is strictly positive,  hence $h$ is strictly
increasing. Thus, we are entitled to
 define
\begin{align}\label{K-T}
\Theta(\w):= \left\{\begin{array}{ll} \ent^{-1}(\w) &\text{if
}\w\ge0,
\\
0              &\text{if }\w<0,\end{array}\right. \qquad\quad
\mathcal{K}(e,\w):= \frac{\mathbb K(e,\Theta(\w))}{c_{\mathrm
v}(\Theta(\w))},
\end{align}
where $\ent^{-1}$  denotes the inverse function to $\ent$.

 With transformation \eqref{hat-c} and the related \eqref{K-T},
 also taking into account the \emph{subdifferential formulations}
 \eqref{subdiff-adhes} and \eqref{e:reactivation-2},
 the PDE
system~\eqref{eq6:adhes-class-form} turns into
\begin{subequations}\label{eqsystem}\begin{align}
\label{eq6:delam-class-trans1} &\left.\begin{array}{ll}
&\hspace{-1.6em} \varrho\DDT{u}-\mathrm{div}\big(\bbD
e(\DT{u}){+}\bbC e(u){-}\bbB\Theta(\w) {-}\mathrm{div}\,\mathfrak
h\big)=\FRM, \quad\mathfrak h=\bbH\nabla e(u)+\bbG \nabla e(\DT{u})
\\[.3em]
 &\hspace{-1.6em}
\DT{\w}-\mathrm{div}\big(\mathcal{K}(e(u),\w)\nabla\w\big)=\bbD
e(\DT{u}){:}e(\DT{u})-\Theta(\w)\bbB {:}e(\DT{u})+\bbG\nabla
e(\DT{u}) {\vdots}\nabla e(\DT{u})
\end{array}\hspace{.0em}\right\}
&&\!\!\!\!\text{in }Q{{\setminus}}\SC,
\\ \label{eq6:delam-class-trans2}
&u=0 &&\!\!\!\!\text{on }\Sdir,&&
\\
\label{eq6:delam-class-trans3} & \left.\begin{array}{ll}
&\hspace{-1.6em}
(\mathcal{K}(e(u),\w)\nabla\w)\nu=g\\[.3em]
&\hspace{-1.6em}
%&&\!\!\!\!\text{on }\Snew,&&\\\label{eq6:delam-class-trans31}&
\mathcal{T}(u,\DT{u},\w)-\mathrm{div}_{{}_{\mathrm{S}}}
(\mathfrak{h} \cdot \nu)=\fRM
\end{array}\hspace{18.8em}\right\}
&&\!\!\!\!\text{on }\Snew,&&
\\
\label{eq6:delam-class-trans32} & \mathfrak{h} \Colon (\nu \otimes
\nu)=0
 &&\!\!\!\!\!\!\!\!\!\!\!\!\!\!\!\!\!\text{on }\Sdir\cup\Snew,&&
\\
\label{class-form-d-trans} &\left.\begin{array}{ll} &\hspace{-1.6em}
  \JUMP{\bbD e(\DT{u})+\bbC e(u)-\bbB\Theta(\w)
-\mathrm{div}\,\mathfrak h}{}\norm - \mathrm{div}_{{}_{\mathrm{S}}}
(\JUMP{\mathfrak{h}}{}\norm) =0
\\[.3em]
&\hspace{-1.6em}
 \mathfrak{h}^+ \Colon
(\norm \otimes \norm)=\mathfrak{h}^- \Colon (\norm \otimes \norm) =0
\\[.3em]
&\hspace{-1.6em}
\partial \ind_{\cone}(\JUMP{u}{}) +
\mathcal{T}(u,\DT{u},\w) -\mathrm{div}_{{}_{\mathrm{S}}}
(\mathfrak{h} \cdot \nu)
 + z\delam\JUMP{u}{}-za_0'(\JUMP{u}{})\ni0
\\[.3em]
&\hspace{-1.6em}
\partial \ind_{(-\infty,0]}(\DT{z}) +
\partial\ind_{[0,1]}(z) +
\mbox{$\frac12$}\bbA\JUMP{u}{}{\cdot}\JUMP{u}{}-a_0(\JUMP{u}{})-a_1(\JUMP{u}{})\ni0
\\[.3em]
&\hspace{-1.6em} \frac12\big(\mathcal{K}(e(u),\w)\nabla\w|_{\GC}^+
\!+\mathcal{K}(e(u),\w)\nabla\w|_{\GC}^-\big){\cdot}\norm +
\eta(\JUMP{u}{},z)\JUMP{\Theta(\w)}{}=0
\\[.3em]
&\hspace{-1.6em} \JUMP{\mathcal{K}(e(u),\w)\nabla\w}{}\cdot\norm
=\minus a_1(\JUMP{u}{})\DT{z}
  \end{array}\hspace{3.1em}\right\}
%    \quad
&&\!\!\!\!\!\text{on }\SC,&&
%\end{align}
%\COM{+B.C. for $\theta$}
\intertext{with}
%\begin{align}
&\mathcal{T}(u,v,\w):= T(u,v,\Theta(\w))= \big[\bbD e(v)+\bbC
e(u)-\bbB\Theta(\w) -\mathrm{div} \big(\bbH\nabla e(u)+\bbG\nabla
e(v)\big)\big]\big|_{\Gamma}\nu\hspace{-9em} &&&&
\end{align}\end{subequations}
where again  we take   as $\nu$ the unit normal to $\GC$ if
$\Gamma=\Gnew$, and   $\nu=\norm$ if $\Gamma=\GC$.
%%%%%%%%%%%%%%%%%%%%%
%%%%%%%%%%%%%%%%%%%%%%%%%%%%5
%%%%%%%%%%%%%%%%%%
\begin{remark}
\label{rmk-on-heat} \upshape The reformulation of the heat equation
\eqref{eq6:adhes-class-form1bis}, viz.\ the second of
\eqref{eq6:delam-class-trans1}, shows the advantage of the enthalpy
transformation \eqref{hat-c}. Indeed, by means of \eqref{hat-c}, the
nonlinear term $c_{\rm v}(\theta) \dot{\theta}$
 in
\eqref{eq6:adhes-class-form1bis} has been replaced by the
\emph{linear} contribution $\dot{\w}$. This  makes the second of
\eqref{eq6:delam-class-trans1} amenable to the time-discretization
procedure developed in Section \ref{s:new}. Such a procedure would
be more troublesome, if directly implemented on the heat equation
\eqref{eq6:adhes-class-form1bis}, also taking into account the
growth  conditions on the heat capacity function $c_{\rm v}$, which
we shall impose in \eqref{30b}
 (for
analogous assumptions, see \cite{tr0,tr1,rr+tr} and
\cite[Sect.5.4.2]{eck} for contact problems in
thermo-viscoelasticity).
\end{remark}
%%%%%%%%%%%%%%%
%%%%%%%%%%%%%%%
\noindent \textbf{Functional setup.} Throughout the paper, we shall
extensively exploit that, in our 3-dimensional case,
\begin{equation}
\label{e:contemb}
\begin{aligned}
& W^{2,2}(\Omega) \subset W^{1,p}(\Omega) \text{ continuously for} \
1 \leq p \leq 6,\text{ and}
\\
 &  u\mapsto u|_\Gamma:W^{2,2}(\Omega) \to C(\Gamma)  \ \text{ compactly,}
\end{aligned}
\end{equation}
with  $\Gamma = \partial \Omega$, or $\Gamma = \GC$, or $\Gamma =
\Gnew$.   Moreover,  for $\gamma \in [2,\infty)$
  we will adopt the notation
\[
\begin{aligned}\label{def-of-W}
&  W_{\Gdir}^{2,\gamma}(\Omega {\setminus} \GC ;\R^3):=\big\{ \testu
\in W^{2,\gamma}(\Omega {\setminus} \GC ;\R^3): \ \ \testu =0 \ \
\text{on $\Gdir$}\big\}\,,
%\\
%& W_{\GC}^{2,2}(\Omega {\setminus} \GC ;\R^3)=\big\{ \testu \in
%W^{2,2}(\Omega {\setminus} \GC ;\R^3)\, : \ \ \testu =0 \ \ \text{on
%$\GC$}\big\}\,.
\end{aligned}
\]
and  denote by   $\pairing{}{}{\cdot}{\cdot}$ the duality pairing
between  the spaces $W_{\Gdir}^{2,\gamma}(\Omega{\setminus}
\GC;\R^3)^*$ and $W_{\Gdir}^{2,\gamma}(\Omega{\setminus} \GC;\R^3)$.
Furthermore, in the case  $\cone(x)$ is  a linear subspace of $\R^3$
for almost all $x \in \GC$, we shall use the notation
\begin{equation}
\label{cone-w}
\begin{aligned}
 W_{\cone}^{2,2}(\Omega {\setminus} \GC ;\R^3):=\big\{ \testu \in
W_{\Gdir}^{2,2}(\Omega {\setminus} \GC;\R^3): \ \ \JUMP{\testu(x)}{}
\in \cone(x)\ \ \foraa\, x \in \GC\big\}.
\end{aligned}
\end{equation}
We will work with the space of measures
$\mathcal{M}(\overline\Omega):=C(\overline\Omega)^* $.
 Finally,  let $X$ be a (separable) Banach space:
we denote by  $\mathcal{M}([0,T];X)$, $B_{\rm w*}([0,T];X)$, and
$BV([0,T];X)$, respectively, the Banach spaces of the measures on
$[0,T]$ with values in $X$, of the   functions from $[0,T]$ with
values in $ X$ that are bounded and
 weakly* measurable (if $X$ has a predual),  and of the functions that  have  bounded variation on
 $[0,T]$.
  Notice that     the   functions in $B_{\rm w*}([0,T];X)$ and $
BV([0,T];X)$   are defined everywhere on $[0,T]$.
%\cong
%$
%\begin{rrnew} $$ \end{rrnew}
%%%%%%%%%%%%%%%%%
%%%%%%%%%%%%%%%%%%%%%%%%%

\noindent \textbf{
%Data
Loading  qualification.}
%%%%%%%%%%%%
Hereafter, the external mechanical and thermal loading  $\FRM$,
$\fRM$, and $\gRM$ will be qualified by
\begin{subequations}\label{hypo-data}
\begin{align}
  \label{eFFe1}
&\FRM \in L^1 (0,T; L^2 (\Omega; \R^3));
\\
& \label{eFFe2}\fRM \in  W^{1,1} (0,T; L^{4/3}(\Gnew;\R^3)),
\\
& \label{posg2} \gRM \in L^1 (\Sigma), \quad  \gRM \geq 0  \  \aein
\ \Sigma\,.
\end{align}
\end{subequations}

\noindent \textbf{Initial conditions  qualification.} As for the
initial data, we impose the following
\begin{subequations}
\label{hyp-init}
\begin{align}
& \label{uzero} u_0 \in W_{\Gdir}^{2,2}(\Omega{\setminus}
\GC;\R^3)\,, \quad \JUMP{u_0}{} \GE 0 \ \ \text{on $\SC$,}
\\
& \label{vzero} \DT{u}_0 \in L^2(\Omega;\R^3) \quad\,\text{ if
$\varrho>0$}\,,
\\
& \label{zzero} z_0 \in L^\infty(\GC), \qquad 0 \leq z_0 \leq 1 \ \
\text{a.e. on}\, \GC\,,
\\
& \label{wzero} \theta_0 \in L^\omega (\Omega)\,,\qquad\ \ \theta_0
\geq 0 \ \ \aein\  \Omega,
\end{align}
\end{subequations}
where $\omega$  is  as in \eqref{30b} below.

%%%%%%%%%%%%%%%%%%%
\noindent \textbf{Weak formulation.} The energetic formulation
associated with system~\eqref{eqsystem} hinges on the following
energy functional $\Phi$, which is in fact the mechanical part of
the  free energy~\eqref{Phi-epos},
 and   on the  dissipation metric  $\calD$:
%%%%%%
\begin{align}\nonumber
&\Phi(u,z):=\frac12\int_{\Omega{\setminus}\GC}\!\!\!\!
%\Big(
\bbC e(u){:}e(u)+ \bbH\nabla e(u){\vdots}\nabla e(u)
%\Big)
\,\d x +\ind_{\cone}(\JUMP{u}{}) +\int_{\GC}\!\!
z\alpha_0\big(\JUMP{u}{}\big)
%+\ind_{\cone(x)}(\JUMP{u}{})}
 + \ind_{[0,1]}(z)\,\d S,
\\&\label{8-1-k}
\hspace{12em}\text{with the abbreviation }\
\alpha_0\big(\JUMP{u}{}\big):=\frac12\dela\JUMP{u}{}{\cdot}\JUMP{u}{}
-a_0(\JUMP{u}{}),
%+\int_{\GC}\!\!\left( \frac
%z2\dela\JUMP{u}{}{\cdot}\JUMP{u}{}+ \ind_{[0,1]}(z) -
%a_0(\JUMP{u}{}) z\right)\,\d S,
\\&
\label{DISS} \calD\big(u,\DT z):=
\begin{cases}
\displaystyle{\int_{\GC}\!a_1(\JUMP{u}{})|\DT z|\,\d S}  & \text{if
$\DT z\le0$ a.e. in $\GC$,}
\\[-.2em] + \infty & \text{otherwise.}
\end{cases}
\end{align}

We are now in the position of introducing the notion of weak
solution to system \eqref{eqsystem} which shall be analysed
throughout this paper.
%%%%%%%%%%%%%%%%
\begin{definition}[Energetic solution of the adhesive contact problem]\label{def4}
\upshape
 Given a quadruple of initial data
$(u_0,\DT{u}_0,z_0,\theta_0)$ satisfying suitable conditions (cf.\
\eqref{hyp-init} later on),
 we call a triple $(u,z,\w)$ an \emph{energetic solution} to the
Cauchy problem for (the enthalpy reformulation of)
system~\eqref{eqsystem}
 if
\begin{subequations}
\label{reguu}
\begin{align}
\label{reguu1} & \!\!\! \!\!\!u \in
W^{1,2}(0,T;W_{\Gdir}^{2,2}(\Omega{{\setminus}}\GC;\R^3)),
\\
\label{reguu2}  & \!\!\! \!\!\! u \in
 W^{1,\infty}(0,T;L^2(\Omega;\R^3)) \qquad \text{if $\varrho>0$,}
\\
\label{reguz} &  \!\!\! \!\!\!z \in L^\infty (\SC)\, \cap\,
BV([0,T];L^1(\GC))\,,\  \text{
$z(\cdot,x)$ nonincreasing on $[0,T]$ for a.a. } x \in \GC, \\
 \label{reguw}  & \left.
 \begin{array}{ll} \!\!\! \!\!\! \!\!\!  \w \in
L^r(0,T;W^{1,r}(\Omega{{\setminus}}\GC))  \,\cap\,
L^\infty(0,T;L^1(\Omega))\,\cap\, B_{\rm
w*}([0,T];\Meas(\overline{\Omega}))
 \\
 \!\!\! \!\!\! \!\!\!
{\w} \in BV([0,T]; W^{1,r'}(\Omega{\setminus} \GC)^*)
\end{array}
\right\}
  \mbox{ for any
$1\le r<\frac{5}{4}$},
\end{align}
\end{subequations}
with $r'$ denoting the conjugate exponent $\frac{r}{r-1}$ of $r$,
and the triple $(u,z,\w)$  complies with:
\\
%\begin{description}
 \ITEM{(i)}{(weak formulation of the) momentum inclusion, i.e.:}
\begin{align}
\label{constraints-delam}
 &  \JUMP{u}{}\GE0 \  \ \text{on
$\SC$,} \quad \text{ and}
\\\nonumber
%\begin{aligned}
 &\varrho \int_\Omega
\DT{u}(T)\cdot\big(\testu(T){-}u(T)\big)\,\d x +\int_Q\big(\bbD
e(\DT{u})+\bbC e(u)-\bbB\Theta(\w)\big){:}e(\testu{-}u)
-\varrho\DT{u}{\cdot}\big(\DT{\testu}{-}\DT{u}\big)
\\\nonumber
&\qquad +\big(\bbH\nabla e(u){+}\bbG \nabla
e(\DT{u})\big){\vdots}\nabla e(\testu{-}u) \,\d x\d t
+\int_{\SC}\!\!
%\left( z\bbA\JUMP{u}{}{-} z a_0'(\JUMP{u}{})\right)
\alpha_0'\big(\JUMP{u}{}\big) {\cdot}\JUMP{\testu{-}u}{}\d S\d t
\\
 &\qquad
\ge\varrho \int_\Omega\!\DT{u}_0{\cdot}\big(\testu(0){-}u(0)\big)\d
x +\int_Q\!\FRM{\cdot}(\testu{-}u)\,\d x\d t
+\int_{\Snew}\!\!\fRM{\cdot}(\testu{-}u)\,\d S\d t
%\end{aligned}
\label{e:weak-momentum-variational}
\end{align}
{for all $\testu$ in
$L^2(0,T;W_{\Gdir}^{2,2}(\Omega{\setminus}\GC;\R^3))$ with
$\JUMP{\testu}{}\GE0$  on~$\SC$ and, if $\varrho>0$, also in
$W^{1,1}(0,T;L^2(\Omega;\R^3))$,}
 \ITEM{(ii)}{total energy balance}
\begin{align}\nonumber
 T_\mathrm{kin}\big(\DT{u}(T)\big)
+\Phi\big(u(T),z(T)\big)+  \int_\Omega\w(T)(\d x) =
T_\mathrm{kin}\big(\DT{u}_0\big) +\Phi\big(u_0,z_0\big)
+\int_\Omega\w_0\,\d x
\\
+\int_Q \FRM {\cdot} (\testu{-}u)  \d x \d t +\int_{\Snew} \fRM
{\cdot} (\testu{-}u)  \d S \d t +\int_{\Sigma}\gRM \, \d S  \d t,
 \label{total-energy-brittle}
\end{align}
 \ITEM{(iii)}{
%irreversibility and
semistability for a.a.~$t\in (0,T)$}
\begin{align}
%\label{irrev} & \forall\,  0 \leq t_1 \leq t_2 \leq T \quad
%z(t_2,\cdot) \leq z(t_1,\cdot) \ \  \text{a.e. in $\GC$,}\\
 \label{semistab} &\forall \tilde{z}\in
L^\infty(\GC):\qquad \Phi\big(u(t),z(t)\big)\le\Phi\big(u(t),\tilde
z\big) +\calD\big(u(t),\tilde z-z(t)\big),
\end{align}
 \ITEM{(iv)}{ (weak formulation of the) enthalpy equation:}
\begin{align}\nonumber
& \int_{\overline{\Omega}}\!\testw(T)\w(T)(\d x)   +
\int_Q\mathcal{K}(e(u),\w)\nabla\w\cdot\nabla \testw-\w\DT{
\testw}\,\d x\d t +\int_{\SC}\!\!\het(\JUMP{u}{},z)
\JUMP{\Theta(\w)}{}\JUMP{ \testw}{}\,\d S\d t
\\&\qquad\nonumber
=\int_Q\!\Big(\bbD e(\DT{u}){:}e(\DT{u})-\Theta(\w)\bbB{:}e(\DT{u})
+\bbG \nabla e(\DT{u}){\vdots}\nabla e(\DT{u}) \Big)\testw\,\d x\d t
\\
&\qquad\label{weak-heat}
-\int_{\overlineSC}\!\!\frac{\testw|_{\GC}^+{+}\testw|_{\GC}^-}2
a_1(\JUMP{u}{})\DT{z}(\d S\d t) +\int_{\Sigma}\gRM\testw\,\d S \d t
+\int_\Omega \w_0\testw(0)\,\d x
\end{align}
\ITEM{}{for all $\testw\in
C([0,T];W^{1,r'}(\Omega{\setminus}\GC))\cap
W^{1,r'}(0,T;L^{r'}(\Omega))$, where $\w_0:= h_0 (\theta_0)$, and
$\w(T)$ and $\DT{z}$ are considered as measures on
$\overline{\Omega}$ and  $\overlineSC$, respectively,
%(=heat produced by rate-independent dissipation) defined by prescribing
%its values for every  closed set of the type
% $A:=[t_1,t_2]{\times}  C \subset[0,T]\times\overline\GC$ as
}
%\begin{align}\label{meash} \varmea(A): =\begin{cases}
%\displaystyle{\int_{C}\!a_1(\JUMP{u}{})\big|z(t_1,x){-}z(t_2,x)\big|\,\d S}
%&\text{if $z(\cdot,x)$ nonincreasing on $[t_1,t_2]$ for a.a.~$x{\in}C$},\\
%+\infty&\text{elsewhere},\end{cases}\end{align}
\ITEM{(v)}{and the remaining initial conditions
%for  $u$ and $z$
(beside $\DT{u}(0)=\dot u_0$, already
%implied
 involved  in \eqref{e:weak-momentum-variational}), i.e.}
\begin{equation}
\label{init} u(0)=u_0 \quad \aein \ \Omega, \qquad z(0)=z_0 \quad
\aein \ \GC, \qquad \w(0)=\w_0 \quad  \aein \ \Omega.
\end{equation}
%\end{description}
\end{definition}
%%%%%%%%%%%%%
%%%%%%%%%%%%%%%%%%%%%%%

\begin{remark}[The weak formulation \eqref{e:weak-momentum-variational} of the
momentum inclusion] \upshape
In order to (partially)
justify \eqref{e:weak-momentum-variational} and its link with the
classical formulation (\ref{eq6:adhes-class-form}a,c,d,f-k) of the
(boundary-value problem for the) momentum inclusion, we may observe
that, upon multiplying \eqref{eq6:adhes-class-form1} by $\testu{-}u$
(with $\testu$ an admissible test function in the sense of
Definition \ref{def4}) and integrating on $Q$, one has to deal with
the term
\[
-\int_Q \mathrm{div}\,\sigma \cdot (\testu{-}u) \, \d x \d t=
-\int_0^T \!\!\int_{\Omega}
\mathrm{div}\,\sigma_{_\mathrm{KV}}{\cdot} (\testu{-}u) \, \d x \d t
+ \int_0^T\!\!\int_{\Omega} \mathrm{div}^2
\mathfrak{h}{\cdot}(\testu{-}u) \, \d x \d t,
\]
where $\sigma_{_\mathrm{KV}}$ is a placeholder for the
``Kelvin-Voigt'' stress  $ \bbD e(\DT{u}){+}\bbC e(u){-}\bbB\theta$.
The treatment of the first integral term on the right-hand side
involves a standard integration by parts. As for the second one, let
us observe (neglecting time-integration and integrating by parts
twice, with the zero Dirichlet condition on $\Gdir$), that
\begin{align}
\nonumber \int_{\Omega} \mathrm{div}^2 \mathfrak{h} \cdot
(\testu{-}u) \, \d x &= -\int_{\Omegaone} \big(\mathrm{div}\mathfrak
h\big){\Colon} \nabla(\testu{-}u)\, \d x -\int_{\Omegatwo}
\big(\mathrm{div}\mathfrak h\big){\Colon} \nabla  (\testu{-}u) \, \d
x
\\\nonumber & \qquad
+ \int_{\Gnew} (\mathrm{div}\mathfrak{h}) {\Colon} ((\testu{-}u)
\otimes \nu) \, \d S + \int_{\GC} (\mathrm{div}\mathfrak{h})^+
{\Colon} ((\testu^+{-}u^+) \otimes \norm) \, \d S
\\\nonumber
&\qquad- \int_{\GC} (\mathrm{div}\mathfrak{h})^- {\Colon}
((\testu^-{-}u^-) \otimes \norm) \, \d S
\\\nonumber&= \int_{\Omegaone} \mathfrak h {\vdots} \nabla^2  (\testu{-}u) \,
\d x +\int_{\Omegatwo} \mathfrak h {\vdots} \nabla^2  (\testu{-}u)\,
\d x
\\
&\qquad  +
 I(\Gnew,u,\testu,\nu) + I(\GC,u^+,\testu^+,\norm) -  I(\GC,u^-,\testu^-,\norm),
\label{from-tr}
\end{align}
where we have used the short-hand notation $ I(\tilde{\Gamma},
\tilde{u},\tilde{v},\tilde{\nu}) := \int_{\tilde{\Gamma}}
(\mathrm{div}\mathfrak{h}) {\Colon} ((\tilde{v}{-} \tilde{u})
\otimes \tilde{\nu}){-}\mathfrak h {\vdots} (\nabla
(\tilde{v}{-}\tilde{u})\otimes \tilde{\nu}) \, \d S $. %for
%$\tilde{\nu}=\nu$ or $n=\tilde{\norm}$.
Then, the calculations developed in
 \cite[2nd ed.,
Sect.2.4.4]{NPDE_roubicek}, \cite{tr-tve}  and based on the
decomposition $ \nabla v= \nabla_{{}_{\mathrm{S}}} v +
\frac{\partial v}{\partial n} n$ yield the following formula
\[
I(\tilde{\Gamma}, \tilde{u},\tilde{v},\tilde{\nu}) =
\int_{\tilde{\Gamma}} \left(
(\mathrm{div}\mathfrak{h}){\cdot}\tilde{\nu} {+}
\mathrm{div}_{{}_{\mathrm{S}}}(\mathfrak{h}{\cdot}\tilde{\nu} ){-}
(\mathrm{div}_{{}_{\mathrm{S}}} \tilde{\nu}) (\mathfrak{h} {\Colon}
(\tilde{\nu} {\otimes} \tilde{\nu})) \right) {\cdot}
(\tilde{v}{-}\tilde{u}) - (\mathfrak{h} {\Colon} (\tilde{\nu}
{\otimes} \tilde{\nu})){\cdot} \frac{\partial
(\tilde{v}{-}\tilde{u})}{\partial \tilde{\nu}}\, \d S,
\]
which we plug in \eqref{from-tr}. Then, we combine the resulting
integrals on $\Gnew$ and $\GC$ with the integrals derived from the
by-part integration of $\int_{\Omega} \mathrm{div}
\sigma_{_\mathrm{KV}}{\cdot}(\testu{-}u) \, \d x$,  rely on the
boundary conditions (\ref{eq6:adhes-class-form}d,f-k), take into
account the enthalpy transformation \eqref{hat-c}, and finally use
that $\mathrm{div}^2 \mathfrak{h} \cdot (\testu{-}u) =
\mathfrak{h}{\vdots} \nabla e(\testu{-}u)$, since $\mathfrak h$ is
symmetric, being so $\bbG$ and $\bbH$. In this way, we  obtain the second and
the third term on the left-hand side of \eqref{e:weak-momentum-variational}.
The remaining terms either follow from an integration by parts in time, or
are trivial.
\end{remark}
%%%%%%%%%%%%%%%

\begin{remark}[The ``weak'' formulation of the flow rule
\eqref{e:reactivation-2}] \upshape In
\cite{Miel05ERIS,MieThe99MMRI,MieThe04RIHM}, a global stability
condition combined with energy conservation was shown to provide the
correct ``weak'' formulation of  rate-independent flow rules
\cite{Miel05ERIS,MieThe99MMRI,MieThe04RIHM,MiThLe02VFRI}. Here, the
concept of energy-preserving solutions (i.e.\ of \emph{energetic
solutions}) is crucial for mathematically treating the full
thermodynamics, cf.\ Step~5 in Section~\ref{ss:4.4} below. We point
out that \eqref{total-energy-brittle} is the integrated version of
the total energy balance~\eqref{total-energy}. Here the energy
conservation involves also the mechanical equilibrium
\eqref{e:weak-momentum-variational}, and the semistability
\eqref{semistab} plays the role of the global stability condition of
\cite{Miel05ERIS,MieThe99MMRI,MieThe04RIHM}. We refer to
\cite[Prop.~3.2]{tr1} for some justification of the
energetic-solution concept in the framework of general
thermomechanical rate-independent processes. In general, energetic
solutions may exhibit unphysical jumps, but this does not occur if
the driving energy $\Phi(u,\cdot)$ is convex, as it is indeed the
case \eqref{8-1-k} considered here. Then there is also a close link
to the conventional weak definition of the flow rule
\eqref{e:reactivation-2}, see \cite{MieThe04RIHM}.
\end{remark}

%%%%%%%%%%%%%%
\begin{remark}[The weak formulation \eqref{weak-heat} of the enthalpy equation]
\upshape \label{rmk:first-comments} A few comments on the first term
on the left-hand side and on the second term on the right-hand side
of \eqref{weak-heat} are in order. First, since $\w\in BV ([0,T];
W^{1,r'}(\Omega{\setminus} \GC)^*)$, then
 for all $t \in [0,T]$ one has
 $\w(t)$ well-defined  as an element of $W^{1,r'}(\Omega{\setminus}
\GC)^*$. Combining this with the fact that   $\w\in
L^\infty(0,T;L^1(\Omega))$ one sees that even
$\w(t)\in\mathcal{M}(\overline\Omega)$ for all $t \in [0,T]$.
However, note that the function $t\mapsto \w(t)$ may jump. Second,
let us observe that
 due to \eqref{reguz},
$\DT{z}$ is a \emph{negative} Radon measure on $\overlineSC$.
 Since
we shall impose that the function $a_1: \R^3 \to \R$ is continuous
(cf.\ \eqref{anot1}), and since the map $(t,x) \mapsto
\JUMP{u(t,x)}{}$ is also continuous because of \eqref{reguu1} and
\eqref{e:contemb}, it turns out that $(t,x) \mapsto a_1 (
\JUMP{u(t,x)}{})$ is a continuous function. Thus, $a_1 (\JUMP{u}{})
\DT{z}$ is a well-defined measure on $\overlineSC$.
\end{remark}
%%%%%%%%%%%%
%%%%%%%%%%
\begin{remark}[Mechanical energy equality]
\upshape \label{rmk:mech-eq}
%Hereafter, we shall use the following
%notation
% Let us introduce the notion of total
%variation associated with the measure $a_1 (\JUMP{u}{}) \DT{z}$ on
%$\overlineSC$, namely
%the $\sup$  taken over all partitions  $t_1= s_0 < \ldots< s_k=t_2 $
%of the interval $[t_1,t_2]$.
Subtracting \eqref{weak-heat} tested by 1 from
\eqref{total-energy-brittle} reveals that energetic solutions comply
with
 the mechanical energy
 equality:
\begin{align}\nonumber
 T_\mathrm{kin}^{\varrho}\big(\DT{u}(T)\big) &
+\Phi\big(u(T),z(T)\big) +\int_Q\!\left(\bbD e(\DT{u}){:}
e(\DT{u}){+} \bbG \nabla e(\DT{u}){\vdots}\nabla e(\DT{u})
\right)\,\d x\d t +\mathrm{Var}_{\mathcal{R}}(u,z;[0,T])
%\int_{\overline{\SC}}\varmea(\d S\d t)
\\ &
=
%\le
T_\mathrm{kin}^{\varrho}\big(\DT{u}_0\big) +\Phi\big(u_0,z_0\big)
+\int_{Q}\FRM{\cdot}\DT{u}+\Theta(\w)\bbB {:} e(\DT{u})\,\d x \d t
+\int_{\Snew}\fRM{\cdot}\DT{u}\,\d S\d t, \label{mech-eq}
\end{align}
where we have used the notation
\begin{equation}
\label{var-R}
\mathrm{Var}_{\mathcal{R}}(u,z;[t_1,t_2]):=\int_{t_1}^{t_2}\int_{\overlineGC}
a_1 (\JUMP{u}{}) |\DT{z}|(\d S\d t) \qquad \text{ for }
[t_1,t_2]\subset [0,T].
\end{equation}
%%%%%%%%%%%%%%%
\end{remark}
%%%%%%%%%%%
%%%%%%%%%%%%%%%
%%%%%%%%%%%%%%%%%
\section{Main result}\label{sec-main-res}
\noindent We now enlist our conditions on the functions $c_{\rm v}$,
$\mathbb{K}$, $\eta$, $a_0$, $a_1$, and  the loading.

 %%%%%%%%%%%
 \noindent \textbf{Assumptions.}
%%%%%%%%%%%%%%
%%%%%%%%%%
We suppose that
\begin{subequations}
\label{C-K}\begin{align}
 &
\label{30a}
 c_{\rm v}:[0,+\infty)\to (0,+\infty)\ \text{
continuous},
\\&
\label{30b} \exists\,\omega_1\ge\omega
>\mbox{$\frac65$},\ c_1\ge c_0>0\
\forall\theta\in (0,+\infty):\quad c_0(1{+}\theta)^{\omega-1}\le
c_{\rm v}(\theta)\le c_1(1{+}\theta)^{\omega_1-1},
%\\\label{30e}& \text{and, further, $\omega>\frac{2d}{d{+}2}$\,,}
\\
\label{30c} &\mathbb K:\R^{3\times 3}\times\R \to\R^{3\times 3}\
\text{is  bounded, continuous, and }
%\text{ positive definite, }
\\&
\label{30dprimo}  \inf_{(e,\w,\xi)\in\R_{\rm sym}^{3\times
3}\times\R\times\R^3,\ |\xi|=1} \mathcal{K}(e,\w)\xi{:} \xi=:
\mathsf{k}>0.
\end{align}
We also require that  $\eta(x,v,\cdot)$ is a non-negative
 affine function of the
delamination parameter $z\in [0,1]$, and, following \cite{rr+tr}, we
assume that
\begin{equation}
\label{eta-affine}
\begin{aligned}
\eta(x,v,z) &=\eta_1(x,v)z +\eta_0(x,v)
% \quad \forall\, z \in [0,1],
\ \text{ for $\eta_1,\eta_0:\GC{\times}\R^3\to [0,+\infty) $
Carath\'eodory s.t.}
\\
& \exists\, C_\eta>0 \quad \forall\, (x,v) \in \GC{\times}\R^3 : \ \
|\eta_0(x,v)|+  |\eta_1(x,v)|  \leq C_\eta (|v|^{4/3} + 1);
\end{aligned}
\end{equation}
%%%%%%%%%
As for the functions $a_0: \R^3 \to \R$ and $a_1 : \R^3 \to \R $, we
suppose that
\begin{align}
 & \label{anot0}
 a_0 \in C^1 (\R^3;\R), \quad \exists\, C_{a_0},\,C_{a_0}'>0
\ \forall\, u \in \R^3:
  \quad \  |a_0 (v)| \leq C_{a_0} |v| +  C_{a_0}'
 \\
\label{anot} &
 \text{the map }  u \mapsto
\frac12\dela u{\cdot}u-a_0 (u)=:\alpha_0(u)\text{ is convex},
 \\
 \label{anot1}
 &
a_1 \in C (\R^3;\R), \quad  \exists\, C_{a_1}>0 \ \forall\, v \in
\R^3:  \quad a_1(v) \geq C_{a_1}>0.
\end{align}
\end{subequations}

\begin{remark}
\label{rmk:on-assumpt} \upshape Let us comment on conditions
\eqref{C-K}. First of all, it is immediate to deduce
from~\eqref{30b} that
\begin{equation}
\label{growthTheta}\exists\, C_{\theta}^1,\, C_{\theta}^2>0\, \ \
\forall\, w \in [0,+\infty)\, : \quad C_{\theta}^1 ( w^{1/\omega_1}
-1) \leq \Theta(w) \leq C_{\theta}^2 (w^{1/\omega} +1)\,.
\end{equation}
It obviously follows from~\eqref{30c} that
\begin{equation}
\label{growthKappa}\exists\, C_{\mathcal{K}}>0\, \ \ \forall\,\xi,\
\zeta \in \R^3\, : \quad \left| \mathcal{K}(e,\w)\xi{:} \zeta\right|
\leq C_{\mathcal{K}} |\xi||\zeta|\,.
\end{equation}
Moreover, let us observe that the functional $u \mapsto \Phi(u,z)$
is convex thanks to \eqref{anot}. Note that $a_0$ itself need not be
concave,  and the possible violation of concavity depends on the
positive-definiteness of $\dela$. Actually, we could even allow for
bigger violation (namely  for $a_0$  \emph{semi-concave}), if the
discretization scheme were slightly modified,  like for example in
\cite[2nd ed., Rem. 8.2.4]{NPDE_roubicek}.
 However, we  have chosen not
to explore this option, since in real-world applications $\dela$ is
large.
\end{remark}
%%%%%%%%%%%%%%%%%%%%%%%%%

%\begin{RCOMM}
%In fact, we might do with $\GRM \in L^1 (0,T;W^{2,2}(\Omega)^*)$ and
%$\gRM \in L^1 (0,T;H^{1/2}(\partial\Omega)^*)$, but I do not know
%whether the calculations leading to the Clausius-Duhem inequality
%would make sense in this case...
%\end{RCOMM}

%\par\noindent\textbf{Existence for the adhesive contact problem.}
%%%%%%%%%%%%%%%%%%%%%%%
%%%%%%%%%%%%%%%%%%%%%%%%%%%
\begin{theorem}[Existence for the adhesive contact problem]\label{th:3.0}
Let us assume \eqref{hypo-data}, \eqref{C-K}, \eqref{hyp-init} and\\
\begin{subequations}
\label{hypo-data-bis} \ITEM{(i)}{if $\varrho=0$,  suppose  also}
\begin{align}
\label{effebis} & \FRM \in
 W^{1,1} (0,T; L^{6/5}(\Omega;\R^3)),
\\
\label{Dirichlet-part} & \mathscr{H}^{2} \left(\partial \Omegaone
\cap \Gdir \right)>0, \quad \mathscr{H}^{2} \left(\partial \Omegatwo
\cap \Gdir \right)>0,
\end{align}
\end{subequations}
\ITEM{}{where $\mathscr{H}^{2}$ denotes the two-dimensional
Hausdorff measure, or} \ITEM{(ii)}{if $\varrho>0$, suppose also
that}
\begin{equation}
\label{varrho2} \text{$\cone(x)$ is a linear subspace of $\R^3$} \
\foraa\, x \in \GC.
\end{equation}
Then, there exists an energetic solution $(u,z,\w)$ to the adhesive
contact problem (in the sense of Definition \ref{def4}),  with the
additional regularity
\begin{equation}\label{additional-cone}
\dot{u}\in W^{1,2}(0,T;W_{\cone}^{2,2}(\Omega{\setminus}
\GC;\R^3)^*) \qquad \text{if $\varrho>0$.}
\end{equation}
Furthermore, in both cases $\varrho>0$ and $\varrho=0$, the
positivity of the initial temperature
\begin{align}\label{strict-pos}
\inf_{x\in\Omega}\theta_0=:\theta^*>0
\end{align}
implies $\inf_{(t,x)\in Q}\theta=\inf_{(t,x)\in
Q}\Theta(\w(t,x))>0$; in particular, $\theta$ is a.e.~positive on
$Q$.
%if \begin{equation}
%\label{strict-pos} \exists\, \theta^*>0   \ \ \foraa\, x \in \Omega
%\, : \ \ \theta_0(x) \geq \theta^*>0
%\end{equation}
%then $\theta(x,t)=\Theta(\w(x,t)) >0$ for almost all $(x,t) \in Q$.
\end{theorem}
%%%%%%%%%%%%%%%%%%%%%%%%%%%%%%%
%%%%%%%%%%%%%%%%%%%%%%%%%
%%%%%%%%%%%%%%%%%%%%%%
\begin{remark}
\upshape \label{rmk:4.3} The analytical reason why in the presence
of inertial terms in the momentum equation we need to restrict to
``linear'' contact conditions on $\GC $ is ultimately that, if
$\varrho>0$, only \eqref{varrho2} makes it possible to test the
(weak formulation of the) momentum equation by the velocity
$\DT{u}$. This is needed for obtaining the mechanical energy
equality \eqref{mech-eq}, which in turn is a crucial step in the
proof of Theorem \ref{th:3.0}.
%
%Let us also emphasize that the analysis of the momentum equilibrium
%equation in which inertia interacts with Signorini boundary
%conditions is remarkably difficult. It has indeed been an open
%problem for a long time. In this connection, we may mention the
%recent results obtained  in~\cite{petrov-schatzman2009,
%petrov-schatzman2010}  for the (uncoupled)
% dynamical viscoelastic equation with Signorini contact conditions
% in the one- and three-dimensional case on unbounded domains. In such
%a context, these existence  results have been proved
%  with refined Fourier analysis techniques.
\end{remark}
%%%%%%%%%%%%%
\noindent In what follows, we shall denote by the symbols $C$, $C'$
most of the (positive) constants occurring in calculations and
estimates.

%%%%%%%%%%%%%%
\section{Semi-implicit time discretization}\label{s:new}
%        ~~~~~~~~~~~~~~~~~~~~~~~~~~~~~~~~~
We perform a semi-implicit time-discretization using an equidistant
partition of $[0,T]$, with time-step $\tau>0$ and nodes
$t_\tau^k:=k\tau$, $k=0,\ldots,K_\tau$. Hereafter, given any
sequence $\{\phi^j\}_{j\ge 1}$, we will denote the {\it backward
difference ope\-ra\-tor} and its iteration by, respectively,
\begin{align}\label{not-interp0}
\dt\phi^k :=\frac{\phi^k{-}\phi^{k-1}\!\!\!}\tau, \qquad \dt^2\phi^k
:={\dt\big(\dt\phi^k\big)=}
\frac{\phi^k{-}2\phi^{k-1}{+}\phi^{k-2}\!\!\!}{\tau^2}.
\end{align}
%%%%%%%%%%%%%%%%%%%
 We approximate the data $\FRM$,
$\fRM$ by local means, i.e. setting for all  $k=1,\ldots,K_{\tau}$
\[
\FRM_{\tau}^k:= \frac{1}{\tau}\int_{t_\tau^{k-1}}^{t_\tau^k}
\FRM(s)\, \d s\,,  \qquad  \fRM_{\tau}^k:=
\frac{1}{\tau}\int_{t_\tau^{k-1}}^{t_\tau^k} \fRM(s)\, \d s\,.
\]
%%%%%%%%%%%%%%%%%%%
Furthermore, we approximate $\gRM$ by suitably constructed discrete
data $\{\gRM_{\tau}^k \}_{k=1}^{K_\tau} \subset H^{1/2}(\partial
\Omega)^*$ such that \eqref{shall-use-later} below holds, and
 the initial datum $u_0$ by a sequence  $\{
u_{0,\tau}\} \subset  W_{\Gdir}^{2,\gamma} (\Omega{\setminus} \GC;
\R^3)$ (with $\gamma>\max\{4,\frac{2\omega}{\omega-1}\}$ as assumed
in Problem~\ref{probk})   such that
\begin{equation}
\label{est-init-data}
\begin{gathered}
 \lim_{\tau \downarrow 0} \sqrt[\gamma]{\tau} \|\nabla
 e(u_{0,\tau})\|_{L^\gamma(\Omega;\R^{3\times 3 \times 3})}=0, \qquad  u_{0,\tau} \to
 u_0 \qquad\text{ in $W^{2,2}(\Omega;\R^3)\ $ as $\ \tau\to0$.}
 \end{gathered}
\end{equation}

We are now in the position of formulating the time-discrete problem,
which we again write in the classical formulation for notational
simplicity.
%%%%%%%%%%%%%%%%%%%%%%%
%%%%%%%%%%%%%%%%%%%%%%
\begin{problem}
\label{probk} \upshape {Let
$\gamma>\max\{4,\frac{2\omega}{\omega-1}\}$.} Given
\begin{align}\label{IC2}
u_{\tau}^0=u_{0,\tau},\qquad u_{\tau}^{-1}=u_{0,\tau}-\tau\DT u_0,
\qquad z_{\tau}^0=z_{0},\qquad\w_{\tau}^0=\w_0,
\end{align}
find $\{(\uk,\wk,\zk) \}_{k=1}^{K_\tau}$ fulfilling,
 for $k=1,...,K_\tau$, the recursive scheme consisting of
%\begin{itemize}
%\item[-]
the \textbf{discrete momentum equation} in $ \Omega {\setminus} \GC
$:
%\begin{equation}
\begin{subequations}\label{semi-implicit}
\begin{align}
 \label{GM1a}
 & \!\! \!\! \!\! \!\!\!\!\!
\left.
  \begin{array}{ll}
 &  \varrho \dt^2 \uk -{\rm div}\Big(\mathbb De\big(\dt
\uk\big) +\bbC e(\uk){-}\bbB\Theta(\wk) +\tau
\big|e(\uk)\big|^{\gamma-2}e(\uk) - \mathrm{div}\hk \Big)
=\FRM_{\tau}^k
\\
&\qquad\qquad\qquad\qquad\text{with}\qquad\hk=\big(\bbH
+\tau\big|\nabla e(\uk)\big|^{\gamma-2}\mathbb{I}\big) \nabla
e(\uk){+}\bbG\nabla(e(\dt\uk))
%\INSERT{+\tau\big|\nabla e(\uk)\big|^{\gamma-2}\DELETE{\bbJ} \nabla e(\uk)}
\!\!\end{array}
 \right\}  \text{ in $\Omega{\setminus}\GC$,}
%\end{equation}
\intertext{where $\mathbb{I}: \R^{3 \times 3 \times 3} \to \R^{3
\times 3 \times 3}$ denotes the 6th-order identity tensor, with the
boundary conditions}
%\label{bcu}
%&\begin{array}{lll} &  \uk = 0 &\text{{on} $\Gdir$}\,,
%\\&\Big(\mathbb De\big(\dt \uk\big) +\bbC e(\uk)-\Theta(\wk)\bbB
%-\mathrm{div} \hk&\\&\qquad
%+\tau \big|\nabla e(\uk)\big|^{\gamma-2} \DELETE{\bbJ} \nabla e(\uk)
%\Big)\nu -\mathrm{div}_{{}_{\mathrm{S}}} (\hk \cdot \nu) =\fRM_{\tau}^k
%&\text{{on} $\Gnew$}\,,\\
%&  \hk \Colon (\nu \otimes \nu)=0 &\text{{on} $\INSERT{\Gdir\cup}\Gnew\,$,}
%\end{array}
\label{bcu1}&\uk = 0 &&\hspace{-12em}\text{{on} $\Gdir$}\,,
\\\nonumber
&\Big(\mathbb De\big(\dt \uk\big) +\bbC e(\uk)-\Theta(\wk)\bbB
&&\hspace{-12em}
\\\label{bcu2}&\qquad\qquad
+\tau \big|e(\uk)\big|^{\gamma-2}e(\uk) -\mathrm{div}\,\hk
%\big(\hk+
%\tau \big|\nabla e(\uk)\big|^{\gamma-2} \DELETE{\bbJ} \nabla e(\uk)\big)
\Big)\nu -\mathrm{div}_{{}_{\mathrm{S}}} (\hk \cdot \nu)
=\fRM_{\tau}^k &&\hspace{-12em}\text{{on} $\Gnew$}\,,
\\\label{bcu3}
&  \hk \Colon (\nu \otimes \nu)=0 &&\hspace{-12em} \text{{on}
$\Gdir\cup\Gnew\,$,}
%\end{equation}
%%%%%%%%%%%%%%%%%
\intertext{and the conditions on the contact boundary}
%\begin{equation}
&\label{bccont} \!\! \!\! \!\! \!\! \left. \!\!\!\!\begin{array}{ll}
   & \JUMP{\bbD e(\dt \uk)+ \bbC
e(\uk){-}\Theta(\wk)\bbB
%+\tau\big|\nabla e(\uk)\big|^{\gamma-2}\DELETE{\bbJ} \nabla e(\uk)
+\tau\big|e(\uk)\big|^{\gamma-2}e(\uk) -\mathrm{div}(\hk)}{}\norm
\\ &
 \hspace{21em}
-\mathrm{div}_{{}_{\mathrm{S}}}  (\JUMP{\hk}{}\norm) =0\,,\!\!
\\
& (\hk)^+ \Colon (\norm \otimes \norm)=(\hk)^- \Colon (\norm \otimes
\norm) =0\,,
\\
 & \!\!\!\!\!\!\!\!
\begin{array}{ll}
    & \zkm \alpha_0'\big(\JUMP{\uk}{}\big)
%\zkm  \mathbb A \JUMP{\uk}{}     -\zkm a_0'(\JUMP{\uk}{})
 + \partial I_K (\JUMP{\uk}{}) +\Big(\bbD e(\dt \uk)+ \bbC
e(\uk){-}\Theta(\wk)\bbB
 \\ &
 \hspace{7.5em}
%+\tau\big|\nabla e(\uk)\big|^{\gamma-2}\DELETE{\bbJ}\nabla e(\uk)
+\tau\big|e(\uk)\big|^{\gamma-2}e(\uk) -\mathrm{div}\hk
 \Big)\norm
 -\mathrm{div}_{{}_{\mathrm{S}}} (\hk  \cdot \norm)
  =0,\!\!
\end{array}
\end{array}
\right\}\text{on $\GC$};
%\end{equation}
%\begin{rcomm}
%it seems to me that, to make the proof of Lemma \ref{lem-first} correct, you need to have $\zkm a_0'(\JUMP{\ukm}{})$
%and not $\zkm a_0'(\JUMP{\uk}{})$, cf. \eqref{concavity}...
%\end{rcomm}
%\item[-]
\intertext{further, the \textbf{discrete enthalpy equation}:}
%\begin{equation}
\nonumber &\dt \wk -\mathrm{div}
\big(\mathcal{K}(\wk,e(\uk))\nabla\wk\big)
 =\mathbb De\big(\dt \uk\big){:}e\big(\dt \uk\big)
\\&\label{GM3a}
\hspace{8em} -\Theta(\wk)\bbB{:}e\big(\dt \uk\big) + \bbG \nabla
e\big(\dt \uk\big){\vdots}\nabla e\big(\dt \uk\big)\quad \text{in
$\Omega{\setminus}\GC$,}
%\end{equation}
\intertext{with the boundary conditions}
%\begin{equation}
&\label{bcteta} \big(\mathcal{K}(\wk,e(\uk))\nabla\wk\big){\cdot}\nu
=\gRM_{\tau}^k \quad \text{{on} $\Gdir\cup\Gnew$,}
%\end{equation}
\intertext{and the conditions on the contact boundary}
%\begin{equation}
\label{bc-contact-teta} &\left.
\!\!\!\!\!\!\!\!\!\!\begin{array}{ll} &
 \displaystyle{\frac12\Big(\mathcal{K}(\wk,e(\uk))\nabla\wk|_{\GC}^+\!\!+
\mathcal{K}(\wk,e(\uk))\nabla\wk|_{\GC}^-\Big){\cdot}\norm}
\\&\hspace{15em}
+\eta\big({\JUMP{\ukm}{},}\zk\big)\JUMP{\Theta(\wk)}{}=0\,,
\\[.4em]
& \JUMP{\mathcal{K}(\wk,e(\uk))\nabla\wk}{}\norm = \minus
a_1\big(\JUMP{\uk}{}\big)\dt{\zk}
\end{array}
\hspace{1.7em}\right\}
\ \text{{on} $\GC$;}
%\end{equation}
%\item[-]
\intertext{ and also the \textbf{discrete flow rule} for the
delamination parameter}
%\begin{align}
 \label{BC-conta}
&\partial_v \mathcal{F}(\zkm;\zk) +\alpha_0\big(\JUMP{\uk}{}\big)
%-a_0(\JUMP{\uk}{})
-a_1\big(\JUMP{\uk}{}\big)
%+\frac1{2}\mathbb{A}\JUMP{\uk}{}\cdot \JUMP{\uk}{}
 \ni0 \quad \text{{on} $\GC$,}
\end{align}
\end{subequations}
where $\mathcal{F}(\zkm;\cdot):\R
%L^1(\GC)
\to[0,+\infty]$ is the convex functional
\begin{equation}
 \label{fkm}
\mathcal{F}(\zkm;v) =
%\int_{\GC}\left(
\ind_{(-\infty,0]}\big(\frac{v{-}\zkm}{\tau}\big)+I_{[0,1]}(v).
%\right)\, \d S.
\end{equation}
\end{problem}

 In the last condition in \eqref{bccont},  traces of the
overall stress either from $\Omegaone$ or from $\Omegatwo$ can be
considered with the same effect, thanks to   the first  boundary
condition in \eqref{bccont}.

%%%%%%%%%%%%%%%%%%%%
\noindent
%%%%%%%%%%%%%%%%%%%%%
%%%%%%%%%%%%%%%%%%%%%%%%%%%%%
\begin{remark}
\upshape \label{rmk:features} Let us highlight the main features of
the time-discrete scheme \eqref{semi-implicit}.

 First,
 the discrete version \eqref{GM3a}  of the enthalpy equation is
\emph{fully implicit}, and in particular on the right-hand side the term
$\Theta(\wk)\bbB{:}e\big(\dt \uk\big) $ appears, instead of
$\Theta(\wkm)\bbB{:}e\big(\dt \uk\big) $. This is crucial to obtain
the positivity of the temperature, i.e.\ $\vartheta_{\tau}^k\ge0$
a.e.\ in $\Omega$, cf.\ Lemma \ref{lem-exist}. Notice that the
 term  $\bbB \Theta(\wkm)$  occurs on the left-hand side of
\eqref{GM1a}: therefore, \eqref{GM1a} and \eqref{GM3a} are coupled.

Second, the boundary conditions \eqref{bc-contact-teta}  on $\GC$
for the discrete enthalpy equation involve $\zk$. In this way,
\eqref{GM3a} is coupled with \eqref{BC-conta}, hence the whole
system is coupled. Nonetheless, the mechanical part of system
\eqref{semi-implicit}   (viz., (\ref{semi-implicit}a-e,i))  can be
reformulated in terms of the subdifferential inclusion
\begin{equation}
\label{discr-subdiff}
\partial T_\mathrm{kin}(\dt^2 \uk)+
\partial_{(\DT u,\DT z)}
\Xi(\uk;e(\dt \uk),\dt \zk)+
\partial_{(u,z)}\Psi(\uk,\zkm,\wk)\ni L_\tau^k
\end{equation}
(cf.\ \eqref{abstract-subdiff}), i.e. it is \emph{semi-implicit}
w.r.t.\ the variable $z$. This will allow for some simplifications
in the a priori estimates, see Lemma \ref{lem-first}.

Third, we  have added the term $\tau\mathrm{div}^2(\tau|\nabla
e(\uk)|^{\gamma -2} \nabla
e(\uk))-\tau\mathrm{div}(\tau|e(\uk)|^{\gamma -2} e(\uk))$ to the
momentum equation in the bulk and to the corresponding
boundary/contact conditions, too. Its role is to compensate the
quadratic growth of the right-hand side  of the enthalpy equation
\eqref{GM3a} when $\gamma$ is chosen large enough,
%\emph{hyperstress}  terms in the momentum equation
cf.\ the proof
of Lemma \ref{lem-exist}.
%Clearly,
 Being premultiplied by the factor $\tau$, this
higher-homogeneity regularization will vanish when passing
$\tau\to0$. Because of this term,  we also need  to regularize the
initial condition $u_0$ in (\ref{IC2}), cf.\ \eqref{est-init-data}.
\end{remark}

\begin{remark}
\upshape \label{rem:comparison-schemes}
% the
 The  time-discrete scheme \eqref{semi-implicit} is simpler than
the one devised in \cite{rr+tr},
 where the first term on the right-hand side
of \eqref{GM3a} was multiplied by the coefficient
$(1{-}\sqrt{\tau})$, and additional terms (featuring monotone
functions of $z$ and $\JUMP{u}{}$) were added to the discrete flow
rule and to the boundary conditions on $\GC$ for $u$. Such terms
were used in the derivation of the discrete a priori estimates (in
particular, of the discrete energy inequalities) via auxiliary
minimization problems. Instead,  here we  adopt a more direct
approach in the proof of the discrete  mechanical  and total energy
inequalities, cf.\ \eqref{disc-energy0} and \eqref{disc-energy}
below. Indeed, we strongly rely on the semi-implicit character of
\eqref{discr-subdiff} and on the convexity of
 $\alpha_0$.
\end{remark}
%%%%%%%%%%%%%%%%%%%%%%%%%%%%%%%%
\begin{lemma}[Existence of weak solutions to Problem~\ref{probk}]
\label{lem-exist} Under the assumptions  of Theorem \ref{th:3.0},
for every $k=1,...,K_{\tau}$, there exists a triple $ (\uk,\zk,
\wk)\in W_{\Gdir}^{2,\gamma}(\Omega{{\setminus}}\GC;\R^3)
 \times
L^\infty(\GC)\times W^{1,2}(\Omega{{\setminus}}\GC)$, fulfilling the
weak formulation of the boundary value problem
\eqref{GM1a}--\eqref{BC-conta}.
 Moreover, $\wk\ge0 $ a.e. in
$\,\Omega$. If, in addition, \eqref{strict-pos} holds, then there
exists some constant $\chi^*>0$ such that, for sufficiently small
$\tau$,
\begin{equation}
\label{posw} \wk\ge\chi^*>0 \quad \aein \ \Omega \ \ \text{for every
$k=1,...,K_{\tau}$.}
\end{equation}
\end{lemma}
\noindent \textbf{Sketch of the proof.} We may argue along the very
same lines as in the proof of \cite[Lemma 7.4]{rr+tr}. Indeed, the
existence of a weak solution to Problem \ref{probk} follows from the
theory of pseudomonotone set-valued operators (see
e.g.~\cite[Chap.~2]{NPDE_roubicek}), and in particular from
  Leray-Lions type theorems, like~\cite[Chap.~5, Cor. 5.17]{NPDE_roubicek}.
  To apply such results, one has to verify
  the strict monotonicity of the main part of the elliptic operator,
  involved in
 the weak formulation of problem \eqref{GM1a}--\eqref{BC-conta}.
 One has also to show that this operator is coercive w.r.t.\
the norm of $W_{\Gdir}^{2,\gamma}(\Omega{{\setminus}}\GC;\R^3)\times
L^\infty(\GC)\times W^{1,2}(\Omega{{\setminus}}\GC)$. For this, the
term  $\tau{\rm div}^2( |\nabla e(\uk)|^{\gamma-2} \nabla e(\uk))
-\tau{\rm div}(|e(\uk)|^{\gamma-2}e(\uk))$ on the left-hand side of
\eqref{GM1a} plays a crucial role, in that it counteracts the
quadratic nonlinearities in $e(\uk)$ and in $\nabla e(\uk)$ on the
right-hand side of \eqref{GM3a}:
 for this , we need $\gamma > \max\{4,\frac{2\omega}{\omega-1}
\}$.
 All the calculations for proving
this  strict monotonicity and coercivity in the present setting, and
also for obtaining the strict positivity \eqref{strict-pos}, are
very similar to those carried out in the proof of  \cite[Lemma
7.4]{rr+tr}.
 Therefore, we prefer to omit all details and directly refer the reader to
\cite{rr+tr}.
 $\hfill\Box$
 %%%%%%%%%%%%%%%%%%
 %%%%%%%%%%%
 \begin{remark}
\upshape Note that the actual discrete version of the flow rule
\eqref{e:reactivation-2} for the delamination parameter would be
\begin{equation}
\label{actual-discrete}
 \partial\ind_{(-\infty,0]}\big(\dt \zk\big)+\alpha_0(\JUMP{\uk}{})
%-a_0(\JUMP{\uk}{})
-a_1(\JUMP{\uk}{})
%+\frac1{2}\mathbb{A}\JUMP{\uk}{}\cdot \JUMP{\uk}{}
+ \partial \ind_{[0,1]}(\zk) \ni0 \quad \text{{on} $\GC$,}
\end{equation}
which in fact yields a solution to \eqref{BC-conta}, since $
\partial\ind_{(-\infty,0]}\big(\dt \zk\big) +
 \partial
\ind_{[0,1]}(\zk)  \subset \partial_v \mathcal{F}(\zkm;\zk)  $,
while the converse inclusion in general does not hold.
 The reason why we have replaced \eqref{actual-discrete} by \eqref{BC-conta}
is that, differently from \eqref{actual-discrete}, the differential
inclusion \eqref{BC-conta}  features only  one nonsmooth unbounded
operator. Hence, we are entitled to directly  apply
  the
aforementioned \cite[Chap. 5, Cor. 5.17]{NPDE_roubicek} to prove
existence of solutions to Problem \ref{probk}. In turn, in view of
the analysis
 which shall be developed later on (cf.\ Lemma \ref{lem-first}),
it is actually sufficient to solve \eqref{BC-conta} in place of
\eqref{actual-discrete}.
\end{remark}
%%%%%%%%%%%%%%%%%%%%%%%%%%%%%%%%%%%
\noindent \textbf{Approximate solutions.} For $\tau>0$ fixed,
 the left-continuous and right-continuous \emph{piecewise constant},  and the
\emph{piecewise linear} interpolants of the discrete solutions
$\{\uk \}_{k=1}^{K_\tau}$
 are respectively the functions
 $\pwc u{\tau} : (0,T)
\to   W_{\Gdir}^{2,\gamma}(\Omega {\setminus} \GC ;\R^3)$,  $\upwc
u{\tau} : (0,T) \to  W_{\Gdir}^{2,\gamma}(\Omega {\setminus} \GC
;\R^3)$,  and
 $\pwl u{\tau} : (0,T) \to  W_{\Gdir}^{2,\gamma}(\Omega {\setminus} \GC ;\R^3)$
  defined by $
\pwc u{\tau}(t) = \uk,$  $\upwc u{\tau}(t) = \ukm,$ $\pwl u{\tau}(t)
=\frac{t-t_\tau^{k-1}}{\tau} \uk + \frac{t_\tau^k-t}{\tau} \ukm $
for $t \in (t_\tau^{k-1}, t_\tau^k]$. In the same way, we shall
denote by $\pwc{\w}{\tau}$, $\upwc{\w}{\tau}$, and $\pwc{z}{\tau}$,
the piecewise constant interpolants of the elements
$\{\wk\}_{k=1}^{K_\tau} $ and $\{\zk\}_{k=1}^{K_\tau} $, and  the
related piecewise linear interpolants by $\wt$ and $\zt$.
Furthermore, we shall use the notation
$\overline{\mathit{t}}_{\tau}$ and $\underline{\mathit{t}}_{\tau}$
for the left-continuous and right-continuous piecewise constant
interpolants associated with the partition, i.e.\
 $\bar{\mathit{t}}_{\tau}(t) = t_\tau^k$ if $t_\tau^{k-1}<t \leq t_\tau^k $
and ${\underline{\mathit{t}}}_{\tau}(t)= t_\tau^{k-1}$ if
$t_\tau^{k-1} \leq t < t_\tau^k $.

We shall also  consider  the interpolants $\pwc \FRM{\tau}$, $\pwc
\fRM{\tau}$, and $\pwl \fRM{\tau}$, of the $K_\tau$-tuples $\{ \FRM_{\tau}^k
\}_{k=1}^{K_\tau}$, $\{\fRM_{\tau}^k \}_{k=1}^{K_\tau}$.
In view of~\eqref{eFFe1}--\eqref{eFFe2} and \eqref{effebis},
the following estimates and strong convergences hold as $\tau \to 0$:
\begin{subequations}
\label{data-converg}
\begin{align}
& \label{data-converg-1}
 \pwc \FRM{\tau} \to \FRM   \ \left\{
\begin{array}{ll}
\text{in $L^1(0,T;L^2(\Omega;\R^3))$}&\text{if $\varrho>0$,}
\\
\text{in $L^p(0,T;L^{6/5}(\Omega;\R^3))$ for all $1\le
p<\infty$}&\text{if $\varrho=0$;}
\end{array}
\right.
\\
& \label{data-converg-2}
\begin{aligned}
& \exists\, C>0 \ \  \forall\, \tau>0\, : \ \  \| \pwc
\fRM{\tau}\|_{L^\infty (0,T; L^{4/3}(\Gnew;\R^3))} \leq C \|
 \fRM\|_{L^\infty (0,T;L^{4/3}(\Gnew;\R^3))}\,,
 \\
& \pwc \fRM{\tau} \to \fRM \quad \text{in $ L^p
(0,T;L^{4/3}(\Gnew;\R^3)) $ for all $1\leq p <\infty$
 as $\tau \to
0$}\,,
\\
& \exists\, C>0 \ \  \forall\, \tau>0\, : \ \   \| \pwl
{\DT{\fRM}}{\tau}\|_{L^1 (0,T;L^{4/3}(\Gnew;\R^3))} \leq 2 \|
\DT{\fRM}\|_{L^1 (0,T;L^{4/3}(\Gnew;\R^3))}\,,
\\
& \exists\, C>0 \ \  \forall\, \tau>0\, : \ \   \| \pwl
{\DT{\FRM}}{\tau}\|_{L^1 (0,T;L^{6/5}(\Omega;\R^3))} \leq 2 \|
\DT{\FRM}\|_{L^1 (0,T;L^{6/5}(\Omega;\R^3))}\,,
 \end{aligned}
\end{align}
\end{subequations}
Finally, we shall construct the discrete data $\{\gRM_{\tau}^k
\}_{k=1}^{K_\tau} \subset H^{1/2}(\partial \Omega)^*$ in such  a way
that the related piecewise constant interpolants $\pwc \gRM{\tau} $
fulfill
\begin{align}
\label{shall-use-later}
 \pwc \gRM{\tau} \to \gRM \qquad \text{ in $L^1 (\Sigma)\ $ as $\ \tau\to0$}\,.
\end{align}
%%%%%%%%%%%%%%%%%%%%%%
%\paragraph{Weak formulation of equations~\eqref{GM1a}, \eqref{GM3a}.}
 Using   the interpolants so far introduced, we now state the
\emph{discrete} versions of the weak formulation
\eqref{e:weak-momentum-variational} of the momentum inclusion, the
total energy balance \eqref{total-energy-brittle}, the semistability
\eqref{semistab}, the weak formulation \eqref{weak-heat} of the
enthalpy equation.
   For the momentum inclusion,
we introduce
  ``discrete test functions'', viz.\ $K_\tau$-tuples
$\{\testu_\tau^k\}_{k=1}^{K_\tau} \subset
W_{\Gdir}^{2,2}(\Omega{\setminus} \GC; \R^3)$, fulfilling
$\JUMP{\testu_\tau^k}{} \GE 0$ on $\GC$,
% and $\{\testw_\tau^k\}_{k=1}^{K_\tau} \subset W^{1,2}(\Omega{\setminus}\GC)$,
and we denote by $\pwc \testu \tau$ and
%,
$\pwl \testu \tau$
%, and $\pwc \testw \tau$
their interpolants.
Furthermore, referring to the definition \eqref{8-1-k} of $\Phi$, we
shall use the notation
\begin{align}
\label{Phi-eps-tau}
\begin{aligned}
\Phi_{\tau}(u,z):  = \Phi(u,z) + \frac\tau\gamma
\int_{\Omega{\setminus}\GC} |e(u)|^\gamma+|\nabla e(u)|^\gamma\,\d
x.
\end{aligned}
\end{align}
Hence, the approximate solutions $(\pwc u\tau, \upwc u\tau,\pwc
\w\tau, \upwc \w\tau, \pwc z\tau,\pwl u\tau,\pwl\w\tau,\pwl z\tau)$
fulfill \textbf{the discrete (weak) momentum inclusion}
\begin{align}
 \nonumber \int_{Q}  & \bigg(  \big(\bbD e({\pwl{\DT{u}}{\tau}})
  +\bbC e(\pwc u{\tau})-\bbB\Theta(\pwc\w{\tau})
+\tau |e(\pwc u {\tau})|^{\gamma-2}e(\pwc u {\tau}) \big){:}
e(\pwc\testu{\tau}{-}\pwc u{\tau})
\\&\quad\nonumber
%+\Big(\tau |\nabla e(\pwc u {\tau})|^{\gamma-2} \DELETE{\bbJ} \nabla e(\pwc u {\tau})
+\big((\bbH +\tau |\nabla e(\pwc u {\tau})|^{\gamma-2}\mathbb{I})
\nabla e({\pwc{u}{\tau}}) + \bbG\nabla e({\pwl{\DT{u}}{\tau}})\big)
{\vdots}\nabla e(\pwc \testu{\tau}{-} \pwc u{\tau}) \bigg)\, \d x\d
t
\\ & \quad\nonumber
+ \int_{\SC} \upwc{z}\tau\alpha_0'( \JUMP{\pwc {u}\tau}{})
%\left(\upwc {z}\tau \mathbb{A} \JUMP{\pwc {u}\tau}{}  -
%\upwc {z}\tau a_0' ( \JUMP{\pwc {u}\tau}{})\right)
{\cdot} \JUMP{\pwc\testu\tau{-} \pwc u\tau}{} \, \d x\d t
-\int_{\tau}^T\!\!\int_{\Omega}
\varrho{\pwl{\DT{u}}{\tau}}(\cdot-\tau) {\cdot}
({\pwl{\DT{\testu}}{\tau}{-}\pwl{\DT{u}}{\tau}}) \, \d x\d t
\\ &\quad \nonumber
+ \int_{\Omega}\varrho
{\pwl{\DT{u}}{\tau}}(T){\cdot}({\pwl{\testu}{\tau}} (T){-}
{\pwl{u}{\tau}} (T))\, \d x
\\& \label{e:discrmom} \geq  \int_{\Omega} \varrho
{\DT{u}_0}{\cdot} (\pwl{\testu}{\tau}(\tau){-}
\pwl{u}{\tau}(\tau))\, \d x + \int_{Q}  \pwc \FRM{\tau} {\cdot}
(\pwc{\testu}{\tau}-\pwc{u}{\tau})  \, \d x\d t + \int_{\Snew}
\pwc\fRM{\tau} {\cdot}
 (\pwc{\testu}{\tau}-\pwc{u}{\tau}) \,\d S  \d t;
\end{align}
\textbf{the discrete total energy inequality}
\begin{align}\nonumber
& T_\mathrm{kin}\big(\pwl{\DT{u}}{{}\tau}(t)\big)
+\Phi_{{}\tau}\big(\pwc{u}{{}\tau}(t),\pwc{z}{{}\tau}(t)\big)+\int_\Omega\pwc{\w}{{}\tau}(t)\,\d
x
%\\& \nonumber
\le T_\mathrm{kin}\big({\DT{u}_0})
+\Phi_{{}\tau}\big(u_{0,\tau},z_{0}) \nonumber
\\
& \hspace{3em}+\int_\Omega\w_0\,\d x
+\int_0^{\bar{\mathit{t}}_{{}\tau}(t)}\!\! \bigg(\int_{\Omega} \pwc
\FRM{\tau}{\cdot} \pwl{\DT{u}}{{}\tau}\,\d x
+\int_{\Gnew}\!\!\pwc\fRM{\tau} {\cdot}\pwl{\DT{u}}{{}\tau}\,\d S
+\int_{\partial \Omega} \pwc \gRM{\tau}\,\d S \bigg)\, \d s;
 \label{disc-energy}
 \end{align}
(see \eqref{extr-load} for the definition of  $T_\mathrm{kin}$),
\textbf{the discrete semistability} for a.a. $t \in (0,T)$
\begin{align}\label{disc-semistability}
\Phi_{\tau}\big(\pwc{u}{{}\tau}(t),\pwc{z}{{}\tau}(t)\big)
\le\Phi_{{}\tau}\big(\pwc{u}{{}\tau}(t),\tilde
z\big)+\calR\left(\JUMP{\pwc u{\tau}(t)}{},\tilde{z}-
\pwc{z}{\tau}(t) \right)
 \quad \text{for all $\tilde z\in
L^\infty(\GC)$};
\end{align}
\textbf{the discrete (weak) enthalpy equation}
\begin{align}
 \nonumber
&\int_\Omega\!\pwl\w{\tau}(T)\pwl\testw{}(T)\, \d x
+\int_{Q}\mathcal{K}(e(\pwc
u{\tau}),\pwc\w{\tau})\nabla\pwc\w{\tau}{\cdot} \nabla
\pwl\testw{}-\pwl\w{\tau}{\pwl{\DT{\testw}}{}}\, \d x\d t
 +\int_{\SC}\!\!\!\!\het({\JUMP{\upwc u{\tau}}{},}
\pwc z{\tau})
\JUMP{\Theta(\pwc\w{\tau})}{}\JUMP{\pwl\testw{}}{}\,\d S\d t
\\ &\qquad\qquad
\nonumber
=\int_Q\!\left(\bbD e(\pwl{\DT{u}}{\tau}) {:}
e(\pwl{\DT{u}}{\tau})-\Theta(\pwc\w{\tau})\bbB {:}
e(\pwl{\DT{u}}{\tau}) + \bbG \nabla e(\pwl{\DT{u}}{\tau}) {\vdots}
\nabla e(\pwl{\DT{u}}{\tau})  \right)\pwl \testw{}\,\d x\d t
\\
 &\qquad\qquad\quad
 \label{weak-heat-discr}
-\int_{\SC}\!\!a_1(\JUMP{\pwc u{\tau}}{})\pwl{\DT{z}}{\tau}\frac{
\pwl \testw{}|_{\GC}^+{+} \pwl \testw{}|_{\GC}^-}2 \, \d S\d t
+\int_\Omega\!\w_0  {\pwl \testw{}}(0)\,\d x
+\int_{\Sigma}\pwc\gRM{\tau} \pwl \testw{}\,\d S \d t\,.
\end{align}
%
%\begin{align}\nonumber
%&\int_\Omega\!\pwl\w{\tau}(T)\pwl\testw{\tau}(T)\, \d x
%+\int_{Q}\mathcal{K}(e(\pwc
%u{\tau}),\pwc\w{\tau})\nabla\pwc\w{\tau}{\cdot} \nabla
%\pwc\testw{\tau} +\int_{\SC}\!\!\!\!\het({\JUMP{\upwc u{\tau}}{},}
%\pwc z{\tau})
%\JUMP{\Theta(\pwc\w{\tau})}{}\JUMP{\pwc\testw{\tau}}{}\,\d S\d t
%\\ &\qquad\nonumber  -\int_\tau^T\!\!\!
%\int_{\Omega}\!\upwc\w{\tau}{\pwl{\DT{\testw}}{\tau}}\, \d x\d t
%=\int_Q\!\left(\bbD e(\pwl{\DT{u}}{\tau}) {:}
%e(\pwl{\DT{u}}{\tau})-\Theta(\pwc\w{\tau})\bbB {:}
%e(\pwl{\DT{u}}{\tau}) + \bbG \nabla e(\pwl{\DT{u}}{\tau}) {\vdots}
%\nabla e(\pwl{\DT{u}}{\tau})  \right)\pwc \testw{\tau}\,\d x\d t
%\\ &\qquad\label{weak-heat-discr}
%-\int_{\SC}\!\!a_1(\JUMP{\pwc u{\tau}}{})\pwl{\DT{z}}{\tau}\frac{
%\pwc \testw{\tau}|_{\GC}^+{+} \pwc \testw{\tau}|_{\GC}^-}2 \, \d S\d
%t +\int_\Omega\!\w_0  {\pwl \testw{\tau}}(\tau)\,\d x
%+\int_{\Sigma}\pwc\gRM{\tau} \pwc \testw{\tau}\,\d S \d t\,.
%\end{align}
 with $w$ qualified as in \eqref{weak-heat}.  Inequality
\eqref{e:discrmom} can be obtained from
 (\ref{semi-implicit}a-e), by using a suitable discrete
``by-part'' summation formula, cf.~\cite[Formula (4.49)]{tr1}.
%; in the same way, \eqref{weak-heat-discr} is a reformulation of
 (\ref{semi-implicit}f-h).
%With Lemma \ref{lem-first} below,
 We now  prove \eqref{disc-energy} and \eqref{disc-semistability}.

\begin{lemma}[Approximate energetics]\label{lem-first}
 Let $\varrho\ge0$. Under the assumptions of Theorem
\ref{th:3.0},  for all $\tau>0$ the approximate solutions $(\pwc
u\tau, \upwc u\tau,\pwc \w\tau, \pwc z\tau, \pwl
u\tau,\pwl\w\tau,\pwl z\tau)$
fulfill the following discrete mechanical energy inequality
\begin{align}
\nonumber &
T_\mathrm{kin}\big(\pwl{\DT{u}}{\tau}(t)\big)
+\Phi_{\tau}\big(\pwc{u}{{}\tau}(t),\pwc{z}{{}\tau}(t)\big)
\\&\quad
\nonumber +\int_0^{\bar{\mathit{t}}_{{}\tau}(t)}\!\!  \bigg(
\int_{\Omega}\bbD e\big(\pwl{\DT{u}}{{}\tau}\big){:}
e\big(\pwl{\DT{u}}{{}\tau}\big) + \bbG \nabla
e\big(\pwl{\DT{u}}{{}\tau}\big){\vdots} \nabla
e\big(\pwl{\DT{u}}{{}\tau}\big) \,\d x
 +\int_{\GC}\!\!\zeta_1\big(\JUMP{\pwc
u\tau}{},\pwl{\DT{z}}{{}\tau}\big)\,\d S\bigg) \d s
\\&\quad
%\nonumber
\le
T_\mathrm{kin}\big({\DT{u}_0}) +\Phi_{{}\tau}\big(u_{0,\tau},z_{0})
%\\
%&\qquad
+\int_0^{\bar{\mathit{t}}_{{}\tau}(t)}\!\!\bigg(\int_\Omega
\Theta(\pwc{\w}{{}\tau})\bbB{:} e\big(\pwl{\DT{u}}{{}\tau}\big) \,\d
x + \int_{\Omega} \pwc \FRM{\tau} {\cdot} \pwl{\DT{u}}{{}\tau}\,\d x
+\int_{\Gnew} \pwc \fRM{\tau} {\cdot} \pwl{\DT{u}}{{}\tau}\,\d S
\bigg)\, \d s, \label{disc-energy0}
\end{align}
(see \eqref{psi-surf-b}  for the definition of $\zeta_1$),
the discrete total energy inequality \eqref{disc-energy}, and
the discrete semistability \eqref{disc-semistability}.
\end{lemma}
\noindent{\it Proof.} Preliminarily, we observe that
\eqref{BC-conta} is the Euler-Lagrange equation for the minimum
problem
\begin{equation}
\label{minz} \zk \in \mathrm{Argmin}_{z \in L^\infty(\GC)}
%\left(
\int_{\GC}
%\bigg(
z\alpha_0\big(\JUMP{\uk}{}\big)
%\frac12 z\mathbb{A}\JUMP{\uk}{}\cdot \JUMP{\uk}{}+
% -a_0(\JUMP{\uk}{})z-a_1(\JUMP{\uk}{})z
%\bigg)
+ \mathcal{F}(\zkm;z)\, \d S,
%\right),
\end{equation}
where $\mathcal{F}(\zkm;z)$ is as in \eqref{fkm}. Therefore, we have
\begin{align}\nonumber
%\begin{equation}\label{conseq-min}\begin{aligned}
\int_{\GC}\!\!
%\left(
I_{(\infty,0]}\Big(\frac{\zk{-}\zkm}{\tau}\Big)
+\zk\alpha_0(\JUMP{\uk}{})
%  -a_0(\JUMP{\uk}{}) \zk
-\zk a_1(\JUMP{\uk}{})
%+ \frac12 \zk\mathbb{A}\JUMP{\uk}{}{\cdot} \JUMP{\uk}{}
+ I_{[0,1]}(\zk)
%\right)
\,\d S
\\\label{conseq-min}
 \leq \int_{\GC}\!\!
%\left(
\zkm\alpha_0(\JUMP{\uk}{})
%-a_0(\JUMP{\uk}{}) \zkm
-\zkm a_1(\JUMP{\uk}{})
%+ \frac12 \zkm \mathbb{A}\JUMP{\uk}{}{\cdot}\JUMP{\uk}{}
%\right)
\,\d S.
%\end{aligned}\end{equation}
\end{align}
%\end{aligned}\end{equation}

To prove \eqref{disc-energy0}, we test  the boundary-value problem
\eqref{GM1a}--\eqref{bccont}  by $\uk-\ukm$. We add to the resulting
inequality the previously observed \eqref{conseq-min}, thus
obtaining
\begin{align}\nonumber
 & \varrho  \int_{\Omega}\dt^2 \uk {\cdot} \dt \uk \, \d x
+ \int_{\Omega} \mathbb De\big(\dt \uk\big)\Colon e\big(\dt \uk\big)
+  \bbG\nabla e\big(\dt \uk\big) {\vdots}\nabla e\big(\dt \uk\big)
\, \d x\\\nonumber & \qquad
 + \int_{\Omega} \bbC e(\uk)\Colon
e\big(\dt\uk \big)  {+}\bbH \nabla e(\uk){\vdots}\nabla e\big(\dt\uk
%-\ukm
\big) \, \d x
\\\nonumber  & \qquad
+ \tau\int_{\Omega}
|e(\uk)|^{\gamma-2}e(\uk){\colon}e\big(\dt\uk\big)+ |\nabla
e(\uk)|^{\gamma-2}\nabla e(\uk){\vdots} \nabla e\big(\dt\uk
%-\ukm
\big) \, \d x
\\\nonumber  & \qquad
%\\& \qquad
+ \int_{\GC}\!\!\zkm\alpha_0'(\JUMP{\uk}{}){\cdot}\JUMP{\dt\uk}{}
-(\dt\zk)\alpha_0(\JUMP{\uk}{})\,\d S -\int_{\GC}\!\!
a_1(\JUMP{\uk}{})(\dt\zkm ) \, \d S
%\left(\zkm \mathbb{A}\JUMP{\uk}{} \JUMP{\dt\uk
%%-\ukm
%}{} + \frac12 (\dt\zk)\mathbb{A}\JUMP{\uk}{}{\cdot}\JUMP{\uk}{}
%%-\frac12\zkm \mathbb{A}\JUMP{\uk}{}\cdot \JUMP{\uk}{}
%\right) \, \d S
%\\\nonumber & \qquad
%-\tau\int_{\GC}\!\!\left(  \zkm a_0'(\JUMP{\uk}{}) \cdot \JUMP{\dt\uk
%%-\ukm
%}{}
%%+ \zkm a_0(\JUMP{\uk}{})- \zk a_0(\JUMP{\uk}{})
%-(\dt\zk a_0)(\JUMP{\uk}{})\right)\, \d S
%\doteq
\\\nonumber & \quad
=:I_1+I_2+I_3+I_4+I_5+I_6
%+I_8
\\ & \quad
 \le
 \tau \int_{\Omega} \mathbb{B} \Theta(\wk) \Colon  e(\dt \uk) \, \d x
 +\tau \int_{\Omega} \FRM_{\tau}^k {\cdot}\dt \uk
\, \d x + \tau \int_{\Gnew} \fRM_{\tau}^k{\cdot}\dt \uk \, \d S.
\label{7.7}
\end{align}
%\end{aligned}\end{equation}
Now, we estimate the terms $I_i$, $i=1,\ldots, 8$. First of all, we observe that
\begin{equation}
\label{7.8} I_1 \geq \frac{\varrho}{2\tau} \int_{\Omega} |\dt
\uk|^2\d x - \frac{\varrho}{2\tau} \int_{\Omega} |\dt \ukm|^2\d x =
\frac{\varrho}2\dt\int_{\Omega} |\dt \uk|^2\d x.
\end{equation}
Clearly, upon summation   $I_2$ yields the third summand on the
right-hand side of \eqref{disc-energy0}, whereas
 we observe that
%%\begin{equation}\label{7.9}\begin{array}{llll}
%\begin{subequations}\label{7.9}\begin{align}
%& I_3 \ge\dt \int_{\Omega}\frac12\bbC e(\uk){:}e(\uk) \, \d x ,
%%-\frac12\int_{\Omega}\bbC e(\ukm){:}e(\ukm) \, \d x
%\\& I_4 \ge
%%  \left( \frac{\tau}\gamma\int_{\Omega}
%%|\nabla e(\uk)|^{\gamma-2}\DELETE{\bbJ} \nabla e(\uk) {\vdots} \nabla e(\uk)   \, \d x   -
%%  \frac{\tau}{\gamma} \int_{\Omega}|\nabla e(\ukm)|^{\gamma-2}  \DELETE{\bbJ} \nabla e(\ukm) {\vdots} \nabla e(\ukm) \, \d x \right)
%\dt\int_{\Omega}\frac{\tau}\gamma
%|\nabla e(\uk)|^{\gamma-2}\DELETE{\bbJ}\nabla e(\uk){\vdots}\nabla e(\uk)\, \d x,
%\\&I_5\ge\dt
%\int_{\Omega} \frac12\bbH \nabla (e(\uk)){\vdots}\nabla (e(\uk))  \, \d x ,
%%-\ \frac12\int_{\Omega} \bbH \nabla (e(\ukm)){\vdots}\nabla (e(\ukm))  \, \d x,
%\end{align}\end{subequations}
%%\end{array}\end{equation}
\begin{equation}\label{7.9}
\left.
\begin{array}{lll}
I_3 &  \ge & \displaystyle{\dt \!\!\int_{\Omega}\frac12\big( \bbC
e(\uk){:}e(\uk){+} \bbH \nabla (e(\uk)){\vdots}\nabla (e(\uk))\big)
\, \d x,}
\\[.9em]I_4 & \ge &  \displaystyle{\dt\!\!\int_{\Omega}
\frac{\tau}\gamma|e(\uk)|^\gamma+ \frac{\tau}\gamma|\nabla
e(\uk)|^\gamma \, \d x,}
\end{array}
\right\}
\end{equation}
where we have used elementary convex-analysis inequalities.
%Similarly, we have
%\begin{equation}
%\label{7.10}
%\begin{aligned}
%I_6\ge\frac12\dt\int_{\GC}\zk\mathbb{A}\JUMP{\uk}{}{\cdot}\JUMP{\uk}{}\, \d S.
%%- \frac12\int_{\GC} \zkm \mathbb{A}\JUMP{\ukm}{}\cdot \JUMP{\ukm}{} \, \d S.
%\end{aligned}
%\end{equation}
Now, it follows from the convexity of $\alpha_0$, cf. \eqref{anot},
that
\[
\int_{\GC}\!\!\zkm \alpha_0'(\JUMP{\uk}{}){\cdot} \JUMP{\dt\uk
%-\ukm
}{}\, \d S \geq  \frac1{\tau}  \int_{\GC}\!\!\zkm
\alpha_0(\JUMP{\uk}{})\, \d S - \frac1{\tau} \int_{\GC}\!\!\zkm
\alpha_0(\JUMP{\ukm}{})\, \d S .
\]
Hence, taking into account the cancellation of the term $\int_{\GC}
\zkm\alpha_0(\JUMP{\uk}{})$, for  $I_5$ we conclude the following
inequality
\begin{equation}
\label{7.11}
\begin{aligned}
I_5 \geq
%\int_{\GC}\zkm a_0(\JUMP{\ukm}{})\, \d S
\dt \int_{\GC}\!\! \zk\alpha_0(\JUMP{\uk}{})\, \d S.
\end{aligned}
\end{equation}
Combining \eqref{7.7}--\eqref{7.11},  rearranging terms and
 multiplying by $\tau$,
we obtain
%\[\begin{aligned}
\begin{align}\nonumber
& \int_{\Omega}\frac{\varrho}2  |\dt \uk|^2 + \tau
%\int_{\Omega}\left(
\mathbb De\big(\dt \uk\big) \Colon e\big(\dt \uk\big) + \tau\mathbb
G \nabla e\big(\dt \uk\big) {\vdots} \nabla e\big(\dt \uk\big)
%\right)
+\frac12\bbC e(\uk){:}e(\uk)
%\, \d x
 \\\nonumber  & \qquad
%+\int_{\Omega}\, \d x
+ \frac{\tau}\gamma|e(\uk)|^\gamma+\frac{\tau}\gamma|\nabla
e(\uk)|^\gamma
%\frac{\tau}\gamma
%%\int_{\Omega}
%|\nabla e(\uk)|^{\gamma-2}\DELETE{\bbJ} \nabla e(\uk) {\vdots} \nabla e(\uk)
% \, \d x
+
%\int_{\Omega}
\frac12\bbH \nabla (e(\uk)){\vdots}\nabla (e(\uk))  \, \d x
 %\\\nonumber  & \qquad
 +\int_{\GC}\!\!\!
%\frac12\zk\mathbb{A}\JUMP{\uk}{}\cdot \JUMP{\uk}{}
%%\, \d S
%-
%%\int_{\GC}
\zk \alpha_0(\JUMP{\uk}{})\, \d S
 \\\nonumber  & \quad \leq
\int_{\Omega}\frac{\varrho}2|\dt \ukm|^2 + \frac12\int_{\Omega}\bbC
e(\ukm){:}e(\ukm)
%\, \d x
+\frac{\tau}\gamma|e(\ukm)|^\gamma +\frac{\tau}\gamma|\nabla
e(\ukm)|^\gamma
% + \frac{\tau}\gamma
%%\int_{\Omega}
%|\nabla e(\ukm)|^{\gamma-2}\DELETE{\bbJ} \nabla e(\ukm) {\vdots} \nabla e(\ukm)   %\, \d x
 \\\nonumber  & \qquad  +
%\int_{\Omega}
\frac12\bbH \nabla (e(\ukm)){\vdots}\nabla (e(\ukm))  \, \d x
 +\int_{\GC}\!\!
%\frac12\zkm\mathbb{A}\JUMP{\ukm}{}{\cdot}\JUMP{\ukm}{}
%%\, \d S
%-
%%\int_{\GC}\!\!
\zkm \alpha_0(\JUMP{\ukm}{})\, \d S
\\ & \qquad
 +   \tau \int_{\Omega} \mathbb{B} \Theta(\wk) \Colon  e(\dt \uk) \, \d x
 +\tau \int_{\Omega} \FRM_{\tau}^k {\cdot} \dt \uk
\, \d x + \tau \int_{\Gnew} \fRM_{\tau}^k {\cdot} \dt \uk \, \d S.
\end{align}
%\end{aligned}\]
Summing over the index $k$,
 we conclude \eqref{disc-energy0}.

In order to obtain \eqref{disc-energy} for a  fixed $t \in (0,T)$,
we test \eqref{weak-heat-discr}, integrated on the time-interval
$(0,\pwc {\mathit{t}}{\tau} (t))$, by $1$, and add the resulting
relation to the mechanical energy equality \eqref{disc-energy0}.
Taking into account all cancellations, we immediately conclude
\eqref{disc-energy}.

%\\

Eventually, from \eqref{conseq-min} and the degree-1 homogeneity of
$\calR(\JUMP{\uk}{},\cdot)$, it also follows that
%\[\begin{aligned}
\begin{align}\nonumber
\Phi_\tau (\uk,\zk)&\le\Phi_\tau (\uk,\tilde{z})
-\int_{\GC}\!\!a_1(\JUMP{\uk}{}) (\tilde{z}{-}\zkm)\, \d S-
\int_{\GC}a_1(\JUMP{\uk}{}) (\zkm{-}\zk)\,\d S
\\&\le
\Phi_\tau (\uk,\tilde{z})-\int_{\GC}\!\! a_1(\JUMP{\uk}{})
(\tilde{z}{-}\zk)\,\d S = \Phi_\tau
(\uk,\tilde{z})+\calR(\JUMP{\uk}{}, \tilde{z}{-}\zk)
\end{align}
%\end{aligned}\]
for all $\tilde{z}\leq \zk$ a.e. on $\GC$. This  is the discrete
version of \eqref{disc-semistability}. $\hfill\Box$
%%%%%%%%%%%%%%%%%%%%%%%%%
%%%%%%%%%%%%%%%%%%%%%%%%%
%%%%%%%%%%%%%%%%%%%%%%%%%

\medskip

%\noindent

We conclude this section with a result collecting all the
a priori estimates on the approximate solutions.

\begin{lemma}[A priori estimates]\label{prop:apriori}
Under the  assumptions of Theorem~\ref{th:3.0},
%there exist constants $S_0>0$ and  $S_r>0$ such that\BT, \ET
for all $\varrho\geq 0$ and $\tau>0$, the
%for all
approximate solutions $(\pwc u {{\eps}\tau}, \pwc \w {{\eps}\tau},
\pwc z {{\eps}\tau}, \pwl u {{\eps}\tau}, \pwl \w {{\eps}\tau}, \pwl
z {{\eps}\tau})$  satisfy
%the following estimates hold
\begin{subequations}\label{a-priori3}
\begin{align}
& \label{a30} \big\|\pwc
u{{\eps}\tau}\big\|_{L^{\infty}(0,T;W_{\Gdir}^{2,2}(\Omega;\R^3))}\le
S_0\,,
\\
& \label{a31} \big\|\pwl u{{\eps}\tau}\big\|_{
W^{1,2}(0,T;W_{\Gdir}^{2,2}(\Omega;\R^3))}\le S_0\,,
\\
& \label{a31bis} \varrho^{1/2} \big\|\pwl u{{\eps}\tau}\big\|_{
W^{1,\infty}(0,T;L^2(\Omega;\R^3)) }\le S_0\,,
\\
&  \label{a35} \big\|\pwc
u{{\eps}\tau}\big\|_{L^\infty(0,T;W_{\Gdir}^{2,\gamma}(\Omega;\R^3))}\le
\frac{S_0}{\sqrt[\gamma]{\tau}},
\\&
\label{a32} \big\|\pwc z{{\eps}\tau}\big\|_{L^\infty(\SC)}
%,\, \big\|\pwl z{{\eps}\tau}\big\|_{L^\infty(\SC)}
\leq S_0\,,
\\
 & \label{a32bis}
 \big\|\pwc z{{\eps}\tau}\big\|_{BV([0,T];L^1(\GC))} \leq S_0,
\\&
\label{a33} \big\|\pwc
\w{{\eps}\tau}\big\|_{L^\infty(0,T;L^1(\Omega))}
%,\, \big\|\pwl \w{{\eps}\tau}\big\|_{ L^\infty(0,T;L^1(\Omega))}
\le S_0,
\\
 & \label{a33bis}
%\big\|\pwc \w{{\eps}\tau}\big\|_{L^r(0,T;W^{1,r}(\Omega))}, \,
\big\|\pwc \w{{\eps}\tau}\big\|_{ L^r(0,T;W^{1,r}(\Omega))} \leq
S_{{r}} \ \ \mbox{  for any $1\le r<\frac5 4$},
%,%
\\&
\label{a34}\big\|\pwl{\DT\w}{{\eps}\tau}\big\|_{L^1(0,T;W^{1,r'}(\Omega)^*)}\le S_0,
\\&
\label{a36} \varrho \,
\big\|\pwl{\DT{u}}{{\eps}\tau}\big\|_{BV([0,T];W_{\Gdir}^{2,\gamma}(\Omega;\R^3)^*)}\le
S_0\,,
\end{align}
\end{subequations}
 for some constants $S_0>0$ and  $S_r>0$
%where neither $S_0$ nor $S_r$  depend on
independent of  $\tau$. Estimates (\ref{a-priori3}e, f) hold for
$\pwl z{{\eps}\tau}$ as well, and so do estimates
(\ref{a-priori3}g,h) for $\pwl \w{{\eps}\tau}$.
\end{lemma}
%%%%%%%%%%%%%%%%%%
%%%%%%%%%%%%%%%%%%%%
\noindent{\it Proof.} We only sketch the calculations for proving
\eqref{a-priori3}, since the argument  closely follows the proof of
\cite[Lemma 7.7]{rr+tr}, to which  we shall systematically refer.

First of all, we use the ``discrete total energy''
balance~\eqref{disc-energy}. Clearly, the first summand on the
left-hand side provides a bound for the quantity $\varrho^{1/2}
\big\|\pwl u{{\eps}\tau}\big\|_{ W^{1,\infty}(0,T;L^2(\Omega;\R^3))
}$. Secondly,
 we observe that (cf. \eqref{8-1-k}, \eqref{Phi-eps-tau}),
%\begin{equation}\label{coerPhi}\begin{aligned}
%\begin{equation}\label{coerPhi}\begin{array}{ll}
\begin{align}\nonumber
\Phi_\tau (u,z) & \geq
  C \Big(\int_{\Omega}|e(u)|^2+ |\nabla e(u)|^2
+\tau  |e(u)|^\gamma+\tau  |\nabla e(u)|^\gamma\, \d x
  \\\nonumber
  & \qquad \qquad \qquad \qquad \qquad \qquad\quad
+ \int_{\GC}\!\!\! z \big|\JUMP{u}{}\big|^2 \, \d S\Big)- C_{a_0} \int_{\GC}\!\!\!  z \big|\JUMP{u}{}\big|\, \d S -C_{a_0}'
  \\
  & \geq  C \Big( \| u\|_{W^{1,2}(\Omega;\R^3)}^2+ \| u\|_{W^{2,2}(\Omega;\R^3)}^2
  + \tau\| u\|_{W^{2,\gamma}(\Omega;\R^3)}^\gamma + \int_{\GC}\!\!\! z
  \big|\JUMP{u}{}\big|^2 \, \d S
  \Big) -C',
\end{align}
%  \end{array}\end{equation}
where we have used the positive-definiteness of $\bbA$, $\bbC$, and
$\bbH$, and the growth condition \eqref{anot0}, to derive the first
inequality.  The second estimate ensues from Korn's inequality,
 and
from absorbing  the term $\int_{\GC}  z |\JUMP{u}{}|\,\d x$ into
$\int_{\GC} z |\JUMP{u}{}|^2\,\d x$, since $z \in [0,1]$ a.e.\ on
$\GC$. Therefore, the second term on the left-hand side
of~\eqref{disc-energy} estimates  $\|\pwc u {{\eps}\tau}(t)
\|_{W^{2,2}(\Omega;\R^3)}^2$ and  $\tau\| \pwc u {{\eps}\tau}(t)
\|_{W^{2,\gamma}(\Omega;\R^3)}^\gamma $ uniformly w.r.t.\ $t \in
[0,T]$. Thirdly,
 $\pwc\w{\tau} \geq
0$ a.e. in $\Omega$ thanks to Lemma \ref{lem-exist}, hence the third
term estimates $\| \pwc\w{\tau} \|_{L^\infty(0,T;L^1(\Omega))}$.
  To deal with the
right-hand side of~\eqref{disc-energy}, we use \eqref{hyp-init},
\eqref{est-init-data} and, for the last integral term,
\eqref{data-converg} and \eqref{shall-use-later},
 arguing  in the very same
way as in the proof of \cite[Lemma 7.7]{rr+tr}. We conclude applying
the discrete Gronwall lemma, and thus obtain estimates \eqref{a30},
 \eqref{a31bis}, and \eqref{a33}. Since  $\pwc z\tau  \in [0,1]$ a.e. on
$\SC$, we obviously have \eqref{a32}.

Secondly, again  arguing as for~\cite[Prop.~4.2]{tr1} and
\cite[Lemma 7.7]{rr+tr}, we make use
 of the technique by %\textsc{Boccardo \& Gallou\"et}
Boccardo and Gallou\"et~\cite{boccardo-gallouet1} with the
simplification devised in~\cite{feireisl-malek}, and we test the
heat equation~\eqref{weak-heat-discr} by $\pi (\pwc {\w}
{{\eps}\tau})$, where $\pi : [0,+\infty) \to [0,1]$ is the map $w
\mapsto \pi (w) = 1-\frac1{(1{+}w)^\varsigma}$, for $\varsigma>0$.
Since  $\pi$ is Lipschitz continuous,
 $\pi (\pwc {\w} {{\eps}\tau})\in W^{1,2}(\Omega{\setminus}
\GC )$ is an admissible test function. We thus have
\begin{align}\nonumber
\varsigma\, \mathsf{k}\! \int_Q \frac{|\nabla \pwc {\w}
{\tau}|^2}{(1{+}\pwc {\w} {\tau})^{1+\varsigma}}\,\d x\d t &\le
\int_Q\mathcal{K}(e(\pwc {u} {\tau}),\pwc {\w} {\tau}) \nabla\pwc
{\w} {\tau}{\cdot}\nabla\pi(\pwc {\w} {\tau})\,\d x\d t
\\\nonumber
&\quad+ \int_{\SC}\!\!\eta(\JUMP{{\underline u_{\tau}}}{}, \pwc
z{\tau}) \JUMP{\Theta(\pwc\w{\tau})}{}\JUMP{\pi(\pwc {\w}
{\tau})}{}\,\d S\d t
 +\int_\Omega\widehat{\pi}(\pwc {\w} {\tau}(T,{\cdot}))\,\d x
\\\nonumber &
 \le\int_\Omega\!\widehat{\pi}(\w_0)\,\d x
 +C\Big( \| \bbD e(\pwl{\DT{u}}{\tau})  {:} e(\pwl{\DT{u}}{\tau})
\|_{L^1 (Q)}\! + \| \bbG \nabla e(\pwl{\DT{u}}{\tau}){\vdots}\nabla
e(\pwl{\DT{u}}{\tau})\|_{L^1 (Q)}
\\ &\label{est:2-mixity}
%\begin{aligned}
\quad + \|\Theta(\pwc\w{\tau})\bbB{:}e(\pwl{\DT{u}}{\tau}) \|_{L^1
(Q)}\!+ \|  \zeta_1 (\JUMP{\pwc{u}{\tau}}{},\pwl
{\DT{z}}{{\eps}\tau}) \|_{L^1 (\SC)} \Big) + \| \pwc\gRM{\tau}
\|_{L^1 (\Sigma))},
\end{align}
where $\widehat{\pi}$ is the primitive function of $\pi$ such that
$\widehat{\pi}(0)=0$. Note that inequality~\eqref{est:2-mixity}
follows from \eqref{30dprimo},  the  fact that $\eta({\underline
u_{\tau},}\pwc z{\tau}) \JUMP{\Theta(\pwc\w{\tau})}{}\JUMP{\pi(\pwc
{\w} {\tau})}{} \geq 0$ a.e. in $\SC$ (by the positivity of $\eta$
and the monotonicity of $\Theta$ and $\pi$), from the ``discrete
chain rule'' \cite[Formula~(4.30)]{tr1} for $\widehat{\pi}$, and
from the fact $0 \leq \pi(\pwc {\w} {\tau}) \leq 1$ a.e.\ in
$\Omega$. Combining \eqref{est:2-mixity} and performing
 the very  same
calculations as in the proof of \cite[Prop.~4.2]{tr1}, with the
Gagliardo-Nirenberg inequality
 we find for all $1\leq r <5/4$, that
\begin{align}
\nonumber \big\|\nabla \pwc {\w} {{\eps}\tau} \big\|_{L^r(Q;\R^3)}^r
&\le C_r\Big(1+\|\bbD
e(\pwl{\DT{u}}{{\eps}\tau}){:}e(\pwl{\DT{u}}{{\eps}\tau}) \|_{L^1
(Q)}+\| \bbG \nabla e(\pwl{\DT{u}}{{\eps}\tau}){\vdots}\nabla
e(\pwl{\DT{u}}{{\eps}\tau}) \|_{L^1 (Q)}
 \\&
\ \ \ +
\|\Theta(\pwc\w{{\eps}\tau})\bbB{:}e(\pwl{\DT{u}}{{\eps}\tau})
\|_{L^1 (Q)}\!+ \|  \zeta_1 (\JUMP{\pwc{u}{\tau}}{},\pwl
{\DT{z}}{{\eps}\tau}) \|_{L^1 (\SC)}\Big) \label{est:3}
\end{align}
for some positive constant $C_r$,  depending on $r$ and also on  the
function $\eta$, cf.\ \eqref{eta-affine}.

Then, we multiply \eqref{est:3} by a  constant  $\rho_1>0$ and add
it to~\eqref{disc-energy0} (in which we set $t=T$). Now, by
positive-definiteness of $\bbD$ and $\bbG$, the third term on the
left-hand side of~\eqref{disc-energy0}
 is bounded from below by
 $c (\|e(\pwl{\DT{u}}{{\eps}\tau}) \|_{L^2(Q;\R^{3\times 3})}^2 +
 \|\nabla e(\pwl{\DT{u}}{{\eps}\tau}) \|_{L^2(Q;\R^{3\times 3\times 3})}^2 )$,
 and it
 controls
$\|\zeta_1(\JUMP{\upwc{u}{\tau}}{},\pwl{\DT{z}}{{\eps}\tau})\|_{L^1
(\SC)}$. Thus, we choose $\rho_1$ small enough in such a way as  to
absorb the second,  the third, and the fifth term on the right-hand
side of~\eqref{est:3} into the left-hand side
of~\eqref{disc-energy0}. Hence, we find
\begin{align}
&c \left( \|e(\pwl{\DT{u}}{{\eps}\tau}) \|_{L^2 (Q; \R^{3\times
3})}^2 + \|\nabla e(\pwl{\DT{u}}{{\eps}\tau}) \|_{L^2 (Q;
\R^{3\times 3 \times 3})}^2 \right) + (1{-}\rho_1) \| \zeta_1
(\JUMP{\pwc{u}{\tau}}{},\pwl {\DT{z}}{{\eps}\tau}) \|_{L^1 (\SC)}
+\rho_1 \big\|\nabla \pwc {\w} {{\eps}\tau} \big\|_{L^r(Q;\R^3)}^r
\nonumber
\\
&\leq   T_\mathrm{kin}\big(\DT{u}_{0,\tau})
+\Phi_{{\eps}\tau}\big(u_{0,\tau},z_{0,\tau})
%\nonumber\\ &
+
%\int_0^T\!\!\!\int_{\Omega}
 \!\int_Q\! \pwc \FRM{\tau}{\cdot} \pwl{\DT{u}}{{\eps}\tau} \,
\d x \, \d t
%\nonumber\\ &
+%\int_0^T\!\!\!\int_{\Gnew}
 \!\int_{\Snew}\! \! \pwc \fRM{\tau} {\cdot}
\pwl{\DT{u}}{{\eps}\tau} \,  \d S \d t
 + (\rho_5 C_r{+}1)
\big\|\Theta(\pwc\w{{\eps}\tau})\bbB e(\pwl{\DT{u}}{{\eps}\tau})
\big\|_{L^1(Q)}. \label{disc-energy01}
\end{align}
The first two summands on the right-hand side
of~\eqref{disc-energy01} are estimated in view of \eqref{hyp-init}
and~\eqref{est-init-data}. Using \eqref{data-converg},  we  handle
the terms $ \int_0^T\int_{\Omega}\pwc \FRM{\tau}{\cdot}
\pwl{\DT{u}}{{\eps}\tau}\, \d x\d t$ and $ \int_0^T\int_{\Gnew} \pwc
\fRM{\tau} {\cdot} \pwl{\DT{u}}{{\eps}\tau} \, \d S\d t$ in the very
same way as in the proof of \cite[Lemma 7.7]{rr+tr}. Finally, we use
\begin{align}
(\rho_1 C_r{+}1)\|\Theta(\pwc\w{{\eps}\tau})\bbB{:}
e(\pwl{\DT{u}}{{\eps}\tau}) \|_{L^1(Q) } & \leq \rho_2 \,\|
e(\pwl{\DT{u}}{{\eps}\tau}) \|_{L^2(Q; \R^{3\times 3}))}^2 +
C_{\rho_2} \|\Theta(\pwc\w{{\eps}\tau}) \|_{L^2(Q)}^2 \nonumber\\ &
\leq \rho_2 \,\|e(\pwl{\DT{u}}{{\eps}\tau})\|_{L^2(Q; \R^{3 \times
3}))}^2 + C_{\rho_2}\big(\|\pwc\w{{\eps}\tau}
\|_{L^{2/\omega}(Q)}^{2/\omega}\!+ 1\big) \nonumber
\\
& \le\rho_2\,\|e(\pwl{\DT{u}}{{\eps}\tau})\|_{L^2(Q; \R^{3 \times 3}))}^2 +
 \rho_3 \int_0^T\!
\big\|\nabla \pwc {\w} {{\eps}\tau} \big\|_{L^r(\Omega;\R^3)}^r \,
\d t+ C_{\rho_3}, \label{est:4}
\end{align}
where the last inequality can be proved, via the Gagliardo-Nirenberg inequality,
by developing the
same calculations as throughout~\cite[Formulae~(4.39)--(4.43)]{tr1},
and using the restriction on $\omega$ in \eqref{30b} and  the
previously proved bound for $\| \pwc \w {{\eps}\tau}\|_{L^\infty
(0,T;L^1 (\Omega))}$.
 Then, we plug~\eqref{est:4}  into
~\eqref{disc-energy01}, and choose $\rho_2$ and $\rho_3$ in such a
way as to absorb the terms $\|e(\pwl{\DT{u}}{{\eps}\tau})\|_{L^2(Q;
\R^{3 \times 3})}^2$ and $ \big\|\nabla \pwc {\w} {{\eps}\tau}
\big\|_{L^r(Q;\R^3)}^r$
 into
the left-hand side of \eqref{disc-energy01}. Thus, we conclude
estimate~\eqref{a31}. We also get
 an estimate for
$\|\zeta_1(\JUMP{\pwc
u\tau}{},\pwl{\DT{z}}{{\eps}\tau})\|_{L^1(\SC)}$, which yields
\eqref{a32bis}, since $a_1$ is bounded from below, cf.\
\eqref{anot1}. Furthermore, we also obtain  a bound for $\nabla \pwc
{\w} {{\eps}\tau}$ in $L^r(Q;\R^3)$. Combining the latter
information with the estimate for $\pwc {\w} {{\eps}\tau}$ in
$L^\infty (0,T; L^1(\Omega))$, we infer~\eqref{a33bis}.

Estimate  \eqref{a34} follows from a comparison
in~\eqref{weak-heat-discr}, and the related calculations are a
trivial adaptation of the ones in the proof of \cite[Lemma
7.7]{rr+tr}.

Finally, for~\eqref{a36} we use that $\pwl{\DDT{u}}{{\eps}\tau}$ is
a measure on $[0,T]$, supported at the jumps of
$\pwl{\DT{u}}{{\eps}\tau}$, and we estimate
 $\varrho \|\pwl{\DDT{u}}{{\eps}\tau}\|_{\mathcal{M}(0,T;W_{\Gdir}^{2,\gamma}(\Omega{{{\setminus}}}\GC;\R^3)^*)}$
 %(where $\mathrm{M}(0,T;W_{\Gdir}^{2,\gamma}(\Omega{{{\setminus}}}\GC;\R^3)^*)$
%denotes the space of Radon measures on $[0,T]$ with values in
% $W_{\Gdir}^{2,\gamma}(\Omega{{{\setminus}}}\GC;\R^3)^*$),
by comparison in~\eqref{e:discrmom}.
%\begin{equation}
%\label{e:desired}
%\end{equation}
 $\hfill\Box$
%%%%%%%%%%%%%%%%%%%%%%%%%
%\section{A priori estimates}\label{ss:6.1}
%        ~~~~~~~~~~~~~~~~~~

%%%%%%%%%%%%%%%%%%%%%%%%%%
%%%%%%%%%%%%%%%%%%%%%%%%%%%%%%%%%%%%

%%%%%%%%%%%%%%%%
%%%%%%%%%%%%%%%%%%%%%%%%%%%%

\section{Proof of Theorem~\ref{th:3.0}}\label{ss:4.4}
%      ~~~~~~~~~~~~~~~~~~~~~~~~~~~~~~~~~~~~~~~~~~~~~~~~~~~~~~~~

\medskip

\noindent  In what follows, we develop the proof of the passage to
the limit in the time-discrete scheme as $\tau \to 0$
%(we shall understand $(\tau)$ as a family of countably many indexes,
%hence we will speak of a (vanishing) \emph{sequence} of time-steps),
unifying  the cases $\varrho>0$ and $\varrho=0$;  we shall take a
\emph{sequence} of time-steps, i.e.\ we understand $(\tau)$ as
countable family of indexes with the accumulation point $0$.

%\DELETE{\par\noindent \textbf{Scheme of the passage to the limit proof.} }
For the  reader's convenience, we  briefly describe the strategy.
After a careful selection of converging subsequences is made in Step
0 below,  by passing to the limit  as $\tau \to 0$ in
\eqref{e:discrmom} we will obtain the weak formulation
\eqref{e:weak-momentum-variational} of the momentum inclusion. Then,
we will proceed to proving the semistability condition
\eqref{semistab}, hence the total energy \emph{inequality} by lower
semicontinuity arguments. By the same tokens we will also obtain the
mechanical energy inequality. We will then show that the latter in
fact holds as an equality, by combining a \emph{chain rule-type}
argument (cf.\ \eqref{e:faith-delam}), with a test of
\eqref{e:weak-momentum-variational} by $\DT u$. To perform the
latter, it will be essential for $\DT u$ to have the regularity
\eqref{relaxed-regu}. This motivates  the dissipative contribution
$\bbG \nabla e(\DT u)$ to the hyperstress.
 Hence, we will
exploit the mechanical energy equality to conclude, via a suitable
comparison argument, the convergence of the quadratic terms in the
right-hand side of \eqref{weak-heat-discr}. This will allow us to
pass to the limit, and conclude the weak formulation
\eqref{weak-heat} of the enthalpy equation and, ultimately, also the
total energy balance \eqref{total-energy-brittle}. %%%%%%%%%%%%%%%%%%%%%

\medskip
\par\noindent {\textbf{Step $0$: selection of convergent subsequences.}}
First of all, it follows from estimates~\eqref{a31}, \eqref{a31bis},
and \eqref{a36}, from the Banach selection principle, the
infinite-dimensional Ascoli and the  Aubin-Lions theorems  (see,
e.g., \cite[Thm.~5, Cor.~4]{simon86}), that there exist a (not
relabeled) sequence $\tau \to 0$
 and a  limit function
$\ue \in W^{1,2}(0,T;W_{\Gdir}^{2,2}(\Omega{{\setminus}} \GC;\R^3))$
such that the following    convergences hold as $\tau\to 0$:
\begin{subequations}
\label{e:convutau}
\begin{align}
 \label{e:convutau1} & \pwl {u}{{\eps}\tau} \weakto
\ue\ \ \text{ in $W^{1,2}(0,T;W_{\Gdir}^{2,2}(\Omega{{\setminus}}
\GC;\R^3)),$}
\\
 \label{e:convutau2}
 & \pwl {u}{{\eps}\tau} \to
\ue\ \ \text{ in
$C([0,T];W_{\Gdir}^{2{-}\epsilon,2}(\Omega{{\setminus}} \GC;\R^3))$}
\ \ \forall\, \epsilon \in
 (0,2],
 \\
\label{e:convutau3} &   \varrho\pwl {u}{{\eps}\tau} \weaksto
\varrho\ue \ \ \text{ in $ W^{1,\infty}(0,T;L^2(\Omega;\R^3))$.}
\end{align}
 Estimate~\eqref{a36}  and a generalization of the  Aubin-Lions theorem
to the case of time derivatives as measures (cf.\ e.g.\
\cite[Cor.~7.9]{NPDE_roubicek}) also yield that $\DT{u}_{\eps} \in
BV ([0,T];W_{\Gdir}^{2,\gamma}(\Omega{{\setminus}} \GC;\R^3)^*)$ and
that
\begin{equation}
\label{strong-for-dot-u} \varrho\pwl {\DT{u}}{\eps \tau} \to
\varrho\DT{u}_{\eps} \quad \text{in
$L^2(0,T;W_{\Gdir}^{2{-}\epsilon,2}(\Omega{{\setminus}} \GC;\R^3))$
for all $\epsilon \in
 (0,2]$.}
\end{equation}
Moreover,
  a
generalization of Helly's principle (see \cite{BarbuPrecupanu86} as
well as \cite[Thm.~6.1]{MieThe04RIHM})  implies that $\pwl
{\DT{u}}{\eps \tau}(t) \weakto \DT{u}_{\eps}(t)$ in
$W_{\Gdir}^{2,\gamma}(\Omega{{\setminus}} \GC;\R^3)^*$ for all $t
\in [0,T]$. In view  of estimate~\eqref{a31bis},  with an elementary
compactness argument we conclude
\begin{equation}
\label{pointiwise-for-u} \varrho \pwl {\DT{u}}{\eps \tau}(t) \weakto
\varrho \DT{u}_{\eps}(t) \quad \text{in $L^2(\Omega;\R^3)$ for all
$t \in [0,T]$.
%, in the case $\varrho>0$.
}
\end{equation}
Next, we observe that
%\begin{equation}\label{e:stability}\begin{aligned}
\begin{align}
%\nonumber
&
%\max\Big(
\|\pwc u\tau{-}\pwl u\tau \|_{L^\infty
(0,T;W_{\Gdir}^{2,2}(\Omega{{\setminus}} \GC;\R^3)) }
%,\, \| \upwc u\tau{-}\pwl u\tau \|_{L^\infty
%(0,T;W_{\Gdir}^{2,2}(\Omega{{\setminus}} \GC;\R^3)) }\Big) \\
 \label{e:stability}
%&\qquad
\le\tau^{1/2}\|{\DT u_\tau} \|_{L^2
(0,T;W_{\Gdir}^{2,2}(\Omega{{\setminus}} \GC;\R^3)) } \leq
S_0\tau^{1/2} \to 0 \quad \text{as $\tau \to 0$}.
\end{align}
%\end{aligned}\end{equation}
Therefore, estimate \eqref{a30} and
\eqref{e:convutau1}--\eqref{e:convutau2} yield for all $\epsilon
\in(0,1]$
%\begin{equation}\label{e:convutau5}\begin{gathered}
\begin{align}\label{e:convutau5}
& \pwc {u}{{\eps}\tau} \weaksto u_{\eps} \ \ \text{ in
 $L^{\infty}(0,T;W_{\Gdir}^{2,2}(\Omega{{\setminus}} \GC;\R^3))$,}
 \quad
 \\\label{e:convutau5+}
&\pwc {u}{{\eps}\tau}
 \to u_{\eps}  \ \ \text{ in
$L^{\infty}(0,T;W_{\Gdir}^{2{-}\epsilon,2}(\Omega{{\setminus}}
\GC;\R^3))$,}
 \\\label{e:convutau5++}
&\pwc {u}{{\eps}\tau}(t)
 \to \ue(t) \ \
\text{in $W_{\Gdir}^{2{-}\epsilon,2}(\Omega{{\setminus}}\GC;\R^3))$
for all $t \in [0,T]$.}
\end{align}
%\end{gathered}\end{equation}
Taking into account the compact embedding \eqref{e:contemb}, from
\eqref{e:convutau2} and \eqref{e:convutau5+} we deduce respectively
\begin{align}
\label{conv2} & \JUMP{\pwl u{{\eps}\tau}}{} \to \JUMP{u_{\eps}}{} &&
\text{in $C(\overlineSC;\R^3)$},
\\&\label{conv2+}
  \JUMP{\pwc u{{\eps}\tau}}{}
  \to \JUMP{u_{\eps}}{} && \text{in
$L^\infty(\SC;\R^3)$.}
\end{align}
%\CHECK{Note that $L^\infty
%(0,T; L^{\infty} (\GC;\R^3))\ne L^\infty(\SC;\R^3)$!!!!!!!}
In fact, the above convergences  are also  in $C([0,T];
C(\overlineGC;\R^3))$ and in $L^\infty(0,T; L^{\infty} (\GC;\R^3))$,
respectively.
 Convergences (\ref{e:convutau}g-i,k)  hold for $\upwc
u{{\eps}\tau}$, too.   Also note that, in view of \eqref{a35}, we
have
\begin{align}
\label{conv1} \left.\begin{array}{l} \tau\big\| |e(\pwc u
{{\eps}\tau})|^{\gamma-2}e(\pwc u
{{\eps}\tau})\big\|_{L^{\gamma/{(\gamma-1)}}(Q; \R^{3\times 3})} \le
S_0{\tau^{1/\gamma}} \to 0\ \text{ and}\\
\tau\big\| |\nabla e(\pwc u {{\eps}\tau})|^{\gamma-2} \nabla e(\pwc
u {{\eps}\tau})\big\|_{L^{\gamma/{(\gamma-1)}}(Q; \R^{3\times
3\times 3})} \le S_0{\tau^{1/\gamma}} \to 0
\end{array}\right\}\quad \text{ as $\ \tau\to0$.}
\end{align}
\end{subequations}

Estimates \eqref{a32} and \eqref{a32bis}, and the very same
compactness arguments as in \cite[Sec. 8]{rr+tr} (based on
\cite[Thm. 6.1, Prop. 6.2]{MieThe04RIHM}), also guarantee
 that there exists a
function $\ze \in L^\infty(\SC) \cap BV ([0,T]; \calZ) $ (where
$\calZ$ is any reflexive space such that $L^1(\GC) \subset \calZ$
with a continuous embedding),   such that, possibly along a
subsequence,
\begin{align}
\label{e:convztau1}
    \pwc {z}{{\eps}\tau}
\weaksto z_{\eps}  \ \ \text{ in $L^{\infty}(\SC)$,} \qquad  \pwc
{z}{{\eps}\tau}(t)\weaksto \ze(t) \ \ \text{ in $L^\infty (\GC)$ for
all $t \in [0,T]$.}
\end{align}
We now prove that
\begin{equation}
\label{e:finally}
\mathrm{Var}_{\mathcal{R}}(u_{\eps},z_{\eps};[s,t])\le \liminf_{\tau
\to 0} \int_s^t\!\int_{\GC}\!\zeta_1
\big(\JUMP{\pwc{u}{\tau}}{},\pwl{\DT{z}}{{\eps}\tau}\big)\,\d S\d r
\quad  \text{for all $0 \leq s \leq t \leq T$.}
\end{equation}
Indeed,
\begin{align}\nonumber
&\int_s^t\!\int_{\GC}\!\zeta_1
\big(\JUMP{\pwc{u}{\tau}}{},\pwl{\DT{z}}{{\eps}\tau}\big)\,\d S\d r
=\int_s^t\!\int_{\GC}\!\!
a_1(\JUMP{\pwc{u}{\tau}}{})|\pwl{\DT{z}}{{\eps}\tau}|\,\d S\d r
\\&\hspace{4em}
=\int_s^t\!\int_{\GC}\!\!a_1(\JUMP{u_\tau}{})|\pwl{\DT{z}}{{\eps}\tau}|\,\d
S\d r
+\int_s^t\!\int_{\GC}\!\!(a_1(\JUMP{\pwc{u}{\tau}}{})-a_1(\JUMP{u_\tau}{}))
|\pwl{\DT{z}}{{\eps}\tau}|\,\d S\d r. \label{7.3}
\end{align}
Now, from  (\ref{e:convutau}j,k) and the continuity of $a_1$ it
follows
\begin{equation}
\label{basic-1}
 a_1\big(\JUMP{\pwc{u}{\tau}}{}) - a_1(\JUMP{\pwl u{\tau}}{})\to 0
 \quad \text{
in $L^\infty (0,T; C(\overlineGC))$.}
\end{equation}
Since $(\pwl{\DT{z}}{{\eps}\tau})_{\tau>0}$ is bounded in
$L^1(\SC)$, we then conclude that the  second term on the right-hand
side of \eqref{7.3} tends to zero as $\tau \to 0$. To pass to the
limit in the first term, we  use  \eqref{conv2} and again the
continuity of $a_1$.
 Since
 $|\pwl{\DT{z}}{{\eps}\tau}|\to |\pwl{\DT{z}}{}|$
weakly* in the sense of  measures on $\overlineSC$, we conclude that
\begin{equation}
\label{measure-convergence}
a_1(\JUMP{u_\tau}{})|\pwl{\DT{z}}{{\eps}\tau}| \to
a_1(\JUMP{u}{})|\pwl{\DT{z}}{}|  \quad \text{ weakly* in the sense
of measures on $\overlineSC$,}
\end{equation}
 where $|\DT z|$ denotes the
variation of the measure $\DT z$; in fact, here simply $|\DT z|=-\DT
z$, since $\DT z\le0$. Then, \eqref{e:finally} follows. Taking into
account that $a_1$ is bounded from below by \eqref{anot1},
\eqref{e:finally} and the definition \eqref{var-R} of
$\mathrm{Var}_{\mathcal{R}}$ imply
 $\ze\in BV([0,T];L^1 (\GC))$.

With the same compactness tools
 as in the above lines,
we conclude from estimates \eqref{a33}, \eqref{a33bis}, and
\eqref{a34} that there exists $\we \in
L^{r}(0,T;W^{1,r}(\Omega{{\setminus}} \GC)) \cap BV ([0,T];
W^{1,r'}(\Omega{{\setminus}} \GC)^*)$ s.t.
\begin{subequations}
\label{e:convwtau}
\begin{align}
 \label{e:convwtau1} & \pwc {\w}{{\eps}\tau}\rightharpoonup \w_{\eps}
\ \ \text{ and }\ \ \pwl {\w}{{\eps}\tau} \rightharpoonup \w_{\eps}
\ \ \text{ in $ L^{r}(0,T;W^{1,r}(\Omega{{\setminus}} \GC))$ for
$1\le r<\frac54$,}
\\
\label{e:convwtau1-bis-serve} &    \pwc {\w}{{\eps}\tau}, \,  \pwl
{\w}{{\eps}\tau} \to \w_{\eps} \ \ \text{ in $
L^{r}(0,T;W^{1-\epsilon,r}(\Omega{{\setminus}} \GC))$   for all
$\epsilon \in (0,1] $}. %\cap L^q (0,T;L^1(\Omega))$  and $1\leq q
%<\infty$
\end{align}
%%%%%%%%%%%%
The latter convergence yields that $\pwc {\w}{{\eps}\tau}, \,  \pwl
{\w}{{\eps}\tau} \to \w_{\eps}$ in $L^r (0,T;
L^{15/7-\epsilon}(\Omega))$ for all $\epsilon \in (0,8/7]$. Taking
into account estimate \eqref{a33} in $L^\infty (0,T;L^1(\Omega))$
and arguing by interpolation, we conclude
\begin{align}
\label{e:new-label} \pwc {\w}{{\eps}\tau}, \,  \pwl {\w}{{\eps}\tau}
\to \w_{\eps} \qquad \text{ in $L^{5/3-\epsilon}(Q)\ $ for all }\
\epsilon \in \big(0,\mbox{$\frac23$}\big].
\end{align}
%\DELETE{for all $\epsilon \in
% (0,r-1]$ and $1 \leq q <\infty$.}
%Convergence \eqref{e:new-label} follows from \eqref{e:convwtau1} and
%an interpolation by exploiting \eqref{a33} and \eqref{a34}.
 Notice
that, under condition \eqref{strict-pos} on $\theta_0$, \eqref{posw}
and  \eqref{e:new-label} imply the strict positivity of $\w$.
 Moreover, Helly's selection principle and  the a priori bound for $(\pwc\w\tau)_\tau$
 in $L^\infty (0,T;L^1(\Omega))$ yield
 \begin{equation}
\label{e:poinwtiwise-w+} \pwl{\w}{{\eps}\tau}(t) \weaksto
\w_{\eps}(t) \quad \text{in $\mathcal{M}(\overline{\Omega})$ for all
$t \in [0,T]$.}
\end{equation}
\end{subequations}
%Furthermore, \eqref{e:convwtau2} yields by interpolation
%\begin{subequations}
% \label{pass-limi-w}
%\begin{equation}
%\label{e:new-label} \pwc {\w}{{\eps}\tau} \to \w_{\eps} \qquad \text{
%in $L^{15/7-\epsilon} (Q)\ $ for all }\ \epsilon \in
%\big(0,\mbox{$\frac87$}\big].
%\end{equation}
 Combining \eqref{e:new-label}
%this information
with~\eqref{growthTheta},
 it is immediate to deduce
%from \eqref{e:new-label} that, for example,
\begin{equation}
\label{conv3} \Theta(\pwc
{\w}{{\eps}\tau})\to\Theta(\w_{\eps})\qquad \text{ in $L^{2} (Q)$}.
\end{equation}
Furthermore, it follows from \eqref{e:convwtau1-bis-serve}, the
trace theorem $\vartheta\mapsto\JUMP{\vartheta}{}:
W^{1-\epsilon,r}(\Omega{{\setminus}} \GC)\to L^{10/7-\epsilon}(\GC)$
for all  $\epsilon \in (0,3/7]$, and \eqref{growthTheta}, that
\begin{equation}
\label{e:interesting}
 \JUMP{\Theta(\pwc {\w}{{\eps}\tau})}{}  \to  \JUMP{\Theta(\w_{\eps})}{}\qquad
\text{ in $L^{r\omega}(0,T;L^{\omega (10/7-\epsilon)}(\GC)) $}  \ \
\forall\, \epsilon \in \big(0,\mbox{$\frac57$}\big]\,.
\end{equation}
%\end{subequations}
%In the end, we are now going to  show that
 Exploiting  \eqref{e:convutau5++}, \eqref{e:convztau1}, and
\eqref{conv2+}, which in particular yields 
\[
\lim_{\tau\to0}
\int_{\GC}\pwc{z}{{\eps}\tau}(t) \alpha_0(\JUMP{\pwc
u{{\eps}\tau}(t)}{})\, \d S =\int_{\GC}z(t)\alpha_0(\JUMP{u(t)}{})\,
\d S
\quad \text{
 for all $t \in [0,T]$,}
 \]
  we conclude by lower semicontinuity
arguments that
\begin{equation}
\label{e:altogether} \Phi_{\eps}(u_{\eps}(t), {z}_{\eps}(t)) \leq
\liminf_{\tau \to 0} \Phi_{{\eps}\tau}(\pwc u{{\eps}\tau}(t), \pwc
z{{\eps}\tau}(t)) \quad \text{for all $t \in [0,T]$}.
\end{equation}
%%%%%%%%%%%%%%%%%%
%%%%%%%%%%%%%%%%%%%%%
%\par\noindent
\textbf{Step $1$: passage to the limit in the momentum equation.}
%\noindent
At  first, we  take the limit as $\tau \to 0$ of the discrete
momentum equation~\eqref{e:discrmom} with  smooth  test functions
$\testu \in C^\infty (Q;\R^3)$, fulfilling $\JUMP{\testu}{}\GE0 $ on
$\SC$. We  approximate them with discrete approximations $\{
\testu_{\tau}^k \}$ such that $\JUMP{\testu_{\tau}^k}{}\GE 0$ on
$\SC$ and the related piecewise constant and linear interpolants
fulfill, as $\tau \to 0$,
\begin{equation}
\label{conve-testu} \left.
\begin{array}{ll}
\pwl \testu{\tau} \to \testu & \text{in $ W^{1,1}(0,T; L^2(\Omega;
\R^3))$,}
\\ \pwc \testu{\tau} \to \testu & \text{in $L^2
(0,T; W^{2,2}(\Omega{{\setminus}} \GC;\R^3))$,}
% \\ \pwc \testu{\tau} \weakto \testu & \text{in $L^2
%(0,T; W^{\INSERT{2},2{+}\upsilon}(\Omega{{\setminus}} \GC;\R^3))$ for some
%$\upsilon>0$,}
\\
\left\|\nabla  e(\pwc \testu{\tau})
\right\|_{L^{\gamma}(Q;\R^{3\times 3\times 3})} \leq C\,.
\end{array}
 \right\}
\end{equation}
In order to pass to the limit in the first integral term on the
left-hand side of \eqref{e:discrmom}, we use \eqref{e:convutau1},
\eqref{e:convutau2}, and \eqref{conv3}, combined with
\eqref{conve-testu}. The regularizing $\gamma$-terms in
\eqref{e:discrmom} vanish in the limit due to \eqref{conv1}.
 Notice that \eqref{conv2+} and the continuity of
$\alpha_0'$ (cf.\ \eqref{anot0}) imply $\alpha_0'(\JUMP{\pwc u
{{\eps}\tau}}{}) \to \alpha_0'(\JUMP{u}{})$ in $L^\infty(\SC)$.
Therefore, taking \eqref{e:convztau1} into account, we deduce that
 ${\upwc z {{\eps}\tau}} \alpha_0'(\JUMP{\pwc u
{{\eps}\tau}}{})  \weaksto z \alpha_0'(\JUMP{u}{})$ in
$L^\infty(\SC)$.
%These convergences
This convergence and \eqref{conve-testu} allow us to pass to the
limit in the second integral term on the left-hand side of
\eqref{e:discrmom}. To take the limit of the third and fourth terms
(in the case $\varrho>0$), we use \eqref{e:convutau3},
\eqref{strong-for-dot-u}, and \eqref{pointiwise-for-u}, as well as
the first of \eqref{conve-testu}. Combining the latter with
\eqref{est-init-data} we also take the limit of the first term on
the right-hand side of \eqref{e:discrmom}. The convergence of the
other two integrals ensues from \eqref{data-converg},
\eqref{e:convutau1}, and \eqref{conve-testu}.  Thus, we have proved
that the triple $(\ue,\ze,\we)$ fulfills
equation~\eqref{e:weak-momentum-variational}  with
 smooth test  functions. With a density
argument, we conclude \eqref{e:weak-momentum-variational} with test
functions $\testu \in L^2 (0,T;W_{\Gdir}^{2,2}(\Omega{{\setminus}}
\GC;\R^3)) \cap W^{1,1}(0,T;L^2 (\Omega;\R^3))$. For later
convenience, let us observe that, in the case $K(x)$ is a linear
subspace for almost all $x \in \GC$ (cf.\ \eqref{varrho2}), taking
test functions $\tilde{v}= u +\lambda \testu$ with $\testu$ any
admissible test function satisfying  $\JUMP{\testu}{}\GE 0$ on $\SC$
and $\lambda$ an arbitrary real number, we obtain
\eqref{e:weak-momentum-variational} in the form
\begin{align}
  &\int_\Omega\!\varrho
\DT{u}(T){\cdot}\testu(T)\,\d x +\int_Q\!\big(\bbD e(\DT{u}){+}\bbC
e(u){-}\bbB\Theta(\w)\big){:}e(\testu) + \big(\bbH\nabla e(u){+}\bbG
\nabla e(\DT{u})\big){\vdots}\nabla e(\testu)
-\varrho\DT{u}{\cdot}\DT{\testu}
 \,\d x\d t \nonumber
\\
 &\qquad\qquad
+\int_{\SC}\!\!\alpha_0'(\JUMP{u}{})
%\left( z\bbA\JUMP{u}{}{-} z a_0'(\JUMP{u}{})\right)
{\cdot}\JUMP{\testu}{}\d S\d t
=\int_\Omega\!\varrho\DT{u}_0{\cdot}\testu(0)\d x +
%\int_0^T\!\left(\int_{\Omega}
\int_Q \!\FRM{\cdot}\testu\, \d x\d t-
%\int_{\Gnew}
\int_{\Snew} \!\!\fRM{\cdot}\testu
%\right)
\,\d S \d t \label{e:weak-mom-identity}
\end{align}
for
%\emph{any}
any $\testu \in L^2(0,T;W_{\cone}^{2,2}(\Omega{{\setminus}} \GC;\R^3)) \cap
W^{1,1}(0,T;L^2 (\Omega;\R^3))$.
 %%%%%%%%%%%%%%%%
\par\noindent
\textbf{Step $2$: passage to the limit in the semistability
condition.}
%It follows from~\eqref{irrev} that $z_{\eps}$ complies
%with the irreversibility constraint~\eqref{irrev}. In order to prove
%that the functions $(u_{\eps},z_{\eps})$ comply with the semistability
%condition,
We consider  a subset $\mathcal{N} \subset (0,T)$ of full measure
such that  for all $t \in \mathcal{N}$  the approximate stability
condition \eqref{disc-semistability} holds for the (countably many)
considered
 $\tau$'s. Then we fix $t \in \mathcal{N}$ and
$\tilde{z} \in L^\infty (\GC)$. We may suppose without loss of
generality that
$\mathcal{R}(u_{\eps}(t),\tilde{z}{-}z_{\eps}(t))<+\infty$, hence
\begin{equation}
\label{e:use} \tilde{z}(x) \leq z_{\eps}(t,x)\qquad \foraa\,x \in
\GC\,.
\end{equation}
Then,  we construct the following recovery sequence
\begin{align}\label{eq5:recov-seq}
   \tilde z_{{\eps}\tau}(t,x):=\begin{cases}
   \displaystyle \pwc z{{\eps}\tau}(t,x)\frac{\tilde{z}(x)}{z_{\eps}(t,x)}
    & \text{where $z_{\eps}(t,x)>0$},
\\
   \displaystyle 0 & \text{where $z_{\eps}(t,x)=0$}.\end{cases}
\end{align}
Now, using~\eqref{e:use} and~\eqref{e:convztau1} one immediately
sees that
\begin{equation}
\label{e:use2} \tilde z_{{\eps}\tau}(t,\cdot) \leq \pwc
z{{\eps}\tau}(t,\cdot) \quad \aein \,\GC, \qquad \tilde
z_{{\eps}\tau}(t) \weaksto \tilde{z} \ \ \text{in $L^\infty (\GC)$.}
\end{equation}
Plugging $\tilde z_{{\eps}\tau}$ in~\eqref{disc-semistability} and
using \eqref{conv2+}, \eqref{e:convztau1}, and \eqref{e:use2}, we
find
\begin{align}
\nonumber 0 & \leq \lim_{\tau \to 0} \Big( \Phi_{{\eps}\tau}(\pwc
u{{\eps}\tau}(t), \tilde z_{{\eps}\tau}{(t)}){+}
\mathcal{R}(\pwc{u}{\tau}(t),\tilde z_{{\eps}\tau}{(t)} - \pwc
z{{\eps}\tau}(t) ){-}\Phi_{{\eps}\tau}(\pwc u{{\eps}\tau}(t), \pwc
z{{\eps}\tau}(t))\Big)
\\\nonumber&
\begin{aligned}=\lim_{\tau \to 0}\int_{\GC}\!\!
%\Big( \frac{1}{2}\left(\tilde z_{{\eps}\tau}\DELETE{(t)}{-}\pwc
%z{{\eps}\tau}(t)\right)\dela\JUMP{\pwc
%u{{\eps}\tau}(t)}{}{\cdot}\JUMP{\pwc u{{\eps}\tau}(t)}{}-
\big(\alpha_0(\JUMP{\pwc{u}{\tau}(t)}{})
{-}a_1(\JUMP{{\pwc{u}{\tau}(t)}}{})\big)\big(\tilde
z_{{\eps}\tau}{(t)} {-} \pwc z{{\eps}\tau}(t)\big) \, \d S
\end{aligned}
\nonumber
\\\nonumber&
\begin{aligned}
= \int_{\GC}\!\!
%\Big(  \frac{1}{2}\left( \tilde {z}(t){-}z_{\eps}(t)\right)
%\dela\JUMP{u_{\eps}(t)}{}{\cdot}\JUMP{u_{\eps}(t)}{}{-}
\big(\alpha_0(\JUMP{u(t)}{}){-}a_1(\JUMP{u(t)}{})\big) \big(\tilde
z_{{\eps}}{(t)} {-}\pwc z{{\eps}}(t)\big)\, \d S
\end{aligned}
\\\label{eq6:delam-recovery}
& = \Phi_{\eps}(u_{\eps}(t), \tilde {z}{(t)}){+}
\mathcal{R}(u_{\eps}(t),\tilde z{(t)}{-}z_{\eps}(t))
{-}\Phi_{\eps}(u_{\eps}(t), {z}_{\eps}(t))\,.
\end{align}
%\par\noindent
\textbf{Step $3$: passage to the limit in the  mechanical and total
energy inequalities.} Using \eqref{e:convutau1},
\eqref{pointiwise-for-u}, \eqref{e:finally}, and
\eqref{e:altogether},
 we pass to the limit on the left-hand side of the
discrete mechanical energy inequality \eqref{disc-energy0}  by weak
lower semicontinuity. To take the limit of the right-hand side, we
employ~\eqref{est-init-data} (which in particular yields
$\Phi_{{}\tau}\big(u_{0,\tau},z_{0}) \to \Phi(u_0,z_0)$), the weak
convergence \eqref{e:convutau1} and the strong
convergence~\eqref{conv3}, which give
\begin{align}
\label{conv-adiab} \Theta (\pwc \w{{\eps}\tau}) \bbB{:}
e(\pwl{\DT{u}}{{\eps}\tau}) \rightharpoonup \Theta (\w_{\eps}) \bbB
{:}e({\DT{u}}_{\eps}) \quad \text{weakly in $L^1 (Q)$.}
\end{align}
 We pass to the limit in the two remaining terms by
\eqref{data-converg-1}--\eqref{data-converg-2}. Hence, the triple
$(\ue,\ze,\we)$ complies  for all $t\in [0,T]$ with
\begin{align}\nonumber
&\! \! \! \! \! \! \! \! T_\mathrm{kin}\big(\DT{u}_{\eps}(t)\big)
+\Phi_{\eps}\big(\ue(t),\ze(t)\big) +\int_0^t \int_{\Omega} \bbD
e\big(\DT{u}_{\eps}(s)\big){:} e\big(\DT{u}_{\eps}(s)\big){+} \bbG
\nabla e\big(\DT{u}_{\eps}(s)\big){\vdots} \nabla
e\big(\DT{u}_{\eps}(s)\big)  \,\d x\d s
\\\nonumber
&\qquad + \mathrm{Var}_{\mathcal{R}}(u_{\eps},z_{\eps};[0,t]) \le
T_\mathrm{kin}\big(\DT{u}_{0})
 +\Phi_{\eps}\big(u_{0},z_{0})
%\nonumber\\& \hspace{.1em}
\\
&\qquad +\int_0^{t}\bigg(\int_\Omega \Theta(\we(s))\bbB
{:}e\big(\DT{u}_{\eps} (s) \big) \,\d x + \int_{\Omega}\FRM(s)
{\cdot}\DT{u}_{\eps}(s)\, \d x   +\int_{\Gnew}\!\!\fRM(s) {\cdot}
\DT{u}_{\eps}(s) \, \d S \bigg) \d s\,. \label{disc-energy0-lim}
\end{align}
By the very  same lower semicontinuity arguments (also
using~\eqref{e:poinwtiwise-w+} and \eqref{shall-use-later}),  we
also pass to the limit in the discrete total energy
inequality~\eqref{disc-energy}.
\par\noindent
 \textbf{Step $4$: mechanical energy equality.}
First of all, we observe that the following chain rule-type
inequality holds for all $t \in [0,T]$
\begin{align}\nonumber
\Phi & \big(\ude(t),\zde(t)\big)  - \Phi\big(u_{0},z_{0}) +
\mathrm{Var}_{\mathcal{R}}(u,z;[0,t]) \geq \int_0^t
\pairing{}{}{\lambda}{\DT{u}}
 \, \d s
 \\\label{e:faith-delam}
 & \text{for any
$\lambda\in L^2(0,T;W^{2,2}(\Omega{\setminus}\GC;\R^3)^*)$ with
$\lambda(t)\in
 \partial_u \Phi (u(t),z(t))$ for a.a. $t \in (0,T)$,}
\end{align}
where $\partial_u \Phi: W^{2,2}(\Omega{{\setminus}} \GC;\R^3)
\rightrightarrows W^{2,2}(\Omega{{\setminus}} \GC;\R^3)^*$ denotes
the subdifferential w.r.t.\ $u$ of the functional $\Phi:
W^{2,2}(\Omega{{\setminus}} \GC;\R^3)\times L^\infty (\GC) \to \R $
defined in \eqref{8-1-k}. Easy calculations show that the operator
$\partial_u \Phi$ is given by
\begin{align}\nonumber
%\begin{gathered}
&\lambda  \in \partial_u \Phi(u,z)
 \ \ \text{if and only if} \ \ \exists\, \ell \in
 \indabs(u);\ \
%\text{such that}
\forall\, v \in W^{2,2}(\Omega{\setminus} \GC;\R^3):
\\\label{subdif-repre}
&\pairing{}{}{\lambda}{v} =\int_{\Omega}\!
%\left(
\bbC e(u){:}e(v){+}\bbH\nabla e(u){\vdots}\nabla e(v)
%\right)
\, \d x + \int_{\GC}\!\!\!
%\left(\delam z \JUMP{u}{} {\cdot} \JUMP{v}{}{-}
z\alpha_0'(\JUMP{u}{}){\cdot}\JUMP{v}{}
%\right)
\, \d S + \pairing{}{}{\ell}{v},
%\end{gathered}
\end{align}
where, for  notational convenience, we have introduced  the
functional $\abs: W^{{2},2}(\Omega{\setminus}\GC;\R^3) \to
[0,+\infty]$ defined by $ \abs(u)= \ind_{\cone}(\JUMP{u}{})$ (cf.\
\eqref{ind_D}), and its subdifferential $\indabs:
W^{{2},2}(\Omega{\setminus}\GC;\R^3) \rightrightarrows
W^{{2},2}(\Omega{\setminus}\GC;\R^3)^* $.
 In order to prove \eqref{e:faith-delam}  for a fixed selection $\lambda(t)\in
 \partial_u \Phi (u(t),z(t))$,
we exploit a technique, combining Riemann sums and the already
proved semistability condition \eqref{semistab}, which is well-known
in the analysis of rate-independent systems and dates back to
\cite{DMFraToa05}. The main difficulty here is to adapt such a trick
to the case of  a Stieltjes integral (cf.\ \eqref{stieltjes}), and
to do so we will mimick the argument in the proof of \cite[Prop.\
3]{tr-dcds}. For any $\neweta>0$, we take a suitable partition $0 =
t_0^\neweta<t_1^\neweta< \ldots < t_{\Nn}^\neweta=T $ with
$\max_{i=1,\ldots,\Nn} (t_{i}^\neweta{-}t_{i-1}^\neweta) \leq
{1/}\neweta$, in such a way that the functions
$\mathcal{A}_\neweta:[0,T]\to L^\infty(\GC)$ given by
$\mathcal{A}_\neweta(t):= a_1 (\JUMP{u(t_{i-1}^\neweta)}{})$  for $t
\in (t_{i-1}^\neweta, t_i^\neweta]$ fulfill
\begin{align}
\label{a-eta} \mathcal{A}_\neweta \to a_1 (\JUMP{u}{}) \ \ \text{ in
$L^\infty(\SC)\ $ as $\neweta \to \infty$.}
\end{align}
%\end{gathered}\end{equation}
The existence of such partitions follows from the fact that
$u:[0,T]\to W^{2,2}(\Omega{{\setminus}} \GC;\R^3)$ is  continuous,
since $u\in W^{1,2}(0,T;W^{2,2}(\Omega{{\setminus}}\GC;\R^3))$.
%as a mapping $[0,T]\to W^{2-\epsilon,2}(\Omega{{\setminus}} \GC;\R^3)$,
Thus it is also uniformly continuous, and so is the mapping
$\JUMP{u}{}:[0,T]\to L^\infty(\GC)$. Then, we use that  uniformly
continuous mappings admit uniform approximation by piecewise
constant interpolants.
%
%\begin{rcomm} WHICH
%REFERENCE?? IS \eqref{a-eta} TRUE?? IS \eqref{a-eta} SUFFICIENT TO
%PASS TO THE LIMIT in \eqref{stieltjes}?? It seems to me that we
%can't hope for convergence in $C^0 (\overline{\SC})$, since
%$ \mathcal{A}_\neweta$ is NOT continuous w.r.t. time. What shall we do?
%\end{rcomm}
%
In fact, \eqref{a-eta} holds for  all partitions of $[0,T]$ whose
fineness tends to $ 0$, therefore we can choose our partition in
such a way that the semistability \eqref{semistab} holds at all
points $\{ t_i^\neweta \,:\, i =0,\ldots,\Nn{-}1,\
\neweta\in\N\}$. Hence, we write \eqref{semistab} at
$t_{i-1}^\neweta$ tested by $\tilde z=z(t_{i}^\neweta)$, thus
obtaining
\[
\begin{aligned}
\Phi  & \big(u(t_{i-1}^\neweta),z(t_{i-1}^\neweta)\big)
\le\Phi\big(u(t_{i-1}^\neweta),z(t_{i}^\neweta)\big)
+\calD\big(u(t_{i-1}^\neweta),z(t_{i}^\neweta){-}z(t_{i-1}^\neweta)\big)
\\
& = \Phi\big(u(t_{i}^\neweta), z(t_{i}^\neweta)\big)+\int_{\GC}\!\!
a_1 (\JUMP{u(t_{i-1}^\neweta)}{})
|z(t_{i}^\neweta){-}z(t_{i-1}^\neweta)|\, \d S
-\int_{t_{i-1}^\neweta}^{t_{i}^\neweta}\!\!\langle\lambda_n(s),\DT
u(s)\rangle\, \d s
\end{aligned}
\]
for any selection $\lambda_\neweta(t)\in
 \partial_u \Phi (u(t),z(t_{i}^\neweta))$
for a.a. $t\in(t_{i-1}^\neweta, t_i^\neweta]$, $i=1,...,\Nn$, where
we have also used the chain rule  for the convex functional  $u
\mapsto \Phi(u,z(t_{i}^\neweta))$,  cf.\
\cite[Prop.~XI.4.11]{visintin96}. In particular, taking into account
formula \eqref{subdif-repre} for $\partial_u \Phi$, we choose
\begin{equation}
\label{particular-lambda-n}
\begin{aligned}
\lambda_\neweta(t)= \lambda(t) - \rho_n(t) \  & \text{ with
$\rho_n(t)\in W^{{2},2}(\Omega{\setminus}\GC;\R^3)^*$  given by}
\\
& \langle \rho_n(t), v \rangle: =\int_{\GC} (
z(t){-}z(t_{i}^\neweta))\alpha_0'(\JUMP{u(t)}{}) \cdot \JUMP{v}{} \,
\d S.
\end{aligned}
\end{equation}
 Summing for
 $i=1,\ldots,\Nn$, we obtain
 \begin{equation}
\label{e:stie-1}
 \Phi\big(u(T), z(T)\big) - \Phi\big(u_0, z_0\big)
 +\sum_{i=1}^{\Nn}  \int_{t_{i-1}^\neweta}^{t_{i}^\neweta} \int_{\overlineGC}
  \mathcal{A}_\neweta  |\DT
 z|(\d S\d t)\geq\sum_{i=1}^{\Nn}\int_{t_{i-1}^\neweta}^{t_{i}^\neweta}
\langle\lambda_\neweta(s), \DT u(s)\rangle \, \d s.
 \end{equation}
Now, reproducing  the calculations
throughout~\cite[Formulae~(4.70)-(4.74)]{tr1},
 it can be shown that
\begin{equation}\liminf_{\neweta \to \infty} \sum_{i=1}^{\Nn}
\int_{t_{i-1}^\neweta}^{t_{i}^\neweta}\!\langle
\lambda{_\neweta}(s), \DT u(s) \rangle\,\d s  \geq
\int_0^t\!\!\pairing{}{}{\lambda(s)}{\DT{u}(s)}
 \, \d s.
\end{equation}
%\CHECK{Here I am somehow lost - will $\lambda\INSERT{_\neweta}$
%somehow converge or what? The Riemann-sum approximations should
%somewhere be mentioned more in detail..... }
%\begin{rcomm} I think \eqref{particular-lambda-n} should be
%sufficient to conclude.. more details?
%\end{rcomm}
 On the other hand, it
follows from \eqref{a-eta} that
\begin{equation}
\label{stieltjes} \lim_{\neweta \to \infty} \sum_{i=1}^{\Nn}
\int_{t_{i-1}^\neweta}^{t_{i}^\neweta}\int_{{\overlineGC}}\mathcal{A}_\neweta|\DT
 z|(\d S\d t) = \int_{{\overlineSC}}
a_1(\JUMP{u}{})\DT{z}(\d S,d t)
 =\mathrm{Var}_{\mathcal{R}}(u,z;[0,T]).
\end{equation}
Indeed,  $\DT{z}\in C(\overlineSC)^*$ can be extended to
$L^\infty(\SC)^*$ by the Hahn-Banach principle, and then tested by
$\mathcal{A}_\neweta-a_1(\JUMP{u}{})\in L^\infty(\SC)$ which
converges to zero by \eqref{a-eta}. Combining
\eqref{e:stie-1}--\eqref{stieltjes}, we obtain
\eqref{e:faith-delam}.

 %Notice
%that $\indabs = \jum^* \circ \partial I_{\cone} \circ \jum$, with
%$\jum$ the jump operator $\jum(u)=\JUMP{u}{}$, and $\jum^*$ its
%adjoint.

 In order to  make \eqref{e:faith-delam} more explicit, we may
observe that
\begin{equation}
\label{chain-rule}
\begin{aligned}
& \int_{0}^t\!\!\big\langle\ell,\DT{u}\big\rangle\,\d s
=\abs(\ude(t))-\abs(\ude(0))=  \ind_{\cone}\big(
\JUMP{\ude(t)}{}\big) - \ind_{\cone}\big( \JUMP{\ude(0)}{}\big)=0
\\
 &\qquad
\text{for all $\ell\in
L^2(0,T;W^{{2},2}(\Omega{\setminus}\GC;\R^3)^*)$ such that
$\ell(s)\in\partial
%\ind_{\cone}
\abs(\ude(s)) \ \foraa\, s \in (0,T)$,}
\end{aligned}
\end{equation}
by the chain rule for the convex functional  $\abs$ (cf.\ again
\cite[Prop.~XI.4.11]{visintin96}), and by~\eqref{uzero}
and~\eqref{constraints-delam}. Therefore, in view of
\eqref{e:faith-delam}--\eqref{chain-rule}, we conclude the following
inequality for all $t\in[0,T]$
\begin{align}
  & \Phi\big(u(t),z(t)\big)  -
\Phi\big(u_{0},z_{0}) + \mathrm{Var}_{\mathcal{R}}(u,z;[0,t])
\nonumber
\\ &
\ge\int_0^t\!\bigg(\int_{\Omega}
%\left(
\bbC e(\ude){:}e(\DT{u}){+}\bbH \nabla e(u){\vdots}\nabla e(\DT{u})
%\right)
\, \d x
%\d s
+\int_{\GC}\!\!\!
%\left(\delam\zde\JUMP{\ude}{}{\cdot}\JUMP{\DT{u}}{}{-}
z \alpha_0'(\JUMP{u}{}){\cdot}\JUMP{\DT{u}}{}\,\d S
%\right)
\bigg)\d s\,. \label{e:other-crucial-inequality}
\end{align}
In order to develop the test of \eqref{e:weak-momentum-variational}
by $\DT{u}$, we need to distinguish the quasistatic case $\varrho=0$
and the dynamical case $\varrho>0$.
\par\noindent
\textbf{Case $\varrho>0$.} First of all, let us observe that,
under~\eqref{varrho2}, the qualification $\testu \in
W^{1,2}(0,T;L^2(\Omega;\R^d))$ for the test functions
in~\eqref{e:weak-momentum-variational} might be relaxed to
\begin{equation}
\label{relaxed-regu} \testu \in
 L^2(0,T;W_{\cone}^{2,2}(\Omega{{\setminus}} \GC;\R^3)) \cap
W^{1,2}(0,T;W_{\cone}^{2,2}(\Omega{\setminus} \GC;\R^d)^*),
\end{equation}
cf.\ notation \eqref{cone-w}. Indeed, thanks
to~\eqref{constraints-delam} and to the linearity of $\cone(x)$ for
almost all $x \in \GC$, the function $u$
fulfilling~\eqref{e:weak-momentum-variational} is such that $ \DT{u}
\in L^2 (0,T;W_{\cone}^{2,2}(\Omega{\setminus} \GC;\R^d)).$ Note
that \eqref{relaxed-regu} is sufficient to give meaning to the term
$\int_{Q} \DT{u}{\cdot}\DT{\testu} \, \d x \d t$, because
 the spaces $L^2 (0,T;W_{\cone}^{2,2}(\Omega{\setminus}
\GC;\R^d)^*)$ and $L^2 (0,T;W_{\cone}^{2,2}(\Omega{\setminus}
\GC;\R^d))$ are in duality.

Now,  a comparison in~\eqref{e:weak-mom-identity}  yields  that
$\DDT{u}\in L^2 (0,T;W_{\cone}^{2,2}(\Omega{{\setminus}}
\GC;\R^3)^*)$. Therefore, \eqref{additional-cone} ensues, and
 $\DT{u}$ is an admissible test function
for the momentum balance
inclusion~\eqref{e:weak-momentum-variational}, since it fulfills
\eqref{relaxed-regu}. Then, upon proceeding with such a test we
conclude for all $t \in [0,T]$ that
\begin{align}\nonumber
 & \frac{\varrho}2\int_\Omega |\DT{u} (t)|^2 \,\d x  + \int_0^t\!\!
\int_{\Omega}
%\left(
\bbD e(\DT{u}){:} e(\DT{u}){+}\bbG \nabla e(\DT{u}){\vdots} \nabla
e(\DT{u})
%\right)
\, \d x \d s
\\ &\qquad \qquad \nonumber
+  \int_0^t\!\!\int_{\Omega}
%\left(
\bbC e(\ude) {:} e(\DT{u}){+} \bbH \nabla e(\ude) {\vdots} \nabla
e(\DT{u})
%\right)
\, \d x \d s + \int_{0}^t\!\!\int_{\GC}\!\!
%\left( \delam \zde \JUMP{\ude}{}{\cdot}\JUMP{\DT{u}}{}{-}
z \alpha_0'(\JUMP{u}{}){\cdot}\JUMP{\DT{u}}{}
%\right)
\,\d S\d s
\\\label{e:moravia-delam}
&\qquad \qquad =\frac\varrho2 \int_\Omega |\DT u_0|^2\, \d x +
\int_0^t\!\!\left(\,
 \int_\Omega \Theta(\wde)\bbB{:}
e\big(\DT{u}\big)\, \d x  + \int_\Omega\!\FRM{\cdot}\DT{u} \, \d x +
\int_{\Gnew}\!\!\!\fRM{\cdot}\DT{u}\, \d S \right) \, \d s\,.
\end{align}
Combining~\eqref{e:moravia-delam}
with~\eqref{e:other-crucial-inequality}, we get  the converse of
inequality  \eqref{disc-energy0-lim}, hence  the desired mechanical
energy equality \eqref{mech-eq} ensues.
\\
 \textbf{Case
$\varrho=0$:} A comparison in \eqref{e:weak-momentum-variational}
with $\varrho=0$ shows  that the functional
\begin{align} \nonumber
\ell:v\mapsto \int_Q\!\big(\bbD e(\DT{u}){+}\bbC
e(u){-}\bbB\Theta(\w)\big) {:}e(v)
%\,\d x\d t
+
%\int_Q\!
\big(\bbG \nabla e(\DT{u}){+}\bbH \nabla  e(u)\big) {\vdots} \nabla
e(v) \,\d x\d t
 \\  +\int_{\SC}\!\!\!
%\left( \dela z\JUMP{u}{} {-}
z \alpha_0'(\JUMP{u}{})
%\right)
{\cdot}\JUMP{v}{} \,\d S\d t -\int_{Q}\!\FRM{\cdot}v\,\d x\d t-
\int_{\Snew}\!\!\!\fRM{\cdot}v\,\d S\d t
\end{align}
is in $L^2(0,T;W^{2,2}(\Omega{\setminus}\GC;\R^3)^*)$, and fulfills
\begin{align}
\int_0^T I_{\cone}\big(\JUMP{v}{}\big)\,\d t \label{use-of-ell} \ge
\int_{0}^T I_{\cone}\big(\JUMP{u}{}\big)\,\d t + \int_0^T
\big\langle\ell,v{-}u\big\rangle \, \d t.
\end{align}
Hence, $\ell(t)\in\partial{\abs}(\ude(t))$  for almost all $t \in
(0,T)$. Thus, \eqref{chain-rule} yields $\int_0^t
\pairing{}{}{\ell}{\DT{u}}\, \d s =0$ for all $t \in [0,T]$, which
is just relation~\eqref{e:moravia-delam} with $\varrho=0$.  Again,
we combine the latter with~\eqref{e:other-crucial-inequality}, and
conclude the mechanical energy equality \eqref{mech-eq}.
\par\noindent
\textbf{Step $5$: passage to the limit in the enthalpy equation.}
First of all,  we observe  the following chain of inequalities  for
all $t\in[0,T]$:
\begin{align}
\nonumber
 &\mathrm{Var}_{\mathcal{R}}(u_{\eps},z_{\eps};[0,t]) +
\int_0^{t}\int_{\Omega} \bbD e(\DT{u}_{\eps}) {:}
e(\DT{u}_{\eps}){+} \bbG \nabla e(\DT{u}_{\eps}) {\vdots} \nabla
e(\DT{u}_{\eps}) \, \d x \d s
\\&\qquad \nonumber
  \leq \liminf_{\tau \to 0} \int_0^t
\int_{\GC}\zeta_1
\left(\JUMP{\pwc{u}{\tau}}{},\pwl{\DT{z}}{{\eps}\tau} \right)\, \d S
\d s + \int_0^{t}\int_{\Omega} \bbD e(\pwl {\DT{u}}{{\eps}\tau}) {:}
e(\pwl {\DT{u}}{{\eps}\tau}){+} \bbG \nabla e(\pwl
{\DT{u}}{{\eps}\tau}) {\vdots} \nabla e(\pwl {\DT{u}}{{\eps}\tau})\,
\d x \d s
 \\&\qquad \nonumber
\le\limsup_{\tau \to 0} T_\mathrm{kin}(\DT{u}_{0,\tau})
+\Phi_{{\eps}\tau}(u_{0,\tau},z_{0}) -
T_\mathrm{kin}(\pwl{\DT{u}}{{\eps}\tau}(t)) -\Phi_{{\eps}\tau}(\pwc
u{{\eps}\tau}(t), \pwc z{{\eps}\tau}(t))
\\\nonumber
&\qquad\qquad+\int_0^{t}\left(\int_\Omega \Theta(\pwc{\w}{{\eps}\tau})
\bbB{:}e\big(\pwl{\DT{u}}{{\eps}\tau}\big) + \int_{\Omega} \pwc
\FRM{\tau} {\cdot} \pwl{\DT{u}}{{\eps}\tau}\, \d x +\int_{\Gnew}\pwc
\fRM{\tau} {\cdot} \pwl{\DT{u}}{{\eps}\tau} \, \d S \right)\, \d s
\\ \nonumber
&\qquad \leq  T_\mathrm{kin}(\DT{u}_{0}) +\Phi_{\eps}(u_{0},z_{0}) -
T_\mathrm{kin}({\DT{u}}_{\eps}(t)) -\Phi_{\eps}( u_{\eps}(t),
z_{\eps}(t))
\\ \nonumber
&\qquad\qquad+\int_0^{t}\bigg( \int_\Omega\Theta({\w}_{\eps})\bbB{:}
e\big({\DT{u}}_{\eps}\big) +\FRM {\cdot} \DT{u}_{\eps}\,\d x   +
\int_{\Gnew}\!\!\fRM {\cdot} \DT{u}_{\eps} \, \d S\bigg)\, \d s
\nonumber
 \\
\label{e:finally-bis} &\qquad
 =
\mathrm{Diss}_{\mathcal{R}}(u_{\eps},z_{\eps};[0,t]) +
\int_0^{t}\int_{\Omega} \bbD e(\DT{u}_{\eps}) {:}
e(\DT{u}_{\eps}){+} \bbG \nabla e(\DT{u}_{\eps}) {\vdots}\nabla
e(\DT{u}_{\eps})  \, \d x \d s.
\end{align}
Indeed,  the first
%passage
inequality ensues from~\eqref{e:convutau1} and~\eqref{e:finally},
the second one from the discrete mechanical energy
inequality~\eqref{disc-energy0}, the third one from
\eqref{est-init-data}, \eqref{pointiwise-for-u},
\eqref{e:altogether}, \eqref{conv-adiab}, and from
\eqref{data-converg-1}--\eqref{data-converg-2}, cf.\ also Step 3.
Finally, the last equality ensues  from~\eqref{mech-eq} proved in
Step $4$. Thus, all of the above inequalities turn out to hold as
equalities. By a standard
%$\liminf/\limsup$
liminf/limsup argument, we find in particular
\begin{align*}
% \label{e:stop1}
\bbD e(\pwl {\DT{u}}{{\eps}\tau}) {:} e(\pwl {\DT{u}}{{\eps}\tau})
\to \bbD e({\DT{u}}_{\eps}) {:} e({\DT{u}}_{\eps})\ \ \ \text{ and
}\ \ \ \bbG \nabla e(\pwl {\DT{u}}{{\eps}\tau}) {\vdots} \nabla
e(\pwl {\DT{u}}{{\eps}\tau}) \to \bbG \nabla e({\DT{u}}_{\eps})
{\vdots} \nabla e({\DT{u}}_{\eps})\ \ \ \text{ strongly in
$L^1(Q)$.}
\end{align*}
%
%We take the limit of~\eqref{weak-heat-discr}, approximating the test
%functions $\testw \in
%C ([0,T];W^{1,r'}(\Omega{{\setminus}} \GC)) \cap W^{1,r'}(0,T;
%L^{r'}(\Omega))$ with discrete approximations
%$\{ \testw_{\tau}^k \}$ such that, for $\tau\to0$, the related
%interpolants fulfill as $\pwc\testw{\tau} \to \testw  $ in $L^\infty
%(0,T; W^{1,r'}(\Omega{{\setminus}} \GC))$  and $ \pwl \testw{\tau}
%\to \testw $ in $C (0,T; W^{1,r'}(\Omega{{\setminus}} \GC))\cap
%W^{1,r'}(0,T; L^{r'}(\Omega))$.
Combining these convergences with \eqref{conv-adiab} we pass to the
limit in the first term on the right-hand side of
\eqref{weak-heat-discr}.  To take the limit of the second
right-hand-side term,
%one,
we observe that
%
%\[\begin{array}{lll}& \lim_{\tau \to 0} \int_{\SC}\!\! a_1(\JUMP{{\pwc
%u\tau}}{})\pwl{\DT{z}}{{\eps}\tau}\pwc \zeta \tau \, \d S\d t &
%\\ & =  \lim_{\tau \to 0} \int_{\SC}\!\! a_1(\JUMP{{\pwc
%u\tau}}{})\pwl{\DT{z}}{{\eps}\tau} (\pwc \zeta \tau {-} \pwl \zeta
%\tau)  \d S\d t & + \lim_{\tau \to 0} \int_{\SC}\!\!
%(a_1(\JUMP{{\pwc u\tau}}{}){-} a_1(\JUMP{{\pwl
%u\tau}}{}))\pwl{\DT{z}}{{\eps}\tau}  \pwl \zeta \tau \d S\d t
%\\&&+ \lim_{\tau \to 0} \int_{\SC}\!\! a_1(\JUMP{{\pwl
%u\tau}}{})\pwl{\DT{z}}{{\eps}\tau}  \pwl \zeta \tau \d S\d t
%\\ & =: l_1+ l_2+l_3 &\end{array}\]#
\begin{align}\nonumber
\lim_{\tau \to 0} \int_{\SC}\!\! a_1(\JUMP{{\pwc
u\tau}}{})\pwl{\DT{z}}{{\eps}\tau}
v
%\frac{\testw|_{\GC}^+{+}\testw|_{\GC}^-}2
\,\d S\d t
&=\lim_{\tau \to 0} \int_{\SC}\!\!
\big(a_1(\JUMP{{\pwc u\tau}}{}){-} a_1(\JUMP{{\pwl
u\tau}}{})\big)\pwl{\DT{z}}{{\eps}\tau}v\,\d S\d t
\\&
+\lim_{\tau \to 0} \int_{\SC}\!\! a_1(\JUMP{{\pwl
u\tau}}{})\pwl{\DT{z}}{{\eps}\tau}v\,\d S\d t
=0+\int_{\SC}\!\! a_1(\JUMP{{u}}{})v\pwl{\DT{z}}{}(\d S\d t)
\end{align}
for any $v\in C(\overlineSC)$, and in particular for $v=
\frac{\testw|_{\GC}^+{+}\testw|_{\GC}^-}2$; here we used
respectively \eqref{basic-1}, \eqref{measure-convergence}, and
\eqref{a32bis}.
%where we have used short-hand notation $\pwc \zeta{\tau}:
%=\frac{\pwc \testw{\tau}|_{\GC}^+{+} \pwc\testw{\tau}|_{\GC}^-}2 $,
%and  $\pwl\zeta\tau$  ($\zeta$, resp.) for the function defined by
%the same formula with $\pwc \testw{\tau}$ replaced by  $\pwl
%\testw{\tau}$ (by $\testw$, resp.). The
%convergences of  $\pwc {\testw}{{\eps}\tau}$ and $\pwl
%{\testw}{{\eps}\tau}$ to $w$ and trace theorems ensure that $\pwc
%\zeta{\tau} \to \zeta$ in $L^\infty(\SC)$, $\pwl \zeta{\tau} \to
%\zeta$ in $C(\overlineSC)$. Then, $l_1 =0$, also in view of
%estimates \eqref{a32bis} and \eqref{a30}. Thanks to \eqref{basic-1}
%and \eqref{a32bis}  again, we have $l_2=0$. Finally, also exploiting
%\eqref{measure-convergence} we find $l_3 =\int_{\SC}\!\! \zeta\,a_1
%(\JUMP{u}{}) \DT{z}(\d S\d t)$.
Then, we pass to the limit in
the left-hand side of \eqref{weak-heat-discr} by exploiting
\eqref{e:convutau2}, \eqref{e:convwtau1},
\eqref{e:convwtau1-bis-serve},
 \eqref{e:poinwtiwise-w+},
\eqref{e:interesting}, as well as properties \eqref{30c} for
$\mathcal{K}$ and \eqref{eta-affine} for $\eta$, and by arguing in
the very same way as in the proof of \cite[Thm.\ 5.1]{rr+tr}, to
which we refer for all details.

%In
At  the end, employing \eqref{shall-use-later}, we take the limit of
the last term on the right-hand side of \eqref{weak-heat-discr},
thus finding that  the triple $(\ue,\ze,\we)$ fulfils the weak
formulation~\eqref{weak-heat} of the enthalpy equation.
%
%for all \begin{rnew} smooth test functions as in~\eqref{testw}.
%Again by a density argument, we conclude that $(\ue,\ze,\we)$
%fulfill~\eqref{weak-heat} for all $ \testw \in $.
\par\noindent
\textbf{Step $6$:  total energy identity.} We test the weak
formulation~\eqref{weak-heat} of the enthalpy equation  by $1$ and
add it to the mechanical energy equality. This gives the total
energy balance \eqref{total-energy-brittle}. $\hfill\Box$

\end{document}